\documentclass[10pt]{article}
\usepackage{amssymb}
\usepackage{graphicx}
\usepackage{xcolor} 
\usepackage{tensor}
\usepackage{fullpage} 
\usepackage{amsmath}
\usepackage{amsthm}
\usepackage{verbatim}
\usepackage{enumitem}

\providecommand{\MR}{\relax\ifhmode\unskip\space\fi MR }

\providecommand{\href}[2]{#2}
\setlist[enumerate]{leftmargin=1.5em}
\setlist[itemize]{leftmargin=1.5em}

\definecolor{green}{rgb}{0,0.8,0} 

\newtheorem{theorem}{Theorem}[section]

\newtheorem{lemma}[theorem]{Lemma}
\newtheorem{proposition}[theorem]{Proposition}

\theoremstyle{definition}

\theoremstyle{remark}
\newtheorem{remark}[theorem]{Remark}

\numberwithin{equation}{section}
\newcommand{\relphantom}[1]{\mathrel{\phantom{#1}}}

\makeatletter
\newcommand{\nrm}{\@ifstar{\nrmb}{\nrmi}}
\newcommand{\nrmi}[1]{\Vert{#1}\Vert}
\newcommand{\nrmb}[1]{\left\Vert{#1}\right\Vert}
\newcommand{\abs}{\@ifstar{\absb}{\absi}}
\newcommand{\absi}[1]{\vert{#1}\vert}
\newcommand{\absb}[1]{\left\vert{#1}\right\vert}
\newcommand{\brk}{\@ifstar{\brkb}{\brki}}
\newcommand{\brki}[1]{\langle{#1}\rangle}
\newcommand{\brkb}[1]{\left\langle{#1}\right\rangle}
\newcommand{\set}{\@ifstar{\setb}{\seti}}
\newcommand{\seti}[1]{\{#1\}}
\newcommand{\setb}[1]{\left\{ #1\right\}}
\makeatother

\newcommand{\tld}[1]{\widetilde{#1}}
\newcommand{\br}[1]{\overline{#1}}

\newcommand{\nnrm}[1]{{\vert\kern-0.25ex\vert\kern-0.25ex\vert #1 
    \vert\kern-0.25ex\vert\kern-0.25ex\vert}}

\newcommand{\sgn}{{\mathrm{sgn}}\,}
\newcommand{\supp}{{\mathrm{supp}}\,}

\let\Re\relax
\DeclareMathOperator{\Re}{Re}
\let\Im\relax
\DeclareMathOperator{\Im}{Im}

\newcommand{\aeq}{\simeq}
\newcommand{\aleq}{\lesssim}

\newcommand{\lap}{\Delta}

\newcommand{\ud}{\mathrm{d}}
\newcommand{\rd}{\partial}
\newcommand{\nb}{\nabla}

\newcommand{\impmi}{\Leftrightarrow}

\newcommand{\bb}{\Big}

\newcommand{\peq}{\relphantom{=}}			

\newcommand{\alp}{\alpha}
\newcommand{\bt}{\beta}
\newcommand{\gmm}{\gamma}

\newcommand{\dlt}{\delta}

\newcommand{\eps}{\epsilon}

\newcommand{\lmb}{\lambda}

\newcommand{\sgm}{\sigma}

\newcommand{\omg}{\omega}

\newcommand{\bfa}{{\bf a}}

\newcommand{\bfeps}{\boldsymbol{\epsilon}}

\newcommand{\bfPhi}{\boldsymbol{\Phi}}

\newcommand{\bbC}{\mathbb C}

\newcommand{\bbN}{\mathbb N}

\newcommand{\bbR}{\mathbb R}
\newcommand{\bbS}{\mathbb S}
\newcommand{\bbT}{\mathbb T}

\newcommand{\calB}{\mathcal B}
\newcommand{\calC}{\mathcal C}

\newcommand{\calL}{\mathcal L}

\newcommand{\calO}{\mathcal O}

\newcommand{\calR}{\mathcal R}

\newcommand{\err}{\boldsymbol{\epsilon}}		

\newcommand{\tphi}{\tilde{\phi}}					

\begin{document}
 
\title{Illposedness for dispersive equations: \\ Degenerate dispersion and Takeuchi--Mizohata condition}
\author{In-Jee Jeong\thanks{Department of Mathematical Sciences and Research Institute of Mathematics, Seoul National University. \newline  E-mail: injee\_j@snu.ac.kr}\and Sung-Jin Oh\thanks{Department of Mathematics, UC Berkeley and School of Mathematics, Korea Institute for Advanced Study. E-mail: sjoh@math.berkeley.edu}}

\date{\today}

\maketitle


\begin{abstract}
We provide a unified viewpoint on two illposedness mechanisms for dispersive equations in one spatial dimension, namely degenerate dispersion and (the failure of) the Takeuchi--Mizohata condition. Our approach is based on a robust energy- and duality-based method introduced in an earlier work of the authors in the setting of Hall-magnetohydynamics. Concretely, the main results in this paper concern strong illposedness of the Cauchy problem (e.g., non-existence and unboundedness of the solution map) in high-regularity Sobolev spaces for various quasilinear degenerate Schr\"odinger- and KdV-type equations, including the Hunter--Smothers equation, $K(m, n)$ models of Rosenau--Hyman, and the inviscid surface growth model. The mechanism behind these results may be understood in terms of combination of two effects: degenerate dispersion -- which is a property of the principal term in the presence of degenerating coefficients -- and the evolution of the amplitude governed by the Takeuchi--Mizohata condition -- which concerns the subprincipal term. We also demonstrate how the same techniques yield a more quantitative version of the classical $L^{2}$-illposedness result by Mizohata for linear variable-coefficient Schr\"odinger equations with failed Takeuchi--Mizohata condition.
\end{abstract}


\section{Introduction} \label{sec:intro}
\subsection{Quasilinear degenerate dispersive equations}
In this paper, we study the issue of illposedness of the Cauchy problem for various quasilinear dispersive equations in one spatial dimension in the presence of \emph{degenerate dispersion}. We consider both Schr\"odinger and KdV-type equations. Examples of Schr\"odinger-type equations we treat include, for instance, the Hunter--Smothers equation \cite{HS}
\begin{equation} \label{eq:DS} 
\begin{aligned}
&i\rd_t \phi + \rd_x(|\phi|^2\rd_x\phi) = 0,
\end{aligned} 
\end{equation}
which was derived from the Majda--Rosales--Schonbeck equation \cite{MRS} describing the resonant reflection of sound waves off a sawtooth entropy wave, as well as its Hamiltonian variant considered by Germain, Harrop-Griffith and Marzuola \cite{GHGM1, HGM}
\begin{equation} \label{eq:GHGM}
	i \rd_{t} \phi + \overline{\phi} \rd_{x} (\phi \rd_{x} \phi) - \mu_{0} \abs{\phi}^{2} \phi = 0 \qquad (\mu_{0} = -1, 0, 1),
\end{equation}
where $\phi:\bbR_+ \times \bbT\rightarrow \bbC$. In the KdV case, our results cover the $K(m, n)$ equation of Rosenau--Hyman \cite{RH} with $n=2$, i.e.,
\begin{equation} \label{eq:Km2}
	\rd_{t} u + (\tfrac{1}{m} u^{m})_{x} + (\tfrac{1}{2} u^{2})_{xxx} = 0 \quad (m \hbox{ is a nonnegative integer}),
\end{equation}
which has been studied extensively in connection with the remarkable nonlinear phenomenon of the existence of \emph{compactons} (solitons with compact spatial support) \cite{Ros94, Ros05, Ros06, RH, ZiRo17, ZiRo18} (see \cite{RoZi} for a recent review), as well as the inviscid surface growth model (see \cite{BR1,BR2,CY,OR} for the full surface growth model, with the dissipation $- \nu h_{xxxx}$ on the right-hand side)
\begin{equation} \label{eq:sgm-0}
	\rd_t h + ((h_{x})^2)_{xx} = 0,
\end{equation}
where $u, h : \bbR_{+} \times \bbT \to \bbR$. Indeed, degenerate KdV-type equations similar in form to \eqref{eq:Km2} appear in various subjects including sedimentation models (\cite{Zum1,BBK}), shoreline problem in water waves (\cite{LaMe}) and magma dynamics (\cite{SPW,SWR}), to name a few. A more extensive list of references on degenerate KdV equations can be found in \cite{GHGM1,GHGM2}.

In each of these equations, observe that the the highest order term is nonlinear -- more specifically, quadratic or cubic -- in the solution. Vanishing of the solution, therefore, leads to some kind of ``degeneracy'' of the highest order term, which in turn gives rise to delicate issues in the (local) well-posedness of the associated Cauchy problem.

 Indeed, for initial data that are uniformly bounded away from $0$ (a property henceforth referred to as \emph{nondegeneracy}), one expects local well-posedness in high-regularity $L^{2}$-based Sobolev spaces $H^{s}(\bbT)$. For example, in the case of \eqref{eq:DS}, for a sufficiently regular solution $\phi$, one has the conservation of the $L^{2}$-norm:
\begin{equation*}
	\frac{\ud}{\ud t} \left(\int_{\bbT} |\phi|^{2} \, \ud x \right) = 0.
\end{equation*}
Obtaining higher regularity apriori estimates is a much more nontrivial task. One can observe the following bound for $n \geq 1$ (in operator notation) at each $t$: \begin{equation*}
\begin{split}
\frac{\ud}{\ud t} \left( \int_{\bbT} |(\rd_{x} |\phi|^{2} \rd_{x})^n \phi|^2\, \ud x \right) \lesssim_{n} \nrm{\phi}_{H^{2n}}^{4n+2}.
\end{split}
\end{equation*} 
Furthermore, as long as the solution stays nondegenerate at $t$, in the sense that $\inf_{x} |\phi(x, t)|^{2} > c$ for some $c > 0$, a standard argument involving the ellipticity of $(\rd_{x} |\phi(x, t)|^{2} \rd_{x})^n$  allows us to bound $\nrm{\phi(\cdot, t)}_{H^{2n}}^{2}$ by the integral on the left-hand side up to errors of the form $O(\nrm{\phi(\cdot, t)}_{H^{2n}}^{4n+2})$. Putting these together, one can establish a short-time $H^{2n}$ apriori estimate for the solution $\phi$ with nondegenerate initial data. However, in the case of \emph{degenerate} initial data (i.e., those without a uniform bound away from $0$), the above scheme for a short-time $H^{2n}$ apriori estimate with $n \geq 1$ clearly breaks down. 

In this paper, we show that this failure of proof of higher derivative apriori estimates is, in fact, a manifestation of genuine illposedness in standard function spaces. Despite the formal conservation of the $L^{2}$-norm, we demonstrate that all of the equations above are rather strongly ill-posed -- in the sense of non-existence of solutions and unboundedness of the data-to-solution map in suitable set-ups -- in a neighborhood of degenerate initial data (e.g., zero data) in high-regularity spaces ($C^{k-1, 1}$, Sobolev or H\"older spaces).

\subsection{Main results for quasilinear degenerate dispersive equations}
\subsubsection{Results for Schr\"odinger-type equations}
To treat the Hunter--Smothers equation \eqref{eq:DS} and the Hamiltonian equation \eqref{eq:GHGM} simultaneously, we shall consider the general equation
\begin{equation} \label{eq:gDS}
	i \rd_{t} \phi + \abs{\phi}^{2} \rd_{xx} \phi + \alp_{1} \phi \abs{\rd_{x} \phi}^{2} + \bt_{1} \br{\phi} (\rd_{x} \phi)^{2} + \mu_{1} \abs{\phi}^{2} \phi = 0,
\end{equation}
where $\phi : \bbR_{t} \times \bbR_{x} \to \bbC$, $\alp_{1}, \bt_{1} \in \bbR$ and $\mu_{1} \in \bbC$. Indeed, the case $\alp_{1} = \bt_{1} = 1$ and $\mu_{1} = 0$ corresponds to \eqref{eq:DS}, while the case $\alp_{1} = 0$ and $\bt_{1} = 1$ corresponds to \eqref{eq:GHGM}. 

For the statement of the main results, we need to introduce the following exponents. Given $\alp_{1}, \bt_{1} \in \bbR$, we introduce the exponent
	\begin{equation} \label{eq:gDS-w}
		\sgm_{c} = - (\frac{\alp_{1}}{2} + \bt_{1} - 1),
\end{equation}
and let $s_{c}$ be the smallest integer greater than $1$ and $\sgm_{c}-\frac{1}{2}$, i.e., 
\begin{equation} \label{eq:gDS-sc}
s_{c} = \max\set{2, \lfloor \sgm_{c} - \tfrac{1}{2} \rfloor + 1}.
\end{equation}
Note that $\sgm_{c} = - \frac{1}{2}$ and $s_{c} = 2$ for \eqref{eq:DS}, while $\sgm_{c} = 0$ and $s_{c} = 2$ for \eqref{eq:GHGM}. For the significance of $\sgm_{c}$ and $\lfloor \sgm_{c} - \tfrac{1}{2} \rfloor + 1$, see Remark~\ref{rem:gmm-c-s-c} and Section~\ref{subsec:discussion}. We note already that the lower bound $s_{c} \geq 2$ is a technical byproduct of our proof, which we have not attempted to optimize.

Our first result is unboundedness of the solution map (i.e., norm inflation) in $C^{s_{c}}$ near any solution with a linear degeneracy. 

\begin{theorem}[Unboundedness of the solution map near a linearly degenerate solution] \label{thm:ill-posed-unbounded-gDS}
Assume that there exists a solution $f \in L^{\infty}([0,\dlt];C^{s_{c}+1,1}(\bbT))$ to \eqref{eq:gDS} with some $\dlt>0$ such that $f(t=0)=f_{0}$ is linearly degenerate; that is, there exists $x_{0}\in\bbT$ with $f(x_{0}) = 0$ and $f_{0}'(x_{0}) \ne 0$. 
	
	Then, for any $\eps > 0$, $s_{0} \geq s_{c}$, and $0<\dlt'\le\dlt $, we can find $\widetilde{\phi}_0 \in C^\infty(\bbT)$ such that $\nrm{\widetilde{\phi}_0}_{C^{s_{0}}} \le \eps$ and one of the following holds: \begin{itemize}
		\item there exists \emph{no} $L^{\infty}([0,\dlt'];C^{s_{c}}(\bbT))$ solution to \eqref{eq:gDS} with initial data $f_{0} + \widetilde{\phi}_{0}$; or
		\item \emph{any} $L^{\infty}([0,\dlt'];C^{s_{c}}(\bbT))$ solution $\phi$ with $\phi(t=0) = f_{0} + \widetilde{\phi}_{0}$ satisfies \begin{equation*}
			\begin{split}
				\sup_{0< t < \dlt' } \nrm{ \phi(t,\cdot) - f(t,\cdot)}_{C^{s_{c}}(\bbT)} > c_{0} (\dlt')^{-\frac12},  
			\end{split}
		\end{equation*} with some $c_{0}>0$ depending only on $f$. 
	\end{itemize}
\end{theorem}

We remark that the norm inflation assertion immediately implies the inflation of any norm that controls $C^{s_{c}}$, such as $H^{\sgm}$ with $\sgm > s_{c} + \frac{1}{2}$. In fact, our proof readily extends to norm inflation in $H^{\sgm}$ for any $\sgm > \sgm_{c}$ in the second alternative, which is expected to be sharp according to Remark~\ref{rem:gmm-c-s-c} and Section~\ref{subsec:discussion} below; see Remark~\ref{rem:sobolev-instab} for further details. We also remark that the statement of Theorem \ref{thm:ill-posed-unbounded-gDS} should extend over to the case of solutions with orders of degeneracy other than 1. For simplicity, however, we restrict ourselves to the linearly degenerate case, which is ``critical'' in some sense; see Section~\ref{subsec:discussion} below.

Our second result is the nonexistence of a regular local-in-time solution in arbitrarily high-regularity $C^{s_{0}}$ spaces.
\begin{theorem} [Nonexistence of regular local-in-time solution] \label{thm:ill-posed-nonexist-gDS}
	For any $\eps >0$ and $s_{0} \geq s_{c}+2$, there exists an initial data $\phi_0\in C^\infty(\bbT)$ satisfying $\nrm{\phi_0}_{C^{s_{0}}}<\eps$ for which there is no corresponding solution to \eqref{eq:gDS} belonging to $L^{\infty}([0,\dlt];  {C}^{s_{c}+2}(\bbT))$ with any $\dlt>0$. 
\end{theorem}

As an immediate corollary of the above, we have that \eqref{eq:DS} and \eqref{eq:GHGM} are ill-posed in the strongest sense of Hadamard in function spaces which contains $C^\infty$ and controls the $C^{4}$--norm (where $4 = s_{c} + 2$): there exists $C^\infty$ initial data without a local solution in $L^\infty_t  W^{s,p} $ with {$s-\frac{1}{p}>4$} and $L^\infty_t C^{k,\alp} $ with $k+\alp\ge 4$. 

\begin{remark}[Exponents $\sgm_{c}$, $s_{c}$, Takeuchi--Mizohata condition and degenerate dispersion] \label{rem:gmm-c-s-c}
Observe that the illposedness results, Theorems~\ref{thm:ill-posed-unbounded-gDS}--\ref{thm:ill-posed-nonexist-gDS}, hold for all possible coefficients $\alp_{1}$, $\bt_{1}$ in front of subprincipal terms, although these possibly \emph{alter} the exponents $\sgm_{c}$ and $s_{c}$. Heuristically, $\sgm_{c}$ is the expected critical $L^{2}$-Sobolev regularity exponent above which the linearization of \eqref{eq:gDS} around a regular linearly degenerate solution is ill-posed. In fact, the negativity of $\sgm_{c}$ already signals $L^{2}$-illposedness of the linearized equation by the classical Takeuchi--Mizohata condition \cite[Chapter~VII]{Mz2}! Even if $\sgm_{c}$ is positive, it turns out to be $L^{\infty}$-ill-posed after taking $k$ many derivatives with $k > \sgm_{c} - \frac{1}{2}$. This consideration motivates the exponent $s_{c}$ and our illposedness results. We shall elaborate on this remark in Section~\ref{subsec:discussion}.
\end{remark}
\subsubsection{Results for KdV-type equations}
To unify our treatment of KdV-type equations, we consider the general equation
\begin{equation} \label{eq:gdkdv}
	\rd_{t} u + u u_{xxx} + \alp_{1} u_{x} u_{xx} + \frac{\mu_{1}}{m} (u^{m})_{x}= 0,
\end{equation}
where $u : \bbR_{t} \times \bbR_{x} \to \bbR$, $\alp_{1} \in \bbR$, $\mu_{1} \in \bbR$ and $m$ is an integer greater than or equal to $2$. Note that there is no need to separately consider the case $m = 1$, as then this term can be easily removed by the change of variables $(t, x) = (t', x' + \mu_{1} t')$. 

Note that $\alp_{1} = 3$ and $\mu_{1} = 1$ corresponds to the $K(m, 2)$ equation \eqref{eq:Km2}. The inviscid surface growth model \eqref{eq:sgm-0} reduces to the case $\alp_{1} = 3$ and $\mu_{1} = 0$ after making the change of variables $u = \sqrt{2} h_{x}$. 

In the present case, the role of linear degeneracy in the Schr\"odinger case is played by cubic degeneracy, see Section~\ref{subsec:discussion}, below. As before, we introduce the constant
\begin{equation*}
	\sgm_{c} = - (\alp_{1} - \frac{3}{2}),
\end{equation*}
and the regularity exponent
\begin{equation*}
	s_{c} = \max\set{5, \lfloor \sgm_{c}-\tfrac{1}{2} \rfloor+1}.
\end{equation*}
Here, $\sgm_{c}$ is the critical $L^{2}$-Sobolev regularity exponent above which the linearization of \eqref{eq:gdkdv} around a regular cubically degenerate solution is ill-posed (see Section~\ref{subsec:discussion} below). The linearized equation is $L^{\infty}$-ill-posed after taking $k$ many derivatives with $k > \sgm_{c} - \frac{1}{2}$; this motivates the exponent $s_{c}$. The lower bound $s_{c} \geq 5$ is again a nonoptimal technical byproduct of our proof; see Proposition~\ref{prop:C41} for where it is used.

\begin{theorem}[Unboundedness of the solution map near a cubically degenerate solution]\label{thm:ill-posed-unbounded-gdkdv}
	Assume that there exists a solution $f \in L^\infty([0,\dlt]; C^{s-1, 1}(I))$ of \eqref{eq:gdkdv} with some $\dlt>0$ and $I=[a,b]$, such that the initial  data $f_0$ is positive on $I \setminus \rd I$, vanishes cubically on $\rd I$ and $f_{0} \in C^{s_{0}-1, 1}$, where $s_{c} \leq s \leq s_{0}$. Then, for any $\eps > 0$, $s \leq m_{0} \leq s_{0}$, and $0 < \dlt' \leq \dlt$, we can find $\phi_{0} \in C^{\infty}(\bbT; \bbR)$ such that $\supp \phi_{0} \subseteq I \setminus \rd I$, $\nrm{\phi_{0}}_{C^{m_{0}}} \leq \eps$, and one of the following holds:
	\begin{itemize}
\item there exists no solution to \eqref{eq:gdkdv} with initial data $f_{0} + \phi_{0}$ that belongs to $L^{\infty}([0, \dlt']; C^{s-1, 1}(I))$; or
\item any solution $u$ with $u(0) = f_{0} + \phi_{0}$ and belonging to $L^{\infty}([0, \dlt']; C^{s-1, 1}(I))$ satisfies, for every $s_{c} \leq s' \leq 2 \lfloor \frac{s}{2} \rfloor$,
\begin{equation*}
	\sup_{0 < t < \dlt'} \nrm{u(t, \cdot) - f(t, \cdot)}_{C^{s'}(I)} > (\dlt')^{-\frac{1}{2}}.
\end{equation*}
\end{itemize}
\end{theorem}
That $s' \leq 2 \lfloor \frac{s}{2} \rceil$ is not essential and expected be replaceable by $s' \leq s$, but assumed here to simplify the proof (see the proof of Theorem~\ref{thm:ill-posed-unbounded-gdkdv} below). 

As in the Schr\"odinger case, the statement of Theorem \ref{thm:ill-posed-unbounded-gdkdv} should extend over to the case of solutions with orders of degeneracy other than 3, provided that $s_{c}$ is modified suitably. We however focus on the cubic degeneracy case for simplicity, which is ``critical''; see Section~\ref{subsec:discussion} below.

The nonexistence result for \eqref{eq:gdkdv} is as follows. 
\begin{theorem}[Nonexistence of regular local-in-time solution]\label{thm:ill-posed-nonexist-gdkdv}
For any $\eps > 0$ and $s_{c} \leq s \leq s_{0}$, where $s$ is an even integer, there exist an initial data $u_{0} \in C^{\infty}(\bbT)$ satisfying $\nrm{u_{0}}_{C^{s_{0}}} < \eps$ for which there is no corresponding solution to \eqref{eq:gdkdv} belonging to $L^{\infty}([0, \dlt]; C^{s}(\bbT))$ for any $\dlt > 0$.
\end{theorem} 

That $s$ is an even integer is not essential, but assumed here to simplify the proof (see the proof of Theorem~\ref{thm:ill-posed-nonexist-gdkdv} below).

\begin{remark}
	We now give a few simple remarks regarding the above.
	\begin{itemize}
		\item In all of the above, the physical domain could be taken to be $\bbR$ instead of $\bbT$. 
		\item Our proof easily extends to norm inflation in $H^{\sgm}$ for any $\sgm > \sgm_{c}$ in the second alternative in Theorem~\ref{thm:ill-posed-unbounded-gdkdv}. Moreover, in contrast to the Sch\"odinger case, we may also easily extend Theorem~\ref{thm:ill-posed-nonexist-gdkdv} to the non-existence of solutions in $H^{\sgm}$ for any $\sgm > \max \set{\sgm_{c}, 5+\frac{1}{2}}$ (see Remark~\ref{rem:numerologies} for why the situations are different). These numerologies are expected to be sharp, as we shall discuss in Section~\ref{subsec:discussion} below. We refer the reader to Remark~\ref{rem:sobolev-instab} for more details on this modification (which is for norm inflation in the Schr\"odinger case, but the overall idea is the same).
		\item As one can expect, the lower order term $(u^m)_x$, for any $m\ge2$, does not play any essential role in the proof of illposedness of \eqref{eq:gdkdv}. Similarly, one can treat versions of \eqref{eq:DS} which have some lower order terms. 
		\item We expect our results to generalize to $K(m, n)$ with $n > 2$, as well as $\calC(m,a,b)$ equations \cite{Ros06} with $n := a+ b > 2$, by considering degeneracies of order $\frac{3}{n-1}$ (which are critical).
		
	\end{itemize}
\end{remark}

\begin{remark} [Comparison with the work of Ambrose--Simpson--Wright--Yang \cite{ASWY}]
In the pioneering paper \cite{ASWY}, the illposedness of $u_{t} = u u_{xxx}$ in the (fairly low regularity) Sobolev space $H^{2}$ has been proved based on the construction of a compactly supported $H^{2}$ (but not smooth) self-similar solution $A$. However, the existence of such a solution is specific to the equation $u_{t} = u u_{xxx}$, and the proof does not extend to more general class of equations \eqref{eq:gdkdv}, nor to higher regularity Sobolev spaces, as in our results. Our approach is distinct from that of \cite{ASWY}: it does not involve self-similar solutions, but is rather based on appropriate smooth wave packet-type approximate solutions traveling towards the degeneracy; see Section~\ref{subsec:discussion} below. While our results (Theorems~\ref{thm:ill-posed-unbounded-gdkdv} and \ref{thm:ill-posed-nonexist-gdkdv}) do not cover Sobolev regularities as low as $H^{2}$ due to technical reasons, our heuristics suggest that our illposedness mechanism should extend to $H^{\sgm}$ with $\sgm > \sgm_{c} = \frac{3}{2}$. 

Nevertheless, we point out that a key heuristic consideration of our approach, namely, the combined effect of degenerate dispersion and sub-principal term, can already be found in \cite{ASWY}, albeit with a different viewpoint.
\end{remark}
\begin{remark} [Comparison with the works of Germain--Harrop-Griffith--Marzuola \cite{GHGM2} and Harrop-Griffith--Marzuola \cite{HGM}]
For solutions to \eqref{eq:DS} and \eqref{eq:GHGM} with degenerate initial data (i.e., initial data with a zero), our proof identifies and exploits, in a nonlinear fashion, a mechanism by which energies in low frequencies are transferred to high frequencies at arbitrarily fast rates, where the frequencies are defined with respect to the original variable $x$. We emphasize, however, that it does \emph{not} rule out the possibility of well-posedness in regularity classes adapted to the degeneracies of the initial data, by working with a renormalized variable and/or suitable weights. Indeed, such positive results have been proved in the interesting works of Germain--Harrop-Griffith--Marzuola for a KdV-type quasilinear dispersive equation \cite{GHGM2}, and Harrop-Griffith--Marzuola for \eqref{eq:GHGM} \cite{HGM}, where Lagrangian-type coordinates adapted to the solution were used to formulate the function spaces. 
\end{remark}

\subsection{Key mechanism: degenerate dispersion and Takeuchi--Mizohata condition} \label{subsec:discussion}

The nonlinear illposedness results in this paper are firmly based on a detailed and quantitative understanding of illposedness for the linearized equation around a background solution $f$ whose initial data contains a degeneracy. For simplicity, in this subsection we shall assume that the linearization takes the form 
\begin{equation} \label{eq:lin-model}
\begin{cases}
	\rd_{t} u - i \rd_{x} (a \rd_{x} u) - i b \rd_{x} u = (\hbox{lower order})   & \hbox{ in the Schr\"odinger case}, \\
	\rd_{t} u + \tfrac{1}{2}(\rd_{x}^{3} (a u) + a \rd_{x}^{3} u) + b \rd_{x}^{2} u = (\hbox{lower order})   & \hbox{ in the KdV case},
\end{cases}
\end{equation}
where $a = a(x)$ is real-valued in both cases, and $b = b(x)$ is also real-valued in the KdV case.\footnote{For Schr\"odinger-type problems, we regard first-order terms of the form $\tld{b}(x) \rd_{x} \br{u}$ as $\hbox{(lower order)}$, as it can be removed by a suitable change of the dependent variable; see the introduction of $\psi$ in \S\ref{subsec:DS-renrm} below. }

\begin{remark}[On time independence of the coefficients in \eqref{eq:lin-model}] \label{rem:f-static}
While we assumed that $a(x)$ and $b(x)$ are time independent, the actual linearization of \eqref{eq:gDS} and \eqref{eq:gdkdv} on a dynamic background solution $f(t, x)$ would, of course, have time-dependent coefficients. Nevertheless, the time scale of the illposedness mechanism is arbitrarily short, and hence we may effectively approximate these coefficients by the initial values for the purpose of our discussion.
\end{remark}

It is conceptually useful to distinguish two intertwined mechanisms for illposedness, \emph{degenerate dispersion} and \emph{Takeuchi--Mizohata instability}, which can be seen from the principal and sub-principal terms, respectively. Both phenomena must be taken into account to obtain a comprehensive picture of the illposedness of \eqref{eq:gDS} and \eqref{eq:gdkdv} in the presence of a degeneracy in the initial data (and, more concretely, to explain the relevance of the exponents $\sgm_{c}$ and $s_{c}$).

\medskip

\noindent \textit{1. Principal term: dynamics of bicharacteristics}. The illposedness of \eqref{eq:lin-model} from degenerate dispersion can be most easily described at the level of the \emph{bicharacteristic ODE system} associated with the principal symbol $p$ of the spatial part of \eqref{eq:lin-model}, which is given by
\begin{equation}\label{eq:bichar}
		\left\{
		\begin{aligned}
			\dot{X} & = \rd_{\xi} p (X, \Xi), \\
			\dot{\Xi} & = - \rd_{x} p (X, \Xi),
		\end{aligned}
		\right.
	\end{equation} 
where $p(x, \xi) = - a(x) \xi^{2}$ in the Schr\"odinger case, and $- a(x) \xi^{3}$ in the KdV case. By geometric optics, the trajectory $(X(t), \Xi(t))$ describes (at least on sufficiently short timescales) wave packets concentrated near $X(t)$ in the physical space and $\Xi(t)$ in the frequency space; see \eqref{eq:HJ-T-sch} and \eqref{eq:HJ-T-kdv} below for further discussion. If \eqref{eq:bichar} admits the growth of $|\Xi|$ by a definite factor (e.g., $2$) in arbitrarily short timescales, we would have a strong indication of illposedness of \eqref{eq:lin-model} in high-regularity Sobolev spaces. In turn, such a growth may come from some degeneracy of $p$ in $X$ -- this phenomenon is what we shall refer to as \emph{degenerate dispersion}. 
	
To be concrete, let us assume that dynamics is given in $(x, \xi) \in \bbR \times \bbR$ and the coefficient $a$ in $p$ is of the degenerate form $a(x) \approx A x^{n}$ $(n > 0)$ for $\abs{x}$ small, so that \begin{equation} \label{eq:p-model}
	\begin{split}
		p(x,\xi) \approx - A x^n \xi^m \quad \hbox{ for $|x|$ small and $|\xi|$ large}.
	\end{split}
\end{equation} 
Note that $m = 2$ in the Schr\"odinger case and $3$ in the KdV case. The associated bicharacteristic ODE system is \begin{equation}\label{eq:bichar-2}
\left\{
\begin{aligned}
	\dot{ X} & \approx -AmX^n \Xi^{m-1}\\
	\dot{ \Xi} & \approx A n X^{n-1} \Xi^{m},
\end{aligned}
\right.
\end{equation}
In view of the fact that the group velocity $\dot{X}$ vanishes at the point $x = 0$ (since $n > 0$), we shall say that $p$ is \emph{degenerate} at $x =0$. 

We shall now describe the illposedness mechanism of degenerate dispersion in this concrete case. (This analysis can be found in the introduction of \cite{GHGM2} as well.) With a change of the time variable, we may take $A = 1$. Assume, for the sake of this heuristic discussion, that $\approx$'s above are exact equalities. Consider the solution to \eqref{eq:bichar-2} with initial conditions $(X_{0}, \Xi_{0})$, where $0<X_0\ll 1$ and $\Xi_0\gg1$. Then, appealing to the fact that $X^{n}\Xi^{m}$ is  conserved in time, we have
\begin{equation*}
\begin{split}
	\Xi(t) = \Xi_{0}\left( 1 + (n-m) X_0^{n-1}\Xi_0^{m-1}t \right)^{\frac{m}{m-n}}
\end{split}
\end{equation*} for $m\ne n$. When $m=n$, which will be referred to as the \emph{critical} case, we have instead \begin{equation} \label{eq:bichar-Xi-exact}
\begin{split}
	\Xi(t) = \Xi_{0} \exp\left(  m X_{0}^{m-1} \Xi_{0}^{m-1} t \right). 
\end{split}
\end{equation}
In all cases, the frequency magnitude doubles (i.e., $\abs{\Xi(\tau_{2})} = 2 \abs{\Xi_{0}}$) at time $\tau_{2} \aeq \abs{\Xi_{0}}^{1-m} \abs{X_{0}}^{1-n}$. If the order of $p$ is greater than $1$ (i.e., $m > 1$), the doubling time $\tau_{2}$ may be taken to be arbitrarily small by choosing $\abs{\Xi_{0}}$ large, as we desired. Such an arbitrarily fast growth of $\Xi$ suggests that high derivatives of the solution following this bicharacteristic flow would also grow arbitrarily fast -- this is what we shall refer to as \emph{illposedness via degenerate dispersion}.

Finally, let us connect the above model case to the equations considered in this work. Recall that the principal coefficient $a$ in the linearized operator is determined by the background solution $f$, where $a = \abs{f}^{2}$ for \eqref{eq:gDS} and $a = f$ for \eqref{eq:gdkdv}. Since the relevant frequency doubling time scale $\tau_{2}$ is arbitrarily small, it is reasonable to make the approximation $f \approx f_{0}$. Assuming that $f_{0}$ is degenerate at $x = 0$, in the sense that $\abs{f_{0}}$ vanishes to some finite order at $0$, we arrive at the ansatz $a(x) \approx A x^{n}$ for some $A \neq 0$ and $n > 0$. 

\begin{remark}[Critical degeneracy] \label{rem:crit-degen}
In this work, for simplicity, we shall consider only background solutions with critical degeneracy $n = m$. The heuristics suggest, however, that a similar arbitrary fast growth of $\abs{\Xi}$ is expected for any order of degeneracy $n > 0$. The techniques in this paper should be generalizable to these cases.
\end{remark}

\medskip

\noindent \textit{2. Sub-principal term: evolution of wave packet amplitude and Takeuchi--Mizohata condition.}
While sub-principal terms do not enter in the dynamics of bicharacteristics, they need to be considered in order to fully understand the well/illposedness issues for \eqref{eq:lin-model}. In fact, the sub-principal term may already cause illposedness in $L^{2}$ even when the principal term is \emph{nondegenerate}! This phenomenon is captured by the classical \emph{Takeuchi--Mizohata condition} (after the works \cite{Ta,Mz} in the Schr\"odinger case); see \eqref{eq:Mz-condition} below. 

To understand this phenomenon, it is instructive to delve a little deeper into the construction of wave packet (approximate) solutions for \eqref{eq:lin-model}. Consider the ansatz\footnote{In the KdV case, we can take the real or imaginary part of $\bfa e^{i \bfPhi}$ at the end to obtain a real-valued wave packet.} $u = \bfa(t, x) e^{i \bfPhi(t, x)}$ with the following properties: (i)~$\bfPhi(t, x)$ is real-valued, $\rd_{x} \bfPhi(0, x) = \Xi_{0}$ on the support of $\bfa(0, x)$, and (ii)~$\bfa(t, x)$ is complex-valued, and $\bfa(0, x)$ is a smooth bump function adapted to a small ball centered at $X_{0}$. With the expectation that the $\rd_{x} \bfPhi(t, x)$ stays large compared to the characteristic frequencies of $\bfa$, $a$ and $b$, we may write
\begin{align*}
	e^{-i \bfPhi}(\rd_{t} - i \rd_{x} a(x) \rd_{x} - i b \rd_{x} ) (\bfa e^{i \bfPhi})
	&= i (\rd_{t} \bfPhi + a (\rd_{x} \bfPhi)^{2}) \bfa  \\
	&\peq + \rd_{t} \bfa + 2 a \rd_{x} \bfPhi \rd_{x} \bfa  + \left( \rd_{x} a + b + a \tfrac{\rd_{x}^{2} \bfPhi}{\rd_{x} \bfPhi} \right) \rd_{x} \bfPhi \bfa + \cdots,
\end{align*}
in the Schr\"odinger case (where we omitted terms that does not involve $\rd_{x} \bfPhi$), and
\begin{align*}
e^{-i \bfPhi}(\rd_{t} + \tfrac{1}{2} (\rd_{x}^{3} a + a \rd_{x}^{3}) + b \rd_{x}^{2}) (\bfa e^{i \bfPhi}) 
&= i (\rd_{t} \bfPhi - a (\rd_{x} \bfPhi)^{3}) \bfa \\
&\peq + \rd_{t} \bfa - 3 a (\rd_{x} \bfPhi)^{2} \rd_{x} \bfa - (b + \tfrac{3}{2} \rd_{x} a - 3 a \tfrac{\rd_{x}^{2} \bfPhi}{\rd_{x} \bfPhi}) (\rd_{x} \bfPhi)^{2} \bfa + \cdots
\end{align*}
in the KdV case (where we omitted terms of order $0$ and $1$ in $\rd_{x} \bfPhi$). To eliminate the main terms on the right-hand sides, we are led to impose the following classical \emph{Hamilton--Jacobi} and \emph{transport equations} for $\bfPhi$ and $\bfa$:
\begin{align}
&\left\{
\begin{array}{l}
	 \rd_{t} \bfPhi + a (\rd_{x} \bfPhi)^{2} = 0 \\
	 \rd_{t} \bfa + 2 a \rd_{x} \bfPhi \rd_{x} \bfa  + \rd_{x} (a \rd_{x} \bfPhi) \bfa = - b \rd_{x} \bfPhi \bfa
\end{array}
\right. & & \hbox{ in the Schr\"odinger case} \label{eq:HJ-T-sch}, \\
&\left\{
\begin{array}{l}
	 \rd_{t} \bfPhi - a (\rd_{x} \bfPhi)^{3} = 0 \\
	 \rd_{t} \bfa - 3 a (\rd_{x} \bfPhi)^{2} \rd_{x} \bfa - \frac{3}{2} \rd_{x} (a (\rd_{x} \bfPhi)^{2}) \bfa = b (\rd_{x} \bfPhi)^{2} \bfa
\end{array}
\right. & & \hbox{ in the KdV case} \label{eq:HJ-T-kdv}.
\end{align}
Observe that $(X(t), \Xi(t))$ solving \eqref{eq:bichar} are precisely the bicharacteristics for the above equations in the method of characteristics \cite[Chapter~3]{Ev}, which explains the relevance of \eqref{eq:bichar}. Moreover, the transport equations show how $b$ influences the evolution of the amplitude $\bfa$. In fact, we may easily check that
\begin{equation} \label{eq:T-L2}
	\frac{1}{2} \frac{\ud}{\ud t} \nrm{\bfa}_{L^{2}}^{2}
	= \begin{cases}
	- \brk{\Re b \, \rd_{x} \bfPhi \, \bfa, \bfa} & \hbox{in the Schr\"odinger case}, \\
	\brk{b (\rd_{x} \bfPhi)^{2} \bfa, \bfa} & \hbox{in the KdV case},
	\end{cases}
\end{equation}
which clearly demonstrates how $b$ influence the evolution of the $L^{2}$-norm (here $\brk{\cdot, \cdot}$ is the standard $L^{2}$-inner product).

We are now ready to give a heuristic derivation of the Takeuchi--Mizohata conditions. By the method of characteristics, we expect, at least for a short time, that $\rd_{x} \bfPhi(t, X(t)) = \Xi(t)$ and $\bfa$ remains a bump function adapted to a ball centered at $X(t)$. Hence, on $\supp \bfa$, we expect
\begin{align*}
	- \brk{\Re b \, \rd_{x} \bfPhi \, \bfa, \bfa} \approx - \Re b(X(t)) \Xi(t) \nrm{\bfa}_{L^{2}}^{2} = - \frac{\Re b(X(t))}{2 a(X(t))} \dot{X}(t) \nrm{\bfa}_{L^{2}}^{2},
\end{align*}
in the Schr\"odinger case, and
\begin{align*}
	\brk{b (\rd_{x} \bfPhi)^{2} \bfa, \bfa} \approx b(X(t)) \Xi(t)^{2} \nrm{\bfa}_{L^{2}}^{2} = - \frac{b(X(t))}{3 a(X(t))} \dot{X}(t) \nrm{\bfa}_{L^{2}}^{2},
\end{align*}
in the KdV case, where we used \eqref{eq:bichar} for the last equalities. Using \eqref{eq:T-L2} and $\nrm{u}_{L^{2}} = \nrm{\bfa}_{L^{2}}$, we arrive at the expectations
\begin{equation} \label{eq:Mz-pre}
	\nrm{u(t)}_{L^{2}} \aeq \begin{cases}
	\exp\left( - \int_{X_{0}}^{X(t)} \frac{\Re b}{2 a} \, \ud x \right) \nrm{u(t=0)}_{L^{2}} & \hbox{ in the Schr\"odinger case}, \\
	\exp\left( \int_{X(t)}^{X_{0}} \frac{b}{3 a} \, \ud x \right) \nrm{u(t=0)}_{L^{2}} & \hbox{ in the KdV case}.
	\end{cases}
\end{equation}
The \emph{Takeuchi--Mizohata conditions} (cf.~\cite{Ta, Mz} in the Schr\"odinger case; see \cite{Ak, AmWr} and Remark~\ref{rem:Mz-KdV} for the KdV case) are simply sufficient conditions for the forward-in-time boundedness of $\nrm{\bfa}_{L^{2}}$ read off from \eqref{eq:Mz-pre}: 
\begin{align} 
\sup_{ x_{0} < x_{1}} \abs*{\int_{x_{0}}^{x_{1}} \frac{\Re b}{2 a} \, \ud x} < + \infty & \quad \hbox{ in the Schr\"odinger case}, \label{eq:Mz-condition} \\
\sup_{x_{0} < x_{1}} \int_{x_{0}}^{x_{1}} \frac{b}{3 a} \, \ud x < + \infty  & \quad \hbox{ in the KdV case}. \label{eq:Mz-condition-kdv}
\end{align}
Conversely, the failure of the Takeuchi--Mizohata conditions \eqref{eq:Mz-condition}, \eqref{eq:Mz-condition-kdv} signals arbitrarily fast growth (i.e., norm inflation) of the $L^{2}$-norm of $u$, since $X(t)$ may travel arbitrarily far from $X_{0}$ in any fixed duration of time if $\Xi_{0}$ is large. In this paper, we shall refer to this norm inflation (or illposedness) mechanism as the \emph{Takeuchi--Mizohata instability}. Below, we shall consider the interaction of degenerate dispersion and the Takeuchi--Mizohata instability, which provides us with a detailed heuristic understanding of the illposedness properties of the linearization of \eqref{eq:gDS} and \eqref{eq:gdkdv} in the presence of a (critical) degeneracy in the initial data.

\begin{remark}[Rigorous results on Takeuchi--Mizohata-type conditions] \label{rem:Mz-sufficiency}
The necessity of \eqref{eq:Mz-condition} for the $L^{2}$-wellposedness of \eqref{eq:lin-model} in the Schr\"odinger case has been known since in the early works \cite{Ta, Mz}; see also \cite{Ak} for the KdV case. On the other hand, whether such a condition alone is sufficient for $L^{2}$ boundedness in general is less clear in general, especially in higher dimensions. Nevertheless, some strengthened form of the Takeuchi--Mizohata condition underlies many works on the wellposedness of the Cauchy problem for linear and even nonlinear Schr\"odinger- and KdV-type equations; see, e.g., \cite{Ak, AkAmWr, AmWr, HG1, HG2, KPV98, KPV04, MMT1, MMT2, MMT3, Mz2}.
\end{remark}
\begin{remark}[Role of $\sgn b$ in the KdV case] \label{rem:Mz-KdV} 
Observe that the absolute value is needed in the Schr\"odinger case \eqref{eq:Mz-condition} since $X(t)$ may travel in both directions, while it is not necessary in the KdV case \eqref{eq:Mz-condition-kdv} since $X(t)$ is \emph{always} decreasing if $a$ is positive (resp.~increasing if $a$ is negative) according to \eqref{eq:bichar}. In particular, in the KdV case, \eqref{eq:Mz-condition-kdv} is always satisfied if $b < 0$, and even when $b$ has some positive parts, it is possible that the Takeuchi--Mizohata condition is still satisfied (e.g., when $b$ oscillates). This phenomenon has been explored by Ambrose--Wright \cite{AmWr}, who prove wellposedness of some variable coefficient linear KdV-type equations in the periodic setting in the presence of the positive part of $b$ (referred to as ``anti-diffusion'' in that paper).
\end{remark}
\begin{remark}
While our main focus is the interaction of Takeuchi--Mizohota instability with degenerate dispersion, the method developed in this paper also provides a new and effective way to rigorously establish the necessity of \eqref{eq:Mz-condition} and \eqref{eq:Mz-condition-kdv} for the $L^{2}$-wellposedness of \eqref{eq:lin-model}. We refer the reader to Section~\ref{subsec:Mz-results} and Appendix~\ref{sec:Mz} below for sample results in the Schr\"odinger case for $a = 1$ (but in arbitrarily dimensions).
\end{remark}

\medskip
\noindent \textit{3.~Combined effect of degenerate dispersion and Takeuchi--Mizohata instability.} We are now ready to discuss the combined effect of the principal and sub-principal terms in \eqref{eq:lin-model} obtained by linearizing around a background solution $f$ with a degeneracy. Keeping Remarks~\ref{rem:f-static} and \ref{rem:crit-degen} in mind, we consider the linearization of \eqref{eq:gDS} and \eqref{eq:gdkdv} around $f(t, x) = x$ and $x^{3}$ (for $\abs{x}$ small), respectively. Then we arrive at \eqref{eq:lin-model} with
\begin{align*}
	& a(x) = x^{2}, \quad b(x) = 3(\tfrac{\alp_{1}}{2} + \bt_{1} - 1) x & & \hbox{ in the Schr\"odinger case}, \\
	& a(x) = x^{3}, \quad b(x) = 3 (\alp_{1} - \tfrac{3}{2}) x^{2} & & \hbox{ in the KdV case}.
\end{align*}
Recall from the above that we are considering bicharacteristics $(X(t), \Xi(t))$ with $X_{0} > 0$ and $X(t)$ traveling to the degeneracy $0$ in both cases. Wave packets corresponding to such bicharacteristics shall be called \emph{degenerating wave packets}. 

The relevant Takeuchi--Mizohata condition (see \eqref{eq:Mz-condition-degen} below with $\sgm =0$) for $\nrm{u}_{L^{2}}$ may or may not hold, meaning that degenerating wave packets may or may not remain bounded in $L^{2}$. Nevertheless, it \emph{always fails for high enough derivatives}, which is consistent with the heuristic $\Xi(t) \to \infty$! Indeed, observe that commutation of \eqref{eq:lin-model} with $\rd_{x}^{\sgm}$ leads to a similar equation for $\rd_{x}^{\sgm} u$ but with the following coefficients:  
\begin{align*}
	& a(x) = x^{2}, \quad b(x) = 3(\tfrac{\alp_{1}}{2} + \bt_{1} - 1 + \sgm) x & & \hbox{ in the Schr\"odinger case}, \\
	& a(x) = x^{3}, \quad b(x) = 3 (\alp_{1} - \tfrac{3}{2} + \sgm) x^{2} & & \hbox{ in the KdV case}.
\end{align*}
In view of $0 < X(t) < X_{0}$, the Takeuchi--Mizohata condition for boundedness of $\nrm{u}_{H^{\sgm}}$ is:
\begin{equation} \label{eq:Mz-condition-degen}
\begin{aligned}
	&\sup_{0 < x_{0} < x_{1} \ll 1} \int_{x_{0}}^{x_{1}} (\tfrac{\alp_{1}}{2} + \bt_{1} - 1 + \sgm) \frac{\ud x}{x} < + \infty & & \hbox{ in the Schr\"odinger case}, \\
	&\sup_{0 < x_{0} < x_{1} \ll 1} \int_{x_{0}}^{x_{1}} (\alp_{1} - \tfrac{3}{2} + \sgm) \frac{\ud x}{x} < + \infty & & \hbox{ in the KdV case},
\end{aligned}
\end{equation}
which fails exactly when $\sgm > \sgm_{c}$ in both cases. Moreover, the preceding heuristic analysis suggests that the $H^{\sgm}$ norm of the degenerating wave packet grows if $\sgm > \sgm_{c}$, stays constant if $\sgm = \sgm_{c}$ and decays if $\sgm < \sgm_{c}$. This consideration explains the relevance of the exponent $\sgm_{c}$.

Working directly with the transport equations for $\bfa$ in \eqref{eq:HJ-T-sch}--\eqref{eq:HJ-T-kdv} in place of \eqref{eq:T-L2}, we may also see that for $\tld{s}_{c} = \sgm_{c} + \frac{1}{2}$, the $W^{s, \infty}$ norm of the wave packet grows if $s > \tld{s}_{c}$, stays constant if $s = \tld{s}_{c}$ and decays if $s < \tld{s}_{c}$. This consideration motivates the integer exponent $s_{c}$ in our results.

\subsection{Discussion of the proof} \label{subsec:discussion-pf}
Our discussion so far has been rather formal; deriving actual nonlinear illposedness in standard function spaces requires more ideas. Our main technical contribution in this work is developing a robust scheme for establishing quantitative illposedness, which is not only able to deduce strong illposedness in quasilinear cases but also yields much stronger statements for linear equations. The scheme largely consists of three parts: 1.~construction of degenerating wave packets for the linearized equation, 2.~duality testing argument and 3.~incorporation of the nonlinearity.

\medskip

\noindent \textit{1. Degenerating wave packets}. 
We first describe the ideas for construction of a degenerating wave packet. Compared to the heuristic discussion above, the actual construction of such an approximate solution $\tld{u}$ to \eqref{eq:lin-model} has to (1)~allow for time-dependent coefficients $a = a(t, x)$ and $b = b(t, x)$ (as the background solution may depend on time, see Part~3 below), and (2)~solve the linearized equation up to an equation error $\bfeps_{\tld{u}}$ of size $\calO(\Xi_{0}^{m-1-\dlt})$ (in a suitable norm) for some $\dlt > 0$ (here, $m=2$ for Schr\"odinger and $3$ for KdV). Property~(2) is necessary to justify the approximation on a longer time scale than $\Xi_{0}^{1-m}$, which is the instability time scale; see \eqref{eq:bichar-Xi-exact}. 

Our idea is to make appropriate changes of the independent and dependent variables from $(x, u)$ to $(y, v)$ to reduce the problem to the constant coefficient case, for which the construction is standard. For time-dependent coefficients $a = a(t, x)$ and $b = b(t, x)$, the transformation $(x, u) \mapsto (y, \check{v})$ is of the form
\begin{equation*}
	\ud x =  (a(t, x))^{\frac{1}{m}} \ud y, \quad u = (w \check{w})^{-1} \check{v}, \quad \hbox{ where } w^{-1} \rd_{x} w = \frac{\Re b}{m a}, \, \check{w}^{-1} \rd_{x} \check{w} = \frac{\rd_{x} a}{2m a}. 
\end{equation*}
Roughly speaking, Takeuchi--Mizohata instability is renormalized by the conjugation of the dependent variable by the weight $w$, in the sense that $v(t, x) := w u(t, x)$ solves \eqref{eq:lin-model} with $b = 0$ (with possibly different lower order terms). Similarly, degenerate dispersion is renormalized by the change of variables $x \mapsto y$ accompanied with the conjugation of the dependent variable by the weight $\check{w}$, in the sense that the $\check{v}(t, y)$ solves the constant coefficient problem $(\rd_{t} + i (i \rd_{y})^{m}) \check{v} = (\hbox{lower order terms})$. Now, starting from a standard wave packet for the constant coefficient problem traveling towards the degeneracy and returning to original variables, we obtain a degenerating wave packet $\tld{u}$ with the desired properties.

In order to make the above heuristic discussion precise, there are several more factors to consider. For instance, we need to make sure that the contribution of $\rd_{t} a$, $\rd_{t} b$ are indeed acceptable, which ultimately relies on the estimates we have on the time derivative of the background solution $f(t, x)$ in applications; see 3.~below. In the Schr\"odinger case, we need the following two additional ideas: (i)~an extra change of dependent variables to treat terms of the form $\tld{b} \rd_{x} \br{u}$, and (ii)~an extra phase rotation $e^{i \lmb S}$ for the wave packet $\check{v}$ to treat terms of the form $-(\Im b + \rd_{t} y) \rd_{y} \check{v}$, both of which are potentially problematic for achieving Property~(2). For more details, see the proofs of Propositions~\ref{prop:DS-WP} and \ref{prop:KdV-deg-wp} for details.

\begin{remark}[Numerologies in the Schr\"odinger vs.~KdV cases] \label{rem:numerologies}
In order to justify the properties of $\tld{u}$ needed for the proof of the $H^{s}$- or $C^{s}$-norm growth, (ii) above forces the technical restriction that $f \in C^{s+1, 1}$ with $s \geq 2$ in the Schr\"odinger case, while $f \in C^{s-1, 1}$ with $s \geq 4$ is sufficient in the KdV case; compare the degeneration bounds in Propositions~\ref{prop:DS-WP} and \ref{prop:KdV-deg-wp}. This point explains the different numerologies in Theorems~\ref{thm:ill-posed-unbounded-gDS} and \ref{thm:ill-posed-nonexist-gDS} in the Schr\"odinger case versus Theorems~\ref{thm:ill-posed-unbounded-gdkdv} and \ref{thm:ill-posed-nonexist-gdkdv} in the KdV case.
\end{remark}

\medskip

\noindent \textit{2. Modified energy estimate and duality testing argument}. 
In order to upgrade the norm growth for a degenerating wave packet $\tld{u}$ to an actual solution $u$ to \eqref{eq:lin-model}, we adapt the \emph{energy estimate and duality method} introduced in our previous work \cite{JO1} on Hall--MHD. Here we shall briefly explain the argument, in the simplest setup. 

Given a degenerating wave packet $\tld{u}$ for \eqref{eq:lin-model}, denote by $u$ any\footnote{Note that it is a-priori possible that uniqueness of the Cauchy problem for \eqref{eq:lin-model} fails. Nevertheless, the method is still applicable and establishes the norm growth of every solution $u$ with the same initial data satisfying \eqref{eq:mod-en-u}.} solution to \eqref{eq:lin-model} with $u(t=0) = \tld{u}_{0}$, where $\tld{u}_{0}:=\tld{u}(t=0)$. In view of the aforementioned fact that $v = w u$ solves \eqref{eq:lin-model} with $b = 0$, the following \emph{modified energy estimates} should hold (at least when $u$ is sufficiently regular):
\begin{equation} \label{eq:mod-en-u}
	\begin{split}
		\nrm{w u}_{L^\infty_t ([0, t_{0}]; L^2)} \lesssim \nrm{w u_0}_{ L^2},\qquad \nrm{w \widetilde{u}}_{L^\infty_t ([0, t_{0}]; L^2)} \lesssim \nrm{w \widetilde{u}_0}_{ L^2},
	\end{split}
\end{equation}
where $0 < t_{0} < 1$. By the same token, the following \emph{generalized (bilinear) energy estimate} should also hold:
\begin{align*}
	\left| \frac{\ud}{\ud t} \brk{ w u, w \widetilde{u} }  \right| \lesssim \nrm{w \err_u}_{L^2} \nrm{w \widetilde{u}}_{L^2} + \nrm{w u}_{L^2} \nrm{w \err_{\widetilde{u}}}_{L^2}. 
\end{align*}
(We remark that $\rd_{t} w$ also arises, but in applications we shall have $\abs{\rd_{t} w} \aleq w$.) Here, $\err_u$ and $\err_{\widetilde{u}}$ are the errors associated with $u$ and $\widetilde{u}$ viewed as approximate solutions to \eqref{eq:lin-model}. Then, as long as the error terms are bounded, we may deduce that $\brk{ w u, w \tld{u} }  \aeq \brk{ w u_0, w \tld{u}_0 } = \nrm{w \tld{u}_{0}}_{L^{2}}^{2}$, which allows us to obtain behavior of $u$ in various norms by simply estimating the degenerating wave packet $\tld{u}$ and using duality.

In actual applications, the errors often contain derivatives and hence $\nrm{w \bfeps_{\tld{u}}}_{L^{2}}$ (resp.~$\nrm{w \bfeps_{u}}_{L^{2}}$) may diverge as $|\Xi_0|\to\infty$. Nevertheless, in view of the fact the instability time-scale is $\aeq \abs{\Xi_{0}}^{1-m}$, for the above argument to work it suffices to have $\int_{0}^{t_{0}} \nrm{w \err_{\tld{u}}}_{L^2} \aleq 1$ (resp.~$\int_{0}^{t_{0}} \nrm{w \err_{u}}_{L^2} \aleq 1$) for $t_{0} > \abs{\Xi_{0}}^{1-m+\dlt}$. For $\tld{u}$, this is precisely Property~(2) in the preceding discussion. For an actual solution $u$ to \eqref{eq:lin-model}, this follows from the fact that $\err_{u}$ does not contain principal nor sub-principal terms (except $\tld{b} \rd_{x} \br{u}$ in the Schr\"odinger case, which may be integrated by parts away).

\medskip

\noindent \textit{3. Incorporation of the nonlinearity}.
The ideas discussed so far explain how to prove the illposedness of \eqref{eq:lin-model} that arises from linearizing \eqref{eq:gDS} and \eqref{eq:gdkdv} around a regular solution $f(t, x)$ whose initial data has a critical degeneracy at $x = 0$. As in \cite{JO1}, the nonlinear norminflation results (Theorems~\ref{thm:ill-posed-unbounded-gDS} and \ref{thm:ill-posed-unbounded-gdkdv}) is derived by assuming the existence of a nonlinear perturbation $u$ around $f$ (i.e., $f+u$ solves the nonlinear equation) without the instability behavior, then applying the above argument. Moreover, the nonlinear non-existence results (Theorems~\ref{thm:ill-posed-nonexist-gDS} and \ref{thm:ill-posed-nonexist-gdkdv}) are proved by superposition of infinitely many configurations exhibiting norm inflation (with unbounded rates of growth), with disjoint supports in physical space. We refer to \cite[Section~1.6]{JO1} for a more detailed summary of the ideas involved, and to Sections~\ref{subsec:DS-unbounded}, \ref{subsec:DS-nonexist}, \ref{subsec:unbounded-gdkdv} and \ref{subsec:Kmn-nonexist} for details. A key new feature of the present paper compared to \cite{JO1}, however, is that the background solution $f$ need not be stationary solutions, and are given as a part of the contradiction hypothesis in the proof of the non-existence theorems.

\subsection{Revisiting $L^{2}$-illposedness \`a la Takeuchi--Mizohata} \label{subsec:Mz-results}
In view of the extensive appearance of the Takeuchi--Mizohata instability in this paper, it is perhaps not surprising that our techniques also apply to the original setting considered by Takeuchi and Mizohata of $L^{2}$-illposedness of linear nondegenerate Schr\"odinger-type equations. In Appendix~\ref{sec:Mz}, we provide a few results concerning the Takeuchi--Mizohata condition obtained through our approach. In particular, we recover the following result of Mizohata \cite[\S VII.2]{Mz2}:
\begin{proposition} \label{prop:Mz-d-norm-infl}
Consider the following linear first-order perturbation of the Schr\"odinger equation on $\bbR^{d}$:
\begin{equation} \label{eq:Mz-d-0}
	i \rd_{t} u + \lap u + b^{j}(x) \rd_{j} u = 0
\end{equation}
where $b \in C^{1, 1}(\bbR^{d})$. Suppose that the Takeuchi--Mizohata condition for \eqref{eq:Mz-d-0} fails, i.e.,
\begin{equation} \label{eq:Mz-d-cond}
\sup \set*{\int_{0}^{T} \Re b^{j}(x - 2 s \omg) \omg_{j} \, \ud s : x \in \bbR^{d}, \, \omg \in \bbS^{d-1}, \, T > 0 } = + \infty.
\end{equation}
Then for any $\dlt > 0$, every solution map $L^{2} \to L^{\infty}_{t} ([0, \dlt]; L^{2})$ for \eqref{eq:Mz-d-0}, if it exists, is unbounded.
\end{proposition}

Note that Proposition~\ref{prop:Mz-d-norm-infl} clearly implies the result proved in \cite[\S VII.2]{Mz2}, namely, the impossibility of having a solution map for the inhomogeneous equation
\begin{equation} \label{eq:Mz-d-f}
	i \rd_{t} u + \lap u + b^{j}(x) \rd_{j} u = f
\end{equation}
satisfying
\begin{equation*}
	\nrm{u}_{L^{\infty}([0, \dlt]; L^{2})} \leq C_{0} \left(\nrm{u_{0}}_{L^{2}} + \nrm{f}_{L^{1}([0, \dlt]; L^{2})} \right)
\end{equation*}
for some $C_{0} < + \infty$. In fact, via Duhamel's principle, this result is equivalent to Proposition~\ref{prop:Mz-d-norm-infl}. Nonetheless, our techniques generalize easily to other situations when such an equivalence is not obvious, e.g., when $b$ depends on time. 

More interestingly, we also provide some new unconditional quantitative lower bounds for \eqref{eq:Mz-d-0}, which are valid way past the trivial $O(\tfrac{1}{\lmb})$ time scale (where $\lmb$ is the initial characteristic frequency), up to a time when the $L^{2}$-norm may grow at a \emph{quantitative} rate depending on $\lmb$; see Propositions~\ref{prop:Mz-1} and \ref{prop:Mz-d} below. These results should be contrasted with the proofs of Proposition~\ref{prop:Mz-d-norm-infl} and \cite[\S VII.2]{Mz2} rely on \emph{qualitative} contradiction arguments up to $O(\frac{1}{\lmb})$ time scales.

\subsection{Organization of the paper} The Schr\"odinger and KdV-type equations are treated respectively in Sections \ref{sec:Schro} and \ref{sec:KdV}. In Appendix~\ref{sec:Mz}, we prove Proposition~\ref{prop:Mz-d-norm-infl} and related results for the linear nondegenerate Schr\"odinger-type equation \eqref{eq:Mz-d-0}.

%
%
%

\subsection*{Acknowledgments}
We thank J. Hunter for bringing the equation \eqref{eq:DS} to our attention. We also thank Benjamin Harrop-Griffiths and Jeremy Marzuola for several insightful discussions on degenerate dispersive equations. I.-J.~Jeong has been supported  by the National Research Foundation of Korea grant (No. 2022R1C1C1011051). S.-J.~Oh was partially supported the Sloan Research Fellowship and the National Science Foundation CAREER Grant under NSF-DMS-1945615. I.-J.~Jeong acknowledges the support from the KIAS associate member program.

\section{Schr\"odinger-type equations}\label{sec:Schro}

This section is organized as follows. To motivate our approach, we analyze in Section~\ref{subsec:approx-linear} a model problem \eqref{eq:lin-approx} derived from \eqref{eq:DS}. In Section~\ref{subsec:DS-bg}, we study the properties of linearly degenerate solutions -- typically denoted by $f$ -- and in Section~\ref{subsec:DS-wp}, we construct degenerating wave packets for the linearized equation around $f$. In Section~\ref{subsec:DS-mee-gee}, we establish a modified and generalized (bilinear) energy estimates for the perturbation (solving the nonlinear difference equation) around $f$. Finally, in Sections~\ref{subsec:DS-unbounded} and \ref{subsec:DS-nonexist}, we prove Theorems~\ref{thm:ill-posed-unbounded-gDS} and \ref{thm:ill-posed-nonexist-gDS}, respectively.

\subsection{Degenerating wave packets for model linear equation}\label{subsec:approx-linear}

To motivate what is to follow, consider the case \eqref{eq:DS} (i.e., \eqref{eq:gDS} with $\alp_{1} = \bt_{1} = 1$ and $\mu_{1} = 0$), which we recall here for convenience: 
\begin{equation*}\tag{\ref{eq:DS}}
i\rd_t \phi + \rd_x(|\phi|^2\rd_x\phi) = 0.
\end{equation*}
 We note that when the domain is taken to be $\bbR$, $f(t,x) = xe^{2it}$ is formally a solution to \eqref{eq:DS}.\footnote{In this section, we take the domain to be $\bbR$ rather than $\bbT$.} Indeed, for initial data $f_0(x)$ which is $x$ in some neighborhood of the origin, it is not difficult to show that any hypothetical smooth solution to \eqref{eq:DS} with initial data $f_0(x)$ should approximate $xe^{2it}$ uniformly for small $|x|$ and $|t|$.
The goal of this section is to sketch a proof of the fact that the linearized equation around the explicit solution $xe^{2it}$ is \textit{ill-posed}. To derive the equation, we write \begin{equation}\label{eq:linear-ansatz}
		\begin{split}
			\phi(t,x) = xe^{2it}+ \tphi(t,x).
		\end{split}
	\end{equation} Plugging in \eqref{eq:linear-ansatz} into \eqref{eq:DS} and dropping quadratic or higher terms in $\tphi$, we obtain that \begin{equation}\label{eq:lin-approx}
\begin{split}
	i\rd_t \tphi + \calL[\tphi] = 0,\quad \calL[\cdot] = \rd_x(x^2 \rd_x(\cdot))+2x\rd_x\Re(\cdot). 
\end{split}
\end{equation} 
An important observation regarding the operator $\calL$ is that for any sufficiently regular $v$, we have the estimate \begin{equation}\label{eq:L-apriori}
\begin{split}
	\left|\brk{|x|^{\frac{1}{2}}\calL[v],|x|^{\frac{1}{2}}v} \right| \le C \nrm{|x|^{\frac{1}{2}}v}_{L^2}^{2},
\end{split}
\end{equation} which suggests that the correct way to measure regularity for solutions of \eqref{eq:lin-approx} is to use $|x|^{\frac{1}{2}}$--weighted spaces\footnote{As we shall see below, the exponents $2$ and $\frac{1}{2}$ in $x e^{2i t}$ and $\abs{x}^{\frac{1}{2}}$, respectively, should be replaced by appropriate constants for \eqref{eq:gDS} in general.}. To this end, we set $\nrm{v}_{L^2_w} = \nrm{|x|^{\frac12}v}_{L^2}$.
The main result in this section is the following.
	\begin{proposition}\label{prop:schro-lin-model}
		The equation \eqref{eq:lin-approx} is ill-posed in $L^{2}$. More specifically, for any profile $g_{0} \in C^\infty(\frac12,1)$, \emph{any} $L^\infty_tL^2_w$-solution $\tphi_{(\lmb)}$ to \eqref{eq:lin-approx} with initial data \begin{equation*}
			\begin{split}
				\tphi_{(\lmb), 0}(x) = g_{0} (x)\exp(i\lmb\ln |x|) , \quad \lmb<0 
			\end{split}
		\end{equation*} satisfies the growth \begin{equation*}
		\begin{split}
			\nrm{\tphi_{(\lmb)}}_{L^{2}}(t) \ge c_{0} \exp(|\lmb|t), \quad  \mbox{for any } 0<t<T,
		\end{split}
	\end{equation*} with constants $c_0, T>0$ depending only on $g_{0}$. 
	\end{proposition}

\begin{remark}
	By an $L^\infty_tL^2_w$ solution to \eqref{eq:lin-approx}, we mean a weak solution $\tphi$ which satisfies the bound \begin{equation*}
		\begin{split}
			\nrm{\tphi}_{L^2_w}(t)\le \exp(Ct)\nrm{\tphi_0}_{L^2_w}
		\end{split}
	\end{equation*} where $C>0$ is the constant from \eqref{eq:L-apriori} and attains the initial data in the weak sense.   Existence of an $L^\infty_tL^2_w$ solution given an $L^2_w$ initial data follows from a standard argument involving Aubin--Lions lemma (see \cite[Appendix A]{JO1} for instance). Note that $\nrm{\tphi_{(\lmb),0}}_{L^{2}}, 	\nrm{\tphi_{(\lmb),0}}_{L^{2}_{w}} \lesssim 1$ uniformly in $\lmb$. While we cannot rule out the possibility of non-uniqueness, the above result applies to \textit{all} $L^\infty_tL^2_w$ solutions. 
\end{remark}

In following the proof, the reader may already find Items 1 (Degenerating wave packets) and 2 (Modified energy estimate and duality testing argument) in Section~\ref{subsec:discussion-pf} expanded in detail. See also the remarks following the proof, which discuss additional ideas that go into the proof of Theorems~\ref{thm:ill-posed-unbounded-gDS} and \ref{thm:ill-posed-nonexist-gDS}.

\begin{proof}
We demonstrate how to construct approximate solutions to \eqref{eq:lin-approx}, from which Proposition \ref{prop:schro-lin-model} naturally follows. To begin with, we make a change of variables $y = \ln x$ for $x\ge 0$. Then using $x\rd_x = \rd_y$, \eqref{eq:lin-approx} transforms into \begin{equation*}
\begin{split}
i\rd_t\tphi + \rd_{yy}\tphi + \rd_y\tphi + 2\rd_y\Re(\tphi) = 0. 
\end{split}
\end{equation*}  
Defining $\varphi = e^{y} \tphi$, 
\begin{equation*}
i \rd_{t} \varphi + \rd_{yy} \varphi + \rd_{y} \overline{\varphi} - 2 \varphi - \overline{\varphi}= 0.
\end{equation*}
We then introduce 
\begin{equation} \label{eq:varphi2psi}
\psi = \varphi + A \overline{\varphi},
\end{equation}
where $A= \frac{1}{2}\rd_y^{-1}$ is (formally) an operator of order $-1$. Then
\begin{align*}
i \rd_{t} \psi
&= - \rd_{yy} \varphi + A \rd_{yy} \overline{\varphi}
- \rd_{y} \overline{\varphi} + A \rd_{y} \varphi \\
&\relphantom{=}
- 2 \varphi - \overline{\varphi} + 2 A \overline{\varphi} + A \varphi \\
&= - \rd_{yy} \psi + A \rd_{yy} \overline{\varphi} + \rd_{yy} A \overline{\varphi}
- \rd_{y} \overline{\varphi}  \\
&\relphantom{=}
+ A \rd_{y} \varphi - 2 \varphi - \overline{\varphi} + 2 A \overline{\varphi} + A \varphi.
\end{align*}
Then we see that 
\begin{equation*}
i \rd_{t} \psi + \rd_{yy} \psi = A \rd_{y} \varphi - 2 \varphi - \overline{\varphi} + 2 A \overline{\varphi} + A \varphi.
\end{equation*}
Since the right hand side is of order zero, it suggests that a degenerating wave packet $\varphi$ may be constructed by taking $\psi$ to be an approximate wave packet solution to the one-dimensional Schr\"odinger equation, and then going back to $\varphi$. More precisely, take
\begin{equation}\label{eq:psi-app}
\psi^{app}_{(\lmb)} (t, y) = e^{i \lmb y - i \lmb^{2} t} a_{0}(y - 2 \lmb t),
\end{equation} where we fix $a_{0}$ to be $C^\infty$--smooth and supported in $\{ -2 < y < -1\}$. We need to take $\lmb < 0$, so that the support of $\phi^{app}(t,\cdot)$ is confined to $\{y<-1\}$ for all $t\ge0$. 
To invert \eqref{eq:varphi2psi}, we wish to take $\varphi \approx \psi - A \overline{\psi}$. Since $A = \frac{1}{2} \rd_{y}^{-1}$ acts like $- \frac{1}{2 i \lmb}$ on $\overline{\psi}$, we are motivated to take
\begin{align}\label{eq:varphi-app}
\varphi^{app}(t,y) &= \psi^{app}_{(\lmb)} (t, y) + \frac{1}{2 i \lmb}  \overline{\psi^{app}_{(\lmb)}}(t, y) 
\end{align} and then set \begin{equation*}
\begin{split}
\tphi^{app} = e^{-y} \varphi^{app} = e^{-y}\left( e^{i\lmb(y-\lmb t)}  a_{0}(y-2\lmb t)+ \frac{1}{2i\lmb}e^{-i\lmb(y-\lmb t)}\overline{a_{0}}(y-2\lmb t)  \right).
\end{split}
\end{equation*} Returning to the $x$-coordinates and defining the error by $\err_{\tphi}=[i\rd_t+\calL] \tphi^{app}$, we have \begin{equation}\label{eq:error-estimate}
\begin{split}
\nrm{ \err_{\tphi}(t)}_{L^2_{w}} \lesssim \nrm{a_{0}}_{H^2_x} ,\qquad t\ge0 ,
\end{split}
\end{equation} uniformly in $\lmb$. In this sense, $\tphi^{app}$ is an approximate solution of \eqref{eq:lin-approx}. Moreover, $\tphi^{app}$ itself satisfies the bound $\nrm{  {\tphi^{app}}(t)}_{L^2_{w}} \lesssim \nrm{a_{0}}_{L^2_x} $. The last key property is degeneration: with a weight higher than $|x|^{\frac{1}{2}}$, $\tphi^{app}(t)$ decays  in the $O(|\lmb|^{-1})$--timescale: for example, with the weight $|x|$, we have \begin{equation}\label{eq:degeneration}
\begin{split}
\nrm{|x|\tphi^{app}(t,x)}_{L^2_x} \lesssim e^{-|\lmb| t} \nrm{a_{0}}_{L^2_x}. 
\end{split}
\end{equation} Interpolating \eqref{eq:degeneration} with $L^2_{w}$--estimate shows that $\nrm{\tphi^{app}(t,\cdot)}_{L^2}\gtrsim e^{|\lmb| t}$. Now, let $\tphi$ be an $L^\infty_{t}L^{2}_w$ solution of \eqref{eq:lin-approx}. Then, with a direct computation, we have the \textit{generalized energy estimate} for the weighted $L^2$--estimate (see \S\ref{subsec:gee}): \begin{equation*}
\begin{split}
	\frac{\ud}{\ud t} \brk{ |x|^{\frac12} \tphi, |x|^{\frac12} \tphi^{app} } =  \brk{ |x|^{\frac12} \tphi, |x|^{\frac12}\err_{\tphi} },
\end{split}
\end{equation*}  which gives, together with \eqref{eq:error-estimate},  \begin{equation*}
\begin{split}
	\Re \,  \brk{ |x|^{\frac12} \tphi, |x|^{\frac12} \tphi^{app} }(t)  \ge \Re \,  \brk{ |x|^{\frac12} \tphi, |x|^{\frac12} \tphi^{app} }(t=0) - C t \nrm{\tphi}_{ L^\infty_{t}L^{2}_w }  \nrm{a_{0}}_{H^2_x} . 
\end{split}
\end{equation*} At the initial time, by choosing $a_{0}$ in a way depending only on $g_{0}$, we can guarantee that $ \Re \,  \brk{ |x|^{\frac12} \tphi, |x|^{\frac12} \tphi^{app} }(t=0) \ge \frac12 \nrm{\tphi_{0}}_{L^2_w} \nrm{\tphi^{app}_{0}}_{L^{2}_{w}}$. Then, for $0<t < C\nrm{a_{0}}_{H^2_x}/(4 \nrm{\tphi_{0}}_{L^2_w} )$, we obtain with \eqref{eq:degeneration} that \begin{equation*}
\begin{split}
	\frac{1}{C} e^{-|\lmb| t} \nrm{a_{0}}_{L^2 } \nrm{\tphi(t)}_{L^{2}}  \ge  \Re \,  \brk{ \tphi, |x| \tphi^{app} }(t)  = \Re \,  \brk{ |x|^{\frac12} \tphi, |x|^{\frac12} \tphi^{app} }(t)  \ge\frac14 \nrm{\tphi_{0}}_{L^2_w} \nrm{\tphi^{app}_{0}}_{L^{2}_{w}},
\end{split}
\end{equation*} which gives the claimed exponential growth of $ \nrm{\tphi(t)}_{L^{2}} $. 
\end{proof}

\begin{remark}[Illposedness in $H^{m}$ for $m > 0$]
In fact, a small modification of the above proof shows that \eqref{eq:lin-approx} is ill-posed  in $H^{m}$ for $m > 0$. More precisely, the following growth occurs:
\begin{equation*}
\nrm{{\rd_{x}^{m}} \tphi_{(\lmb)}}_{L^{2}}(t) \ge c_{0} \exp({(1+2m)} |\lmb|t), \quad  \mbox{for any } m \ge 0  \mbox{ and } 0 < t < T,
\end{equation*}
with $c_{0}, T > 0$ depending only on $g_{0}$ and $m$.

We now sketch the needed modification; see Section~\ref{subsec:DS-unbounded} for the complete proof. We would like to modify the last part of the proof of Proposition~\ref{prop:schro-lin-model} using ``differentiation by parts'': under the assumption that $\rd_{x}^{-m} ( |x| \tphi^{app}) \in L^{2}$, \begin{equation*}
\begin{split}
 \nrm{\rd_{x}^{m} \tphi}_{L^2 } \nrm{ \rd_{x}^{-m} ( |x| \tphi^{app}) }_{L^{2}}  \ge  \Re \,  \brk{ \rd_{x}^{m} \tphi,  \rd_{x}^{-m} ( |x| \tphi^{app}) }(t)  = \Re \,  \brk{ |x|^{\frac12} \tphi, |x|^{\frac12} \tphi^{app} }(t) .
\end{split}
\end{equation*} Now the point is that $\rd_{x}^{-1} = x\rd_{y}^{-1} = e^{y}\rd_{y}^{-1}$ and $y \simeq -2|\lmb|t$ on the support of $\tphi^{app}(t)$, which gives a faster rate of degeneration $\nrm{ \rd_{x}^{-m} ( |x| \tphi^{app}) }_{L^{2}} \lesssim e^{- (2m+1)|\lmb| t }$. This gives the claimed lower bound for $ \nrm{\rd_{x}^{m} \tphi}_{L^2 }$. In general, there could be some low frequency part of $|x|\tphi^{app}$ which does not degenerate, and for this reason we introduce a decomposition of $|x|\tphi^{app}$ into high and low frequency parts in the actual proof in Section~\ref{subsec:DS-unbounded}.
\end{remark}

\begin{remark}[Additional ideas in the proof of Theorems~\ref{thm:ill-posed-unbounded-gDS} and \ref{thm:ill-posed-nonexist-gDS}]
For a general equation of the form \eqref{eq:gDS}, we do not have access to a stationary solution with a linear degeneracy in general (furthermore, we shall also require that $f_{0}$ be compactly supported, which rules out $x e^{2 i t}$, too). Hence, we shall carry out the above analysis (degenerating wave packet construction, modified energy estimate and duality) where the background solution $f$ is merely a regular (most likely) \emph{time-dependent} solution to \eqref{eq:gDS}, which has compactly supported initial data $f_{0}$ with a linear degeneracy, in place of $x e^{2 i t}$.

Theorem~\ref{thm:ill-posed-unbounded-gDS} is proved by considering a perturbation $f + \tld{\phi}$ of such an $f$, and arguing that if $f + \tld{\phi}$ exists as a regular solution (i.e., if we are in the second case in Theorem~\ref{thm:ill-posed-unbounded-gDS}), then the above growth mechanism for $\tld{\phi}$ can be justified. To prove Theorem~\ref{thm:ill-posed-nonexist-gDS}, we consider initial data $\tld{\phi}_{0}$ consisting of a superposition of an infinitude of configurations as above (i.e., $\sum_{k} (f_{k, 0} + \tld{\phi}_{k, 0})$, where $f_{k, 0}$ has a linear degeneracy and $\tld{\phi}_{k, 0}$ is a degenerating wave packet adapted to $f_{k, 0}$) with unbounded rates (i.e., the initial frequencies of the degenerating wave packets are unbounded), disjoint supports (i.e., $\set{\supp f_{k, 0} \cup \supp \phi_{k, 0}}_{k}$ is pairwise disjoint), yet with an $\eps$-small $C^{m_{0}}$-norm. Then we perform a contradiction argument: if a regular solution $\phi$ to such initial data exists, then we may justify the growth mechanism (as in Proposition~\ref{prop:schro-lin-model}), which is absurd. For details, see Sections~\ref{subsec:DS-unbounded} and \ref{subsec:DS-nonexist} below.
\end{remark}

\subsection{Properties of a regular linearly degenerate solution} \label{subsec:DS-bg}

We shall assume that there exists a smooth solution to \eqref{eq:gDS} which is linearly degenerate, and analyze its properties. To be precise, let $f: [0,\dlt]\times [-x_{1},x_{1}]\rightarrow\bbC$ be a  $L^{\infty}([0,\dlt];C^{3,1}([{x_{0}}-x_{1},{x_{0}+}x_{1}]))$ solution to  {\eqref{eq:gDS}} with some $x_1, \dlt>0$ satisfying  \begin{equation*}
\begin{split}
f_{0} \in C^{3,1}([{x_{0}}-x_{1}, {x_{0}}+x_{1}]),\quad f_0(x_0) =0, \quad |f'_0(x_0)| > 0  
\end{split}
\end{equation*} at $t=0$ for some $x_0\in\bbT$. 

Owing to the symmetries of  {\eqref{eq:gDS}} (translation and phase rotation), as well as its behavior under the transformation $\phi \mapsto c \phi$, we may assume without loss of generality that $x_0 = 0$ and $f'_{0}(0) = 1$. Then, from the equation it is easy to see that, on the time interval $[0,\dlt]$, \begin{equation*}
\begin{split}
f(t,0) = 0
\end{split}
\end{equation*} and   \begin{equation*}
\begin{split}
i \frac{\ud}{\ud t} f'(t,0) = -(\alp_{1} + \bt_{1}) |f'(t,0)|^2 f'(t,0),
\end{split}
\end{equation*} which implies in particular that 
\begin{equation*}
|f'(t,0)| = 1, \quad \hbox{ and } |f(t,x)| = x + O(|x|^2) \hbox{ uniformly in $t$.}
\end{equation*}
More generally, we have the following lemma.
\begin{lemma} \label{lem:gDS-vanishing-propagate}
Let $s \geq 2$ be an integer and let $f \in C_{t}([0, \dlt]; C^{s-1, 1}(\bbT))$ be a solution to \eqref{eq:gDS}. Then 
\begin{enumerate}
\item The zero set of $f(t, x)$ is preserved in time, i.e., $a \in \bbT$ is a zero of $f(0, x)$ if and only if it is a zero for $f(t, x)$ for all $t \in [0, \dlt]$.
\item Let $a \in \bbT$ be a zero of $f(0, x)$. Then $\set{\rd_{x}^{k} f(t, a)}_{k=0}^{s-1}$ is determined by the initial data at $x = a$, i.e., $\set{\rd_{x}^{k} f(0, a)}_{k=0}^{s-1}$.
\end{enumerate}
\end{lemma}
Here, the important point is that, thanks to the regularity assumption, $f(0, x)$ vanishes at least linearly at each zero $x = a$, which is \emph{critical} for \eqref{eq:gDS} in the senses discussed in \S\ref{subsec:discussion}.

\begin{proof}
By the regularity assumption (in particular, that $s \geq 2$), it follows from \eqref{eq:gDS} that $\abs{\rd_{t} f(t, x)} \leq C \abs{f(t, x)}$; hence the first statement follows. To prove the second statement, consider a zero $a$ of $f(0, x)$. Without any loss of generality, we may assume that $a = 0$. By the assumption and Taylor expansion, we have
\begin{align*}
	f(t, x) &= \sum_{k=1}^{s-1} \frac{1}{k!} \rd_{x}^{k} f(t, 0) x^{k} + O(\abs{x}^{s}), \\
	f_{x}(t, x) &= \sum_{k=0}^{s-2} \frac{1}{k!} \rd_{x}^{k+1} f(t, 0) x^{k} + O(\abs{x}^{s-1}), \\
	f_{xx}(t, x) &= \sum_{k=0}^{s-3} \frac{1}{k!} \rd_{x}^{k+2} f(t, 0) x^{k} + O(\abs{x}^{s-2}),
\end{align*}
where the implicit constants depend only on $\nrm{f}_{L^{\infty}_{t} C^{s-1, 1}_{x}}$ (and in particular are independent of $(t, x)$). Plugging this into \eqref{eq:gDS} and matching the coefficients of $x^{k}$ for $k = 1, \ldots, s-1$, we formally obtain a determined system of first-order ODE's for $\set{\rd_{x}^{k} f(t, 0)}_{k=1}^{s-1}$; here, the fact that $f(t, x)$ vanishes at least linearly is crucially used to ensure that no $\rd_{x}^{k} f(t, 0)$ with $k > s-1$ arises. Indeed, these ODEs may be justified in the sense of distributions by testing \eqref{eq:gDS} against a test function of the form $\eta(t) (-1)^{k} \rd_{x}^{k} \chi_{\eps}(x)$, where $\eta \in C^{\infty}_{c}(0, \dlt)$, $\chi \in C^{\infty}_{c}(-\frac{1}{2}, \frac{1}{2})$ with $\int \chi = 1$ and $\chi_{\eps}(x) = \eps^{-1} \chi(\eps^{-1} x)$. By the uniqueness of this ODE system, the desired statement follows.
\end{proof}


From now on, given $f_{0}$ which are linearly degenerate at $x =0$ and $f_{x}(0,0)$ positive real, we are going to take $0<x_{1}<1$ smaller if necessary, so that \begin{equation}\label{eq:gDS-x1}
	\begin{split}
		\left(\sup_{ x \in [-x_{1},x_{1}] } |f_{xx}(0,x)|\right) x_{1} < \frac{ f_{x}(0,0) }{2}
	\end{split}
\end{equation} holds. In particular, we have that \begin{equation*}
\begin{split}
	\frac12 f_{x}(0,0) <|f_{x}(0,x)| < 2f_{x}(0,0),\quad \mbox{for all} \quad  x\in [-x_{1},x_{1}]. 
\end{split}
\end{equation*}

\begin{proposition}\label{prop:C31}
	Let $f \in L^{\infty}_{t}([0, \dlt]; C^{3, 1}([-x_1,x_1])$ be a solution to \eqref{eq:gDS}, and set $M = \nrm{f}_{L^\infty([0,\dlt];C^{3,1})}$. Then, we have the pointwise bounds \begin{equation}\label{eq:gDS-pointwise}
		\begin{split}
			|f(t,x)| \le |f_{0}(x)| \exp( CM^{2}t ) 
		\end{split}
	\end{equation} and 
	\begin{equation}\label{eq:obs}
		\begin{split}
				\left| \rd_{t} (|f(t,x)|^{2})  \right| \le CM\exp( C M^{2} \dlt)\left( 1 + (f_{x}(0,0))^{-1} M \right)^{3} \left( |f_{0}(x)|^{3} + t M^{3}|f_{0}(x)|^{2} \right)
		\end{split}
	\end{equation} for all $t\in [0,\dlt]$ and $|x|\le x_{1}$. 
\end{proposition}

\begin{proof} We first note directly from \eqref{eq:gDS} that $\left| \rd_{t} |f(t,x)| \right| \le C \nrm{f}^{2}_{L^\infty_t C^{1,1}_x} |f(t,x)|$ holds, which gives \eqref{eq:gDS-pointwise}. Now note that $f\in L^{\infty}([0,\dlt];C^{3,1}([-x_{1},x_{1}]))$ implies, {via \eqref{eq:gDS}}, that \begin{equation}\label{eq:gDS-time-der}
\begin{split}
	|\rd_{t} f(t,x)|  \le C|f(t,x)|\nrm{f}^{2}_{L^\infty_t C^{1,1}_{x}} , \qquad 	|\rd_{tt} f(t,x)| &\le C|f(t,x)|^{2} \nrm{f}^{3}_{L^\infty_{t} C^{3,1}_{x}}. 
\end{split}
\end{equation} 
Then, the Taylor expansion in time of $f(t,x)$ gives \begin{equation*}
\begin{split}
	|f(t,x)|^{2} & = |f_{0}(x)|^{2} + 2 \Re\left( \overline{f_{0}(x)}  \int_{0}^{t} (\rd_{t}f)(t',x)\,\ud t'  \right) + \left|  \int_{0}^{t} (\rd_{t}f)(t',x)\,\ud t' \right|^{2} . 
\end{split}
\end{equation*} Taking the time derivative, \begin{equation}\label{eq:rd-t-sq}
\begin{split}
	\rd_{t} (|f(t,x)|^{2}) = 2\Re\left(  \overline{f_{0}(x)}  (\rd_{t}f)(t,x) \right) + 2 \Re\left(  \overline{(\rd_{t}f)(t,x)}\int_{0}^{t} (\rd_{t}f)(t',x)\,\ud t'  \right). 
\end{split}
\end{equation} Using \eqref{eq:gDS-pointwise}, the last term in \eqref{eq:rd-t-sq} is bounded by \begin{equation*}
\begin{split}
	C\nrm{f}^{4}_{L^\infty_t C^{3,1}_x} |f(t,x)|\int_0^{t} |f(t',x)| \, \ud t' \le CM^{4}t\exp(CM^{2}\dlt) |f_{0}(x)|^{2}. 
\end{split}
\end{equation*} For the other term in the right hand side of \eqref{eq:rd-t-sq}, we further rewrite as 
\begin{equation*}
	\begin{split}
		 2\Re\left(  \overline{f_{0}(x)}  (\rd_{t}f)(t,x) \right)  =  2\Re\left(  \overline{f_{0}(x)}  (\rd_{t}f)(0,x) \right)  +  2\Re\left(  \overline{f_{0}(x)}  \int_{0}^{t} (\rd_{tt}f)(t',x) \, \ud t' \right),  
	\end{split}
\end{equation*} and note that the last term is bounded using \eqref{eq:gDS-time-der} by $C\exp(CM^{2}\dlt)  M^{3} t|f_{0}(x)|^{3}.$ 
On the other hand, the first term on the right hand side equals \begin{equation*}
\begin{split}
	&\Im \left(  \overline{f_{0}(x)} \left(  \abs{f_{0}(x)}^{2} \rd_{xx} f_{0}(x) + \alp_{1} f_{0}(x) \abs{\rd_{x} f_{0}(x)}^{2} + \bt_{1} \br{f_{0}(x)} (\rd_{x} f_{0}(x))^{2} + \mu_{1} \abs{f_{0}(x)}^{2} f_{0}(x) \right) \right) \\
	&\qquad =  \bt_{1}\Im \left(  \overline{f_{0}(x)}^{2} (\rd_{x} f_{0}(x))^{2}  \right) + O(\nrm{f_{0}}_{C^{1,1}})|f_{0}(x)|^{3},
\end{split}
\end{equation*} and we see that the leading term in the Taylor expansion of $\overline{f_{0}(x)}^{2} (\rd_{x} f_{0}(x))^{2} $ is purely real, with remainder bounded by \begin{equation*}
\begin{split}
	C( f_{x}(0,0)^3 |x|^3 \nrm{f_0}_{\dot{C}^{1,1}} + f_{x}(0,0)^{2}  |x|^{4} \nrm{f_0}_{\dot{C}^{1,1}}^{2}  + f_{x}(0,0)^{3} |x|^{5} \nrm{f_0}_{\dot{C}^{1,1}}^{2} ) \le C(1+ f_{x}(0,0)^{-1} M)^{3}|f_{0}(x)|^{3} ,
\end{split}
\end{equation*} where we have used $|x|<x_{1}$ and the smallness of $x_{1}$ from \eqref{eq:gDS-x1}. Collecting the bounds, we obtain the proposition. \end{proof}

\subsection{Degenerating wave packets for the linearized equation}\label{subsec:DS-wp}

In this subsection, our goal is to construct approximate solutions, called degenerating wave packets, for the linearization of \eqref{eq:gDS} around a (possibly hypothetical) regular linearly degenerate solution, which possess the desired degeneration property that is responsible for the illposedness of \eqref{eq:gDS}; see Proposition~\ref{prop:DS-WP} below.

\subsubsection{Properties of degenerating wave packets}

Given a smooth solution $f$ to \eqref{eq:gDS}, let us write  $\phi = f + \tphi$ where $\phi$ is another smooth solution to \eqref{eq:gDS}. The equation for $\tphi$ is given by 
\begin{equation}\label{eq:deg-Schro-pert}
\begin{split}
&i \rd_{t} \tphi + \abs{f}^{2} \rd_{xx} \tphi + \alp_{1} f (\br{\rd_{x} f} \rd_{x} \tphi + \rd_{x} f \rd_{x} \br{\tphi}) + 2 \bt_{1} \br{f} \rd_{x} f \rd_{x} \tphi
+ V_{f} \tphi + W_{f} \br{\tphi} = Q_{f}[\tphi], 
\end{split}
\end{equation} 
with \begin{equation}\label{eq:Qf}
\begin{split}
	V_{f} &=  \br{f} \rd_{xx} f + \alp_{1} \abs{\rd_{x} f}^{2} + 2 \mu_{1} \abs{f}^{2}, \\
	W_{f} &= f \rd_{xx} f + \bt_{1} (\rd_{x} f)^{2}  {+} \mu_{1} f^{2}, \\
	Q_f[\tphi] &= - (\br{f} \tphi + f \br{\tphi}) \rd_{xx} \tphi
	- \alp_{1} \tphi (\rd_{x} \br{f} \rd_{x} \tphi + \rd_{x} f \rd_{x} \br{\tphi}) - \alp_{1} f \abs{\rd_{x} \tphi}^{2}
	- 2 \bt_{1} \br{\tphi} \rd_{x} f \rd_{x} \tphi - \bt_{1} \br{f} (\rd_{x} \tphi)^{2} \\
	&\relphantom{=}
	- \abs{\tphi}^{2} \rd_{xx} \tphi - \alp_{1} \tphi \abs{\rd_{x} \tphi}^{2} - \bt_{1} \br{\tphi} (\rd_{x} \tphi)^{2}  {- \mu_0 |\tphi|^{2}\tphi - 2\mu_{1}f |\tphi|^2 - \mu_{1} \br{f} (\tphi)^{2}}. 
\end{split}
\end{equation} Note that $Q_f[\tphi]$ is at least quadratic in $\tphi$ and its derivatives. 
Dropping the right hand side, we obtain the linearized equation around $f$: 
\begin{equation}\label{eq:deg-Schro-genlin}
\begin{split}
i \rd_{t} \tphi + \abs{f}^{2} \rd_{xx} \tphi + \alp_{1} f (\br{\rd_{x} f} \rd_{x} \tphi + \rd_{x} f \rd_{x} \br{\tphi}) + 2 \bt_{1} \br{f} \rd_{x} f \rd_{x} \tphi
+ V_{f} \tphi + W_{f} \br{\tphi} = 0.
\end{split}
\end{equation} We now state the key proposition of this section, which shows properties of degenerating wave packets for \eqref{eq:deg-Schro-genlin}. We introduce the notation
\begin{equation*}
	\nrm{g}_{W^{s,p}_{(L)}} = \sum_{j = 0}^{s} \nrm{(L \rd_{x})^{j} g}_{L^{p}_{x}}
\end{equation*} and write $H^s_{(L)}$ when $p=2$. 

\begin{proposition}\label{prop:DS-WP}
	Let $f \in L^{\infty}([0,\dlt]; {C}^{s_{0}-1,1}([0,x_1]))$ be a solution to \eqref{eq:gDS} with $s_{0}\ge 4$ satisfying \begin{equation}\label{eq:f-assumption-lin-deg}
		\begin{split}
			f(0,0) = 0, \quad f'(0,0) ={ A},
		\end{split}
	\end{equation} for some $A>0$. By taking $x_{1}<1$ small if necessary, assume that \eqref{eq:gDS-x1} holds. Then, to any $\lmb\le -1$ and a $C^\infty$-smooth complex-valued profile $g_0$ supported in $(\frac{1}{2}x_1,x_1)$, we may associate a function $\tphi^{app}_{(\lmb)} = \tphi^{app}_{(\lmb)}[g_0,f]$ defined in $[0,\dlt]\times\bbR$ satisfying the following properties: \begin{itemize}
		\item (linearity) the map $g_0 \mapsto \tphi^{app}_{(\lmb)}$ is (real) linear; 
		\item (support property) $\supp(\tphi^{app}_{(\lmb)}[g_0]) \subset (0, e^{-|\lmb| A^{2} t} x_1)$;
		\item (initial data) for any $1 \le p \le \infty$, \begin{equation*}
			\begin{split}
				\frac{1}{C}\nrm{g_{0}}_{L^2}\le A^{\sgm_{c}-1} \nrm{ |f|^{-\sgm_{c}} \tphi^{app}(0,x) }_{L^{2}} & \le C\nrm{g_{0}}_{L^2}, \\
				A^{\sgm_{c}-1} \nrm{ |f|^{\sgm_{c}} \tphi^{app}(0,x) }_{L^{p}} & \le  C x_{1}^{\frac1{p}-\frac12} \nrm{g_{0}}_{L^{p}} ; 
			\end{split}
		\end{equation*}
		\item (regularity) for $0\le n \le s_{0}-2$, we have \begin{equation}\label{eq:ds-wp-reg}
		\begin{split}
			\nrm{{\abs{f}}^{-\sgm_{c}} ({\abs{f}} \rd_x)^n \tphi^{app}_{(\lmb)}(t,x) }_{L^2} \le C_{f,\dlt} A^{-\sgm_{c} + n + 1}|\lmb|^{n}\nrm{g_0}_{ H^{n}_{(x_{1})}  }, \quad t\le\min\{ A^{-2}|\lmb|^{-\frac12},\dlt \} ;
		\end{split}
		\end{equation}
		\item (degeneration) 
for any $1 \leq p \leq 2$, $ {0 \le s \le s_{0}-2}$, and $\gmm' \geq -s-\frac{1}{p}+\frac{1}{2}$, we have
		\begin{equation} \label{eq:ds-wp-deg}
		\abs{f}^{-\sgm_{c}+\gmm'} \tphi^{app}_{(\lmb)} = \rd_{x}^{s}\left(\frac{\abs{f}^{\gmm'+s-\frac{1}{2}} }{i^{s} \lmb^{s} (1+\abs{f} \rd_{x} S)^{s}}  \psi^{app}_{(\lmb)}\right) + \abs{f}^{-\sgm_{c}+\gmm'} \tphi^{small}_{(\lmb)},
\end{equation}
		for some $\psi^{app}_{(\lmb)}$, $\tphi^{small}_{(\lmb)}$, and $S$,  where $\psi^{app}_{(\lmb)}$ is independent of $p$, $s$ and $\gmm'$, and 
		\begin{align}
		\left\| \frac{\abs{f}^{\gmm'+s-\frac{1}{2}}}{\lmb^{s} (1+\abs{f} \rd_{x} S)^{s}}\psi^{app}_{(\lmb)}(t, x) \right\|_{L^{p}} &\leq C_{f,\dlt}^{1+\gmm'} A^{\gmm'+s+\frac12} {|\lmb|^{-s}} e^{-2 \abs{\lmb}(\gmm'+s+\frac{1}{p}-\frac{1}{2}) A^{2} t} \nrm{g_{0}}_{L^{2}},  \label{eq:ds-wp-deg-main} \\
		\nrm{\abs{f}^{-\sgm_{c}} \tphi^{small}_{(\lmb)}(t, x)}_{L^{2}} &\leq C_{f,\dlt} A^{-\sgm_{c} + 1} {|\lmb|^{-1}} \nrm{g_{0}}_{ H^{s}_{(x_{1})}  } \label{eq:ds-wp-deg-small}
		\end{align} for $ t \leq \min\set{A^{-2}\abs{\lmb}^{-\frac{1}{2}}, \dlt}$, after taking $\dlt>0$ smaller in a way that $\dlt \nrm{f}^{2}_{L^\infty([0,\dlt];C^{1,1})}$ is small in terms of $A^{-1}\nrm{f}_{L^{\infty} ([0,\dlt];C^{3,1})}$; 
		\item (error estimate) defining the error $\err[\tphi^{app}_{(\lmb)}]$ by the left hand side of \eqref{eq:deg-Schro-genlin} with $\tphi = \tphi^{app}_{(\lmb)}$, we have the estimate \begin{equation}\label{eq:ds-wp-err}
		\begin{split}
		\nrm{\abs{f}^{-\sgm_{c}} \err[\tphi^{app}_{(\lmb)}](t)}_{L^2} \le C_{f,\dlt} A^{-\sgm_{c}+3} \nrm{g_0}_{ H^{2}_{(x_{1})}  },\quad t\le\min\{ A^{-2}|\lmb|^{-\frac12},\dlt \}.
		\end{split}
		\end{equation} 
	\end{itemize} 
In the above estimates, the constant $C_{f,\dlt}$ satisfies \begin{equation}\label{eq:Cf-bound}
		\begin{split}
			C_{f,\dlt} \le C_{0} \left( 1 + A^{-1}\nrm{f}_{L^\infty_{t}C^{s_{0}-1,1}} \right)^{N_{0}} \exp(C_{0} \nrm{f}_{L^\infty_{t}C^{s_{0}-1,1}}^{2} \dlt )
		\end{split}
	\end{equation} for some $C_{0}, N_{0}>0$ depending on $\alp_{1}, \bt_{1}, \mu_{1}$ and $s_{0}$ but not on $f$ and $x_{1}$.
\end{proposition}

We fix $A = 1$ and prove Proposition~\ref{prop:DS-WP} in the remainder of this subsection. In the general case, given $f$ we can define $\tilde{f}(t,x) := A^{-1} f(A^{-2}t,x) $ which is another solution $L^\infty_t C^{s-1,1}$ solution to \eqref{eq:gDS} satisfying $\tilde{f}_{x}(0,0)=1$. Then, we simply define \begin{equation*}
\begin{split}
	\tilde{\phi}^{app}_{(\lmb)}[g_{0},f] (t,x) :=  \tilde{\phi}^{app}_{(\lmb)}[ A g_{0}, \tilde{f}](A^{2}t,x),
\end{split}
\end{equation*} and verify the claimed properties of $\tilde{\phi}^{app}_{(\lmb)}[g_{0},f]$ using those for $\tilde{\phi}^{app}_{(\lmb)}[A g_{0},\tilde{f}]$. In the proof, it will be seen that $|f|\rd_{x}S$ remains invariant under this rescaling. 

\subsubsection{Renormalization and wave packet construction} \label{subsec:DS-renrm}
With $x_1>0$ given in Proposition \ref{prop:DS-WP}, we define the variable $y$ for $t \in [0,\dlt]$ and $x \in (0,x_1]$ by \begin{equation*}
\begin{split}
y(t,x)= -\int_{x}^{x_1} \frac{1}{|f(t,x')|}\,  \ud x' \le 0.
\end{split}
\end{equation*} For each $t\ge0$, the inverse of $x \mapsto y(t,x)$ is denoted by $x=x(t,y)$. From $|f(t,x)| = x + O(|x|^{2})$, we have that \begin{equation} \label{eq:sch-xy-compare}
\begin{split}
{y(t,x)}- {\ln\frac{x}{x_1}} = B(t,x), \qquad x(t,y) = x_{1} e^{y-B}, \qquad |B(t,x)|\le Cx_{1} \nrm{f}_{L^{\infty}_t C^{1,1}}. 
\end{split}
\end{equation}  
Using $|f|\rd_x = \rd_y$, we rewrite \eqref{eq:deg-Schro-genlin} in $(t,y)$-coordinates: 
\begin{equation}\label{eq:deg-Schro-genlin-prime}
i \rd_{t} \tphi + i h \rd_{y} \tphi + \rd_{y}^{2} \tphi +  \frac{\alp_{1} f \br{\rd_{y} f} + 2 \bt_{1} \br{f} \rd_{y} f - \abs{f} \rd_{y} \abs{f}}{\abs{f}^{2}} \rd_{y} \tphi + \alp_{1} \frac{f \rd_{y} f}{\abs{f}^{2}} \rd_{y} \br{\tphi}
+ V_{f} \tphi + W_{f} \br{\tphi} = 0.
\end{equation} Here, we have introduced $h(t,y) = \rd_{t}y$ so that $\rd_t \tphi(t,x) = \rd_t \tphi(t,y)+  h(t,y)\rd_y \tphi(t,y).$ Now define \begin{equation*}
\begin{split}
  {G(t,y) = (-\sgm_{c}+\frac{1}{2}) \ln \abs{f} (t, y)}, \quad (\rd_{y}G)(t,y) =  {(-\sgm_{c}+\tfrac{1}{2})} \frac{\Re(\overline{f}\rd_yf)}{|f|^2} (t,y)
\end{split}
\end{equation*} and introducing the conjugation $\varphi = e^{G}\tphi,$
we obtain  {(recall from \eqref{eq:gDS-w} that $\sgm_{c} = - (\frac{\alp_{1}}{2} +\bt_{1} - 1)$)}
\begin{equation}\label{eq:deg-Schro-genlin-renorm}
\begin{split}
&i\rd_t\varphi + \rd_{yy} \varphi +  {\alp_{1}} \frac{f\rd_yf}{|f|^2} \overline{\rd_y \varphi} +  \left( {(-\alp_{1} + 2 \bt_{1})}\frac{\Im(\overline{f}\rd_y f)}{|f|^2} + h\right) i \rd_y\varphi = \calB_0[\varphi] 
\end{split}
\end{equation} with \begin{equation*}
\begin{split}
\calB_0[\varphi] =  {i} (\rd_tG)\varphi + ( \rd_{yy}G + (\rd_{y}G)^2 )\varphi 
 {+ \alp_{1} \frac{f \rd_{y} f}{\abs{f}^{2}} \rd_{y}G \br{\varphi} + \left( (-\alp_{1} + 2 \bt_{1})\frac{\Im(\overline{f}\rd_y f)}{|f|^2} + h \right) i (\rd_{y}G) \varphi - V_{f} \varphi - W_{f} \br{\varphi}.}
\end{split}
\end{equation*} Note that  {the} terms in $\calB_0$ do not contain derivatives of $\varphi$. To handle the term containing $\overline{\rd_y\varphi}$ in the left hand side of \eqref{eq:deg-Schro-genlin-renorm}, we make yet another change of variables: introducing formally \begin{equation*}
\begin{split}
\psi = \varphi +  {\frac{\alp_{1}}{2}}\frac{ f\rd_y f}{|f|^2} \rd_y^{-1} \overline{\varphi} , 
\end{split}
\end{equation*} we have that \eqref{eq:deg-Schro-genlin-renorm}  turns into \begin{equation}\label{eq:deg-Schro-genlin-renorm2}
\begin{split}
i\rd_t\psi + \rd_{yy}\psi + \left( {(-\alp_{1} + 2 \bt_{1})}\frac{\Im(\overline{f}\rd_y f)}{|f|^2} +  {h}\right) i \rd_y\psi = \cdots ,
\end{split}
\end{equation} where  {the} terms  {on} the right hand side do not contain any derivatives of $\psi$. Motivated by this computation, we construct a wave packet approximate solution for \eqref{eq:deg-Schro-genlin-renorm} by starting with a wave packet for the preceding equation for $\psi$, then coming back to $\varphi$. More precisely, given $g_{0}(x)$ as in Proposition \ref{prop:DS-WP}, we take \begin{equation*}
	\begin{split}
		a_{0}(y) = x_{1}^{\frac12} g_{0}(x(0,y)),
	\end{split}
\end{equation*} which is supported in $y\in {(-\frac{1}{2} \ln 2, 0)}$ by \eqref{eq:sch-xy-compare}. For each $\lmb < 0$, we define  \begin{equation}\label{eq:psi-app-def}
\begin{split}
\psi^{app}_{(\lmb)}(t,y) := e^{i\lmb(y-\lmb t)} a_{(\lmb)}(t,y) 
\end{split}
\end{equation} where $a_{(\lmb)}(t,y)$ is the unique solution to \begin{equation}\label{eq:WKB-a}
\begin{split}
\rd_t a_{(\lmb)} + 2\lmb\rd_y a_{(\lmb)} = \frac{\lmb}{i} \left(  {(-\alp_{1} + 2 \bt_{1})} \frac{\Im(\overline{f}\rd_y f )}{|f|^2} + h \right) a_{(\lmb)}
\end{split}
\end{equation} with initial data $a_{(\lmb)}(0,y)=a_0(y)$. Given $\psi^{app}$, set \begin{equation}\label{eq:phi-app-def}
	\begin{split}
		\varphi^{app}_{(\lmb)} = \psi^{app}_{(\lmb)} + \frac{ {\alp_{1}}}{2i\lmb} \frac{f\rd_y f}{|f|^2} \overline{\psi^{app}_{(\lmb)}} , 
	\end{split}
\end{equation} which will be shown to be an approximate solution to \eqref{eq:deg-Schro-genlin-renorm}. Finally, the degenerating wave packet is defined by  \begin{equation}\label{eq:tphi-app-def}
\begin{split}
	\tphi^{app}_{(\lmb)}[g_0,f] = e^{-G} \varphi^{app}_{(\lmb)} = |f|^{\sgm_{c} - \frac12} \varphi^{app}_{(\lmb)} .  
\end{split}
\end{equation}

\subsubsection{Proof of Proposition \ref{prop:DS-WP}}

Now that we have defined the wave packet solution, let us proceed to confirm the properties stated in Proposition \ref{prop:DS-WP}.

\medskip

\noindent \textit{Linearity and support property}. From the definition, linearity is clear. Furthermore, note from \eqref{eq:tphi-app-def}, \eqref{eq:phi-app-def}, \eqref{eq:psi-app-def} that the support of $\tphi^{app}_{(\lmb)}(t,\cdot)$ coincides with that of $a_{(\lmb)}(t,\cdot)$. (From now on, we shall suppress from writing out the subscript $\lmb$.) On the other hand, note the following formulae for $a$:
\begin{align}
	a(t, y) &= e^{i \lmb S(t, y)} a_{0}(y-2\lmb t), \label{eq:WKB-a-expr} \\ 
	S(t, y) &= \int_{0}^{t} \left((-\alp_{1} + 2 \bt_{1}) \frac{\Im (\br{f} \rd_{y} f)}{\abs{f}^{2}} + h \right)(t', y - 2 \lmb(t-t')) \, \ud t'. \label{eq:WKB-S-expr}
\end{align}
Since $\lmb<0$, the support of $a(t,\cdot)$ is contained in the interval  {$(-\frac{1}{2} \ln 2 +2\lmb t, 2\lmb t) \subseteq (-\infty,0)$} for $t\ge0$, which verifies the support property of $\tphi^{app}$ via \eqref{eq:sch-xy-compare}.

\medskip

\noindent \textit{Regularity estimates}. To begin with, we obtain estimates on $h:= \rd_{t}y$. Recalling \eqref{eq:rd-t-sq}, we have  \begin{equation*}
		\begin{split}
			h & = -\int_{x}^{x_1} \rd_{t} \left( \frac{1}{|f(t,x')|} \right) \,  \ud x' = \int_{x}^{x_{1}} \frac{\rd_{t} (|f|^{2})}{|f|^{3}} \, \ud x' , \qquad \rd_{y}h = |f|\rd_{x}h = -  \frac{\rd_{t} (|f|^{2})}{|f|^{2}}. 
		\end{split}
	\end{equation*} Applying \eqref{eq:gDS-pointwise} and \eqref{eq:obs}, we obtain the pointwise estimates \begin{equation}\label{eq:h-est}
	\begin{split}
		|h| \le C_{f,\dlt}(1 + t\ln\frac{1}{x})x_{1}, \qquad |\rd_{y}h| \le C_{f,\dlt}(x+t). 
	\end{split}
\end{equation}
We now estimate $a$. Observing that the right hand side in \eqref{eq:WKB-a} is purely imaginary, \begin{equation*}
\begin{split}
\frac{1}{2}\frac{\ud}{\ud t} |a|^2(t,y) =   -2\lmb \Re(\rd_{y} a \overline{a}) (t,y), \quad \mbox{which gives} \quad \frac{\ud}{\ud t} \nrm{a}_{L^2_y}^{ {2}}=0.
\end{split}
\end{equation*} Now, taking a $y$-derivative, we see that \begin{equation*}
\begin{split}
\frac{1}{2}\frac{\ud}{\ud t} \nrm{\rd_{y} a_{(\lmb)}}_{L^2}^{2}   \le C_{f,\dlt} |\lmb| (e^{-|\lmb| t} + t)\nrm{a}_{L^2_y}  \nrm{\rd_{y}a}_{L^2_y} 
\end{split}
\end{equation*} where we have used \begin{equation}\label{eq:rd-y-S}
\begin{split}
	\left| \rd_{y}  \left(  {(-\alp_{1} + 2 \bt_{1})} \frac{\Im(\overline{f}\rd_y f )}{|f|^2} + h \right) \right| \le C_{f,\dlt}  (x(t,y) + t ) \le C_{f,\dlt}  ( \exp(-|\lmb|t) + t )
\end{split}
\end{equation} on the support of $a(t,\cdot)$. This estimate follows from \eqref{eq:h-est} and \begin{equation*}
\begin{split}
	\left| \rd_{y} \frac{\Im(\overline{f}\rd_y f )}{|f|^2} \right| \le 	\left| \rd_{x} \rd_{y} \frac{\Im(\overline{f}\rd_y f )}{|f|^2} \right| x . 
\end{split}
\end{equation*} Therefore, by integrating in time, we obtain $$\nrm{\rd_{y}a(t)}_{L^2_y} \le C_{f,\dlt} \nrm{a_0}_{H^{1}_{y}}$$ uniformly in $\lmb$, for $(t,\lmb)$ satisfying $t\le {|\lmb|^{-\frac12}} $. A similar argument applies to the estimate of $\rd_{y}^{k}a$, as long as $k \le s_{0}-2$; one can proceed by an induction in $k$, using the bound   \begin{equation*}
\begin{split}
	\left| \rd_{y}^{k}  \left(  {(-\alp_{1} + 2 \bt_{1})} \frac{\Im(\overline{f}\rd_y f )}{|f|^2} + h \right) \right| \le C_{f,\dlt}  (x(t,y) + t ) \le C_{f,\dlt}  ( \exp(-|\lmb|t) + t )
\end{split}
\end{equation*} on the support of $a(t,\cdot)$. The estimate for $|\rd_{y}^{k}h|$ readily follows from the explicit decomposition $ \rd_{y} h = h_1 + t h_2$, where $h_1$ and $h_2$ are $L^{\infty}_t C^{s_{0}-2,1}$-smooth functions defined by \begin{equation*}
\begin{split}
	h_1(t,x) = 2|f|^{-2}\Re\left(  \overline{f_{0}(x)}  (\rd_{t}f)(t,x) \right), \quad h_{2}(t,x)= 2 |f|^{-2} \Re\left(  \overline{(\rd_{t}f)(t,x)} \, \frac{1}{t} \int_{0}^{t} (\rd_{t}f)(t',x)\,\ud t'  \right). 
\end{split}
\end{equation*} Hence we conclude \begin{equation}\label{eq:wp-a}
\begin{split}
\nrm{a(t)}_{H^{k}_y} \le C_{f,\dlt}\nrm{a_0}_{H^{k}_y}, \quad 0\le t\le\min\{{|\lmb|^{-\frac12}},\dlt \}  .
\end{split}
\end{equation} In what follows, we shall restrict the variable $t$ to $[0,\min\{|\lmb|^{-\frac12},\dlt \}]$.

\medskip

\noindent \textit{Initial data and regularity estimates}. At the initial time, from $a_{0}(y) = x_{1}^{\frac12} g_{0}(x(0,y))$ we have that \begin{equation*}
	\begin{split}
		\int |a_{0}(y)|^{2} \, \ud y = \int x_{1} |f_{0}(x)|^{-1} |g_{0}(x)|^{2} \, \ud x
	\end{split}
\end{equation*} and note that the right hand side is equivalent up to constants with $\nrm{g_{0}}_{L^2_{x}}^{2}$. This gives the claimed initial data estimate in the case $p = 2$, and the case of general $p$ can be proved similarly. Next, with $\rd_{y} = |f_{0}(x)|\rd_{x}$ at the initial time, we note the bound \begin{equation*}
\begin{split}
	|\rd_{y}^{k}a_{0} (x) | \le C_{k} x_{1}^{\frac12} \left( \sum_{j=1}^{k} \nrm{f_{0}}_{C^{k-2,1}}^{k-j} |f_{0}(x)|^{j} |\rd_{x}^{j}g_{0}(x)| \right)
\end{split}
\end{equation*} which gives \begin{equation}\label{eq:a-g}
\begin{split}
	\nrm{a_{0}}_{H^{k}_{y}} \le C_{k} ( 1 + \nrm{f_{0}}_{C^{k-2,1}})^{k-1} \nrm{g_{0}}_{H^{k}_{(x_{1})}}, \quad k \le s_{0}-1.
\end{split}
\end{equation} 
Let us now check the regularity estimate \eqref{eq:ds-wp-reg} in the case $n = 0$: using $|f|^{- \sgm_{c} +\frac{1}{2}}  {=} e^{G}$, \begin{equation*}
	\begin{split}
		\nrm{ \abs{f}^{-\sgm_{c}} \tphi^{app}_{(\lmb)} (t,x) }_{L^2}^2 &= \int_0^{x_1} |\tphi^{app}_{(\lmb)} (t,x)|^2  |f(t, x)|^{- 2\sgm_{c}}\, \ud x = \int_{-\infty}^{0} |\tphi^{app}_{(\lmb)} (t,y)|^2  |f(t,y)|^{- 2 \sgm_{c}+1}  \, \ud y\\
		&=\int_{-\infty}^{0} |\varphi^{app}_{(\lmb)} (t,y)|^2 \, \ud y \le C(1+|\lmb|^{-1})\int_{-\infty}^{0} |\psi^{app}_{(\lmb)} (t,y)|^2 \, \ud y \\
		&\le C\nrm{a(t)}_{L^2_y}^2 \le C\nrm{a_0}_{L^2_y} \le C \nrm{g_0}_{L^2}. 
	\end{split}
\end{equation*} The cases $1 \le n \le s_{0}-2$ can be handled similarly, using \eqref{eq:wp-a} and \eqref{eq:a-g}. 

\medskip

\noindent \textit{Degeneration estimate}. Next, we check the degeneration property \eqref{eq:ds-wp-deg}. To simplify the notation, we introduce the notation 
\begin{equation*}
	H = q O_{k}(a_{0}) \quad \impmi \quad \sup_{t \in [0, \min\set{\lmb^{-\frac{1}{2}}, \dlt}]}\nrm*{\abs{f}^{\frac{1}{2}} \frac{H}{q}}_{L^{2}_{y}} \leq C_{f,\dlt} \nrm{a_{0}}_{H^{k}_{y}}.
\end{equation*}
Note  that $\nrm{\abs{f}^{\frac{1}{2}}(\cdot)}_{L^{2}_{y}} = \nrm{\cdot}_{L^{2}_{x}}$ for each $t$. 
The terms that are abbreviated as $\frac{1}{\lmb} O_{k}(a_{0})$ (for $k \leq s$) will constitute $\abs{f}^{-\sgm_{c}} \tphi^{small}_{(\lmb)}$; the desired estimate \eqref{eq:ds-wp-deg-small} would be an immediate consequence of the $L^{2}$ norm estimate embedded in the $O_{k}(\cdot)$ notation. Recalling the definitions of $\tphi^{app}_{(\lmb)}$, $\varphi^{app}_{(\lmb)}$, and $\psi^{app}_{(\lmb)}$, and arguing as in the proof of the  regularity estimate, we have
\begin{align*}
	\abs{f}^{-\sgm_{c}+\gmm'}\tphi^{app}_{(\lmb)} 
	&= \abs{f(t, y)}^{\gmm' - \frac{1}{2}} \psi^{app}_{(\lmb)} + \frac{\abs{f}^{\gmm'}}{\lmb} O_{0}(a_{0}).
\end{align*}
For the first term, we have
\begin{align*}
	\nrm{\abs{f(t, x)}^{\gmm' - \frac{1}{2}} \psi^{app}_{(\lmb)}(t, x)}_{L^{p}}^{p}
	&\leq C \int_{-\infty}^{0} \abs{\psi^{app}_{(\lmb)}(t, y)}^{p} \abs{f(t, y)}^{p\gmm' - \frac{p}{2} + 1} \, \ud y \\
	&\leq C \left( \int_{-\infty}^{0} \abs{\psi^{app}_{(\lmb)} (t, y)}^{2} \, \ud y\right)^{\frac{p}{2}} \left( \int_{\supp \psi^{app}_{(\lmb)}(t, \cdot)} \abs{f(t, y)}^{\frac{p}{1-\frac{p}{2}} \gmm'+1} \, \ud y \right)^{1-\frac{p}{2}} \\
	&\leq C \nrm{a_{0}}_{L^{2}_{y}}^{p} \left( \int_{\supp \psi^{app}_{(\lmb)}(t, \cdot)} \abs{f(t, y)}^{\frac{p}{1-\frac{p}{2}} \gmm'+1} \, \ud y \right)^{p(\frac{1}{p}-\frac{1}{2})},
\end{align*}
so it remains to estimate the last factor. Note that, since $\abs{f}^{-1} \rd_{y} \abs{f} = \rd_{x} \abs{f} = 1 + O(x)$ and $x \leq x_{1} e^{\frac{y}{2}}$ for $y \in (-\infty, 0)$, we have
\begin{equation*}
	\abs{f(t, y)} \leq C e^{y} \quad \hbox{ for } y \in (-\infty, 0). 
\end{equation*}
Using the support property $\supp \psi^{app}_{(\lmb)}(t, \cdot) \subseteq (-\infty, - 2 \abs{\lmb} t)$, we see that
\begin{align*}
	\left( \int_{\supp \psi^{app}_{(\lmb)}(t, \cdot)} \abs{f(t, y)}^{\frac{p}{1-\frac{p}{2}} \gmm'+1} \, \ud y \right)^{\frac{1}{p}-\frac{1}{2}}
	\aleq e^{- 2 \abs{\lmb} (\gmm' + \frac{1}{p} - \frac{1}{2})t}.
\end{align*}
Hence, the desired estimate \eqref{eq:ds-wp-deg} in the case $s = 0$ now follows.

To treat the cases $s > 0$, we begin by recalling that $\psi^{app}_{(\lmb)} = \psi^{app} = e^{i \lmb (y - \lmb t + S(t, y))} a_{0}(y - 2 \lmb t)$. Note the identity
\begin{align*}
	e^{i \lmb (y - \lmb t + S)} = \frac{\abs{f}}{i \lmb (1 + \rd_{y} S)}\left(\frac{1}{\abs{f}}\rd_{y}\right) e^{i \lmb (y - \lmb t + S)}.
\end{align*}
For the expression $\rd_{y} S$ in the denominator, recalling \eqref{eq:WKB-S-expr} and \eqref{eq:rd-y-S}, we have
\begin{equation*}
	\abs{\rd_{y} S} \le C_{f,\dlt} t {x},
\end{equation*} and in particular note that {$1 + \rd_{y} S \geq \frac{1}{2}$} when $t$ is sufficiently small, which can be arranged by taking $\dlt>0$ smaller. 
Commuting $\frac{1}{\abs{f}}\rd_{y}$ (which equals $\rd_{x}$ in the $(t, x)$--coordinates) outside, we have
\begin{align*}
	\abs{f(t, y)}^{\gmm' - \frac{1}{2}} \psi^{app}_{(\lmb)}
	&= \frac{1}{\abs{f}} \rd_{y} \left( \frac{\abs{f(t, y)}^{\gmm' + 1 - \frac{1}{2}}}{i \lmb (1+\rd_{y} S)} \psi^{app}_{(\lmb)} \right) + \frac{\abs{f}^{\gmm'}}{\lmb} O_{1} (a_{0}).
\end{align*}
By arguing as in the case of $s = 0$, the expression inside the parentheses can be shown to obey the degeneration bound \eqref{eq:ds-wp-deg-main}. The cases $s > 1$ are handled similarly.

\medskip

\noindent \textit{Error estimate}. To begin with, at the level of $\psi^{app}_{(\lmb)}$, the point of choosing $a(t,  {y})$ as the solution of \eqref{eq:WKB-a} is to have \begin{equation}\label{eq:O-nota-1}
\begin{split}  i\rd_t\psi^{app}_{(\lmb)} + \rd_{yy}\psi^{app}_{(\lmb)} + \left( {(-\alp_{1} + 2 \bt_{1})} \frac{\Im(\overline{f}\rd_y f)}{|f|^2} + h\right)i \rd_y\psi^{app}_{(\lmb)}  = O_{2}(a_{0}), 
\end{split}
\end{equation} which can be checked with a direct computation using \eqref{eq:deg-Schro-genlin-renorm2}. We now see that $\varphi^{app}_{(\lmb)}$ is an approximate solution to \eqref{eq:deg-Schro-genlin-renorm}, which is motivated by the following heuristics: \begin{equation*}
\begin{split}
\varphi^{app}_{(\lmb)} \simeq \psi^{app}_{(\lmb)} -  {\frac{\alp_{1}}{2}} \frac{f\rd_y f}{|f|^2} \rd_y^{-1} \overline{\psi^{app}_{(\lmb)}} \simeq  \psi^{app}_{(\lmb)} + \frac{ {\alp_{1}}}{2i\lmb} \frac{f\rd_y f}{|f|^2} \overline{\psi^{app}_{(\lmb)}} . 
\end{split}
\end{equation*}  To this end, \eqref{eq:O-nota-1} gives \begin{equation*}
\begin{split}
	-i\overline{\rd_t\psi^{app}_{(\lmb)}} + \overline{\rd_{yy}\psi^{app}_{(\lmb)}} - \left( {(-\alp_{1} + 2 \bt_{1})}\frac{\Im(\overline{f}\rd_y f)}{|f|^2} + h\right) i \overline{\rd_y\psi^{app}_{(\lmb)}} =O_{2}(a_{0})
\end{split}
\end{equation*} and from this it is not difficult to see that \begin{equation}\label{eq:O-nota-2}
\begin{split}
	\frac{1}{2i\lmb}[i\rd_t + \rd_{yy}] \overline{\psi^{app}_{(\lmb)}} = \frac{1}{i\lmb}\rd_{yy}\overline{\psi^{app}_{(\lmb)}} + O_{2}(a_0) = -\rd_y \overline{\psi^{app}_{(\lmb)}} + O_{2}(a_0),
\end{split}
\end{equation}
 {so that}\begin{equation*}
\begin{split}
	(i \rd_{t} + \rd_{yy}) \left(\frac{\alp_{1}}{2 i \lmb} \frac{f \rd_{y} f}{\abs{f}^{2}} \br{\psi^{app}_{(\lmb)}} \right)
	= - \alp_{1} \frac{f \rd_{y} f}{\abs{f}^{2}} \rd_{y} \br{\psi^{app}_{(\lmb)}} + O_{2}(a_{0}).
\end{split}
\end{equation*}
 Using \eqref{eq:O-nota-1}, \eqref{eq:O-nota-2}, and \begin{equation*}
\begin{split}
	\left( {(-\alp_{1}+2\bt_{1})}\frac{\Im(\overline{f}\rd_y f)}{|f|^2} + h\right)i \rd_y\left(\frac{ {\alp_{1}}}{2i\lmb} \frac{f\rd_y f}{|f|^2} \overline{\psi^{app}_{(\lmb)}}  \right) = O_{2}(a_{0}), \quad	\frac{ {\alp_{1}^{2}} f\rd_y f}{2i\lmb|f|^2} \rd_y\left(   \frac{\overline{f\rd_y f}}{|f|^2} \psi^{app}_{(\lmb)} \right)  = O_{2}(a_0), 
\end{split}
\end{equation*} we simplify \begin{equation*}
\begin{split}
	&\left[i\rd_t +\rd_{yy}+\left( {(-\alp_{1}+2\bt_{1})} \frac{\Im(\overline{f}\rd_y f)}{|f|^2} + h\right) i \rd_y\right]\varphi^{app}_{(\lmb)}  +  {\alp_{1}}\frac{f\rd_yf}{|f|^2} \overline{\rd_y \varphi^{app}_{(\lmb)}} \\
	&\quad = \left[i\rd_t +\rd_{yy}+\left( {(-\alp_{1}+2\bt_{1})} \frac{\Im(\overline{f}\rd_y f)}{|f|^2} + h\right) i \rd_y\right]\left( \psi^{app}_{(\lmb)} + \frac{ {\alp_{1}}}{2i\lmb} \frac{f\rd_y f}{|f|^2} \overline{\psi^{app}_{(\lmb)}}  \right) \\
	&\quad \relphantom{=} +  {\alp_{1}} \frac{f\rd_yf}{|f|^2} \overline{\rd_y\psi^{app}_{(\lmb)}} - \frac{ {\alp_{1}^{2}} f\rd_y f}{2i\lmb|f|^2} \rd_y\left(   \frac{\overline{f\rd_y f}}{|f|^2} \psi^{app}_{(\lmb)} \right) \\
	&\quad = -  {\alp_{1}} \frac{f\rd_yf}{|f|^2} \rd_y \overline{\psi^{app}_{(\lmb)}} +  {\alp_{1}} \frac{f\rd_yf}{|f|^2} \rd_y \overline{\psi^{app}_{(\lmb)}} +O_{2}(a_0)=O_{2}(a_0). 
\end{split}
\end{equation*} Moreover, it is easy to see that $\calB_0[\varphi^{app}_{(\lmb)}] = O_{2}(a_0)$ and finally, the error estimate \eqref{eq:ds-wp-err} follows from \begin{equation*}
\begin{split} \left[i\rd_t +\rd_{yy}+\left( {(-\alp_{1} + 2 \bt_{1})} \frac{\Im(\overline{f}\rd_y f)}{|f|^2} + h\right)i \rd_y\right]\varphi^{app}_{(\lmb)}  +  {\alp_{1}} \frac{f\rd_yf}{|f|^2} \overline{\rd_y \varphi^{app}_{(\lmb)}}  - \calB_0[\varphi^{app}_{(\lmb)}] = O_{2}(a_{0})
\end{split}
\end{equation*} and \eqref{eq:a-g}. This completes the proof of Proposition \ref{prop:DS-WP}. \hfill \qedsymbol

\subsection{Modified and generalized energy estimates} \label{subsec:DS-mee-gee}

In \S\ref{subsec:modified-energy-estimate} and \S\ref{subsec:gee}, we establish the modified and generalized (or bilinear) energy estimates, respectively, that we shall need the proofs of Theorems~\ref{thm:ill-posed-unbounded-gDS} and \ref{thm:ill-posed-nonexist-gDS}.

\subsubsection{Modified energy estimate}\label{subsec:modified-energy-estimate}
Assume that $f$ and $\phi = f+\tphi$ are solutions to \eqref{eq:DS} on some time interval. Then, recall that $\tphi$ solves \begin{equation}\label{eq:plugin2}
	\begin{split}
		 {i \rd_{t} \tphi + \abs{f}^{2} \rd_{xx} \tphi + \alp_{1} f (\br{\rd_{x} f} \rd_{x} \tphi + \rd_{x} f \rd_{x} \br{\tphi}) + 2 \bt_{1} \br{f} \rd_{x} f \rd_{x} \tphi
+ V_{f} \tphi + W_{f} \br{\tphi} = Q_{f}[\tphi],}
	\end{split}
\end{equation}
where  {$V_{f}$, $W_{f}$ and $Q_f[\cdot]$ are} defined in \eqref{eq:Qf}. For a solution $\tphi$ of \eqref{eq:plugin2}, we have the following \begin{proposition}\label{prop:mee-Schro}
	 Let $f \in L^{\infty}([0,\dlt];C^{ {s_{c}-1},1}))$ be a solution to \eqref{eq:gDS} and  $\tphi \in L^{\infty}([0,\dlt'];C^{ {s_{c}-1},1}))$ to \eqref{eq:plugin2} for some $0<\dlt' \le \dlt$. When $\sgm_{c} \geq \frac{1}{2}$, assume furthermore that at every zero $a$ of $f_{0}$, we have $\rd_{x}f_{0}(a) \neq 0$ and $\tphi_{0}(x)$ vanishes up to order $\lfloor \sgm_{c} - \frac{1}{2} \rfloor$ at $a$. Then, on $t\in[0,\dlt']$, we have \begin{equation}\label{eq:mee}
	\begin{split}
	\nrm{|f(t)|^{-\sgm_{c}}\tphi(t)}_{L^2} \le \nrm{|f_0|^{-\sgm_{c}} \tphi_0}_{L^2}\exp\left( C(\nrm{f}_{L^\infty_tC^{1,1}}^2 + \nrm{\tphi}_{L^\infty_tC^{1,1}}^2)t \right)
	\end{split}
	\end{equation}  {where $\sgm_{c}$ is as in \eqref{eq:gDS-w} and} $C>0$ is an absolute constant. 
\end{proposition}

\begin{proof}
We first present a formal computation without worrying about the finiteness of the modified energy and the validity of integration by parts, and discuss its justification below. We  compute that \begin{equation*}
\begin{split}
\frac{\ud}{\ud t}\int |\tphi|^2 \, |f|^{ - 2 \sgm_{c}}\, \ud x = \int|\tphi|^2\, \rd_t|f|^{- 2 \sgm_{c}} \, \ud x + \int  {\rd_{t} |\tphi|^{2}} \,|f|^{- 2 \sgm_{c}}\, \ud x.
\end{split}
\end{equation*} The first term can be bounded using the pointwise inequality \begin{equation*}
\begin{split}
\left| \frac{\rd_t|f|}{|f|} \right| \lesssim \nrm{f}_{L^\infty_tC^{1,1}}^2.
\end{split}
\end{equation*} To handle the second term, we write \begin{equation}\label{eq:tphi-sq}
\begin{split}
\rd_t|\tphi|^2 &= \Re\left( i \abs{f}^{2} \rd_{xx} \tphi \br{\tphi}\right) + \alp_{1} \Re\left(i f (\rd_{x} \br{f} \rd_{x} \tphi + \rd_{x} f \rd_{x} \br{\tphi}) \br{\tphi} \right) + 2 \bt_{1} \Re\left(i \br{f} \rd_{x} f \rd_{x} \tphi \br{\tphi} \right) \\
&\relphantom{=} + \Re\left( i V_{f} \tphi \br{\tphi}\right) + \Re\left( i W_{f} (\br{\tphi})^{2}\right) - \Re \left( iQ_{f}[\tphi] \br{\tphi}\right).
\end{split}
\end{equation}
We multiply both sides by $|f|^{- 2\sgm_{c}}$ and integrate in $x$. From the first term on the right-hand side, we obtain, after an integration by parts,
\begin{equation} \label{eq:tphi-sq-1}
\begin{aligned}
\int \Re\left( i \abs{f}^{2} \rd_{xx} \tphi \br{\tphi}\right) \abs{f}^{- 2 \sgm_{c}} \, \ud x
&= -\int \abs{f}^{2 - 2\sgm_{c}}\Re(i \rd_{x} \tphi \rd_{x} \br{\tphi}) \, \ud x - (2 - 2\sgm_{c})\int \abs{f}^{1 - 2\sgm_{c}} \rd_{x} \abs{f} \Re(i \rd_{x} \tphi \br{\tphi}) \\
&= - (2 - 2\sgm_{c}) \int \abs{f}^{- 2\sgm_{c}} \Re(\br{f} \rd_{x} f) \Re(i \rd_{x} \tphi \br{\tphi}).
\end{aligned}
\end{equation}
From the second and third terms on the right-hand side of \eqref{eq:tphi-sq}, we have
\begin{align*}
& \int \left[ \alp_{1} \Re\left(i f (\rd_{x} \br{f} \rd_{x} \tphi + \rd_{x} f \rd_{x} \br{\tphi}) \br{\tphi} \right) + 2 \bt_{1} \Re\left(i \br{f} \rd_{x} f \rd_{x} \tphi \br{\tphi} \right)\right] \abs{f}^{- 2 \sgm_{c}} \, \ud x \\
&= (\alp_{1} + 2 \bt_{1}) \int \abs{f}^{- 2\sgm_{c}} \Re(\br{f} \rd_{x} f) \Re(i \rd_{x} \tphi \br{\tphi}) \, \ud x
+ \frac{\alp_{1} - 2 \bt_{1}}{2} \int \abs{f}^{- 2\sgm_{c}} \Im(\br{f} \rd_{x} f) \rd_{x}\abs{\tphi}^{2} \, \ud x \\
&\relphantom{=}
+ \frac{\alp_{1}}{2} \int \abs{f}^{- 2 \sgm_{c}} \Re\left(i  {f \rd_{x} f} \rd_{x} (\br{\tphi})^{2} \right) \, \ud x \\
&= (\alp_{1} + 2 \bt_{1}) \int \abs{f}^{- 2 \sgm_{c}} \Re(\br{f} \rd_{x} f) \Re(i \rd_{x} \tphi \br{\tphi}) \, \ud x
- \frac{\alp_{1} - 2 \bt_{1}}{2} \int \left( \rd_{x} (\abs{f}^{- 2 \sgm_{c}} \Im(\br{f} \rd_{x} f))\right) \abs{\tphi}^{2} \, \ud x \\
&\relphantom{=}
- \frac{\alp_{1}}{4} \int \left(i \rd_{x} (\abs{f}^{- 2 \sgm_{c}}  {f \rd_{x} f}) \right) (\br{\tphi})^{2} \, \ud x
+ \frac{\alp_{1}}{4} \int \left( i \rd_{x} (\abs{f}^{- 2 \sgm_{c}}  {\overline{f \rd_{x} f}})\right)  (\tphi)^{2} \, \ud x.
\end{align*}
By our choice of $\sgm_{c}$ in \eqref{eq:gDS-w}, the first term on the right-hand side cancels exactly with \eqref{eq:tphi-sq-1}. The remaining terms are estimated from the above by $C \nrm{f}_{L^\infty_tC^{1,1}}^2 \int \abs{\tphi}^{2} \abs{f}^{- 2 \sgm_{c}} \, \ud x$. Next, it is easy to see that
\begin{align*}
\left| \int \Re (i V_f \tphi \br{\tphi}) \abs{f}^{- 2 \sgm_{c}} \,\ud x \right|
+ \left|  \int \Re (i W_f \br{\tphi} \br{\tphi}) \abs{f}^{- 2 \sgm_{c}} \,\ud x \right|
\lesssim \nrm{f}_{L^\infty_tC^{1,1}}^2 \int |\tphi|^2 \, |f|^{- 2 \sgm_{c}}\, \ud x.
\end{align*}
It remains to estimate $\int \Re(i Q_{f}[\tphi] \br{\tphi}) \abs{f}^{- 2 \sgm_{c}} \, \ud x$. The contribution of any term with at least one factor of $\tphi$ (without any derivatives) may be easily estimated by $(\nrm{f}_{C^{1,1}}^2 + \nrm{\tphi}_{C^{1,1}}^2) \int |\tphi|^2 \, |f|^{- 2 \sgm_{c}}\, \ud x$. Recalling the expression for $Q_{f}$ from \eqref{eq:Qf}, we may estimate
\begin{align*}
	\left| \int \Re(i Q_{f}[\tphi] \br{\tphi}) \abs{f}^{- 2 \sgm_{c}} \, \ud x \right|
	&\aleq (\nrm{f}_{C^{1,1}}\nrm{\tphi}_{C^{1,1}} + \nrm{\tphi}_{C^{1,1}}^2) \int |\tphi|^2 \, |f|^{- 2 \sgm_{c}}\, \ud x \\
	&\relphantom{\aleq} + \left( \int \abs{\rd_{x} \tphi}^{4} \abs{f}^{2 - 2\sgm_{c}}\right)^{\frac{1}{2}} \left(\int |\tphi|^2 \, |f|^{- 2 \sgm_{c}}\, \ud x\right)^{\frac{1}{2}}.
\end{align*}
Integrating by parts and using H\"older's inequality, we have
\begin{align*}
\int (\rd_{x} \tphi)^{4} \abs{f}^{2 - 2\sgm_{c}}
&= \int \tphi (\rd_{x} \tphi)^{2} (- 3\rd_{xx} \tphi \abs{f} - (2 - 2 \sgm_{c}) \rd_{x} \tphi \rd_{x} \abs{f} )\abs{f}\abs{f}^{- 2 \sgm_{c}} \, \ud x \\
&\leq C \nrm{\tphi}_{C^{1, 1}} \nrm{f}_{C^{0, 1}} \left( \int \tphi^{2} \abs{f}^{- 2 \sgm_{c}} \, \ud x \right)^{\frac{1}{2}} \left( \int (\rd_{x} \tphi)^{4} \abs{f}^{2 -  2 \sgm_{c}} \, \ud x \right)^{\frac{1}{2}} .
\end{align*}
Hence,
\begin{align*}
	\left| \int \Re(i Q_{f}[\tphi] \br{\tphi}) \abs{f}^{- 2 \sgm_{c}} \, \ud x \right|
	&\aleq (\nrm{f}_{C^{1,1}}\nrm{\tphi}_{C^{1,1}} + \nrm{\tphi}_{C^{1,1}}^2) \int |\tphi|^2 \, |f|^{- 2 \sgm_{c}}\, \ud x.
\end{align*} 
Collecting all the terms, we conclude that \begin{equation*}
\begin{split}
\left| \frac{\ud}{\ud t}\int |\tphi|^2 \, |f|^{ - 2 \sgm_{c}}\, \ud x  \right| \lesssim \left(\nrm{f}_{L^\infty_tC^{1,1}}^2 +\nrm{\tphi}_{L^\infty_tC^{1,1}}^2\right)\int |\tphi|^2 \, |f|^{- 2 \sgm_{c}}\, \ud x .
\end{split}
\end{equation*} Integrating in time gives \begin{equation}\label{eq:mee2}
\begin{split}
	\nrm{|f(t)|^{- \sgm_{c}} \tphi(t)}_{L^2} \le \nrm{|f_0|^{- \sgm_{c}} \tphi_0}_{L^2}\exp\left( C(\nrm{f}_{L^\infty_tC^{1,1}}^2 +\nrm{\tphi}_{L^\infty_tC^{1,1}}^2)t \right),
\end{split}
\end{equation}  {which is the desired conclusion.}

We now sketch the observations needed to make the above computation rigorous. Note that, in order for \eqref{eq:mee} to be nontrivial, the right-hand side must be finite, i.e., $\nrm{\abs{f_{0}}^{- \sgm_{c}} \tphi_{0}}_{L^{2}} < + \infty$. When $\sgm_{c} \geq \frac{1}{2}$, this implies the vanishing of $\tphi_{0}$ at each zero $a$ of $f$ (which is isolated by the assumption in this case) up to order $\lfloor \sgm_{c} -\frac{1}{2} \rfloor $. Applying Lemma~\ref{lem:gDS-vanishing-propagate} to the $L^{\infty}_{t} ([0, \dlt]; C^{s_{c}-1, 1}))$ solutions $f$ and $f + \tphi$, it follows that the zero set of $f(t, x)$, as well as the nonvanishing of $f'(t, a)$ and the vanishing of $\tphi(t, x)$ up to order $\lfloor \sgm_{c} -\frac{1}{2} \rfloor$  at each zero $a$ of $f$, is preserved in $t \in [0, \dlt]$. As a consequence, $\nrm{\abs{f}^{- \sgm_{c}} \tphi}_{L^{2}} < + \infty$ for every $t \in [0, \dlt]$ as well. Using the vanishing properties of $f$ and $\tphi$ (the latter is needed only when $\sgm_{c} \geq \frac{1}{2}$), the above computation can then be justified.
\end{proof}

\medskip

\noindent The above proposition motivates us to define the following hermitian product and norm (which will be referred to as the modified energy): given some $f$, we shall consider \begin{equation*}
	\begin{split}
				\brk{v,u}_{f}  (t) = \int  |f(t,\cdot)|^{- 2 \sgm_{c}} v(t,\cdot)\overline{u(t,\cdot)} \,\ud x., \qquad \nrm{v}_{f}^2 (t) = \int |f(t,\cdot)|^{- 2 \sgm_{c}} |v(t,\cdot)|^2\, \ud x .
	\end{split}
\end{equation*}

\subsubsection{Generalized (bilinear) energy estimate}\label{subsec:gee}
We proceed to prove the generalized energy estimate. 

\begin{proposition}\label{prop:gDS-gee}
	
	Let $\tphi$ to be a solution of \begin{equation*}
		\begin{split}
			[i\rd_t + \calL_{f}]\tphi = Q_f[\tphi],
		\end{split}
	\end{equation*} where $[i \rd_{t} + \calL_{f}] \tphi$ denotes the left-hand side of \eqref{eq:deg-Schro-pert}, and let $\tphi^{app} = \tphi^{app}[g_{0},f]$ be the degenerating wave packet constructed in Proposition \ref{prop:DS-WP}. Then, we have the following estimate on $t \in [0,\min\{|\lmb|^{-\frac12},\dlt\}]$: \begin{equation}\label{eq:DS-gei}
		\begin{split}
			\left| \frac{\ud}{\ud t} \Re\brk{\tphi,\tphi^{app}}_f  \right| \le \left(C (\nrm{f}_{L^\infty_t C^{1,1}_x}^2+\nrm{\tphi}_{L^\infty_t C^{1,1}_x}^2) \nrm{\tphi^{app}}_{L^2_f} + C_{f,\dlt}A^{- \sgm_{c}+3} \nrm{g_{0}}_{H^{2}_{(x_{1})}} \right) \nrm{\tphi}_{L^2_f}. 
		\end{split}
	\end{equation}
\end{proposition}
\begin{proof} In the proof, the time variable $t$ will be restricted to the interval $[0,\min\{|\lmb|^{-\frac12},\dlt\}]$.
Before we proceed, let us recall that $\calL_{f}$ is given by \begin{equation*}
\begin{split}
	\calL_{f}[\tphi] &:= \abs{f}^{2} \rd_{xx} \tphi + \alp_{1} f (\br{\rd_{x} f} \rd_{x} \tphi + \rd_{x} f \rd_{x} \br{\tphi}) + 2 \bt_{1} \br{f} \rd_{x} f \rd_{x} \tphi
+ V_{f} \tphi + W_{f} \br{\tphi} \\
	& \relphantom{:}=
	\abs{f}^{2} \rd_{xx} \tphi + (\alp_{1} + 2 \bt_{1}) \Re(\br{f} \rd_{x} f) \rd_{x} \tphi 
	+ (-\alp_{1} + 2 \bt_{1}) i \Im (\br{f} \rd_{x} f) \rd_{x} \tphi
	+ \alp_{1} f \rd_{x} f \rd_{x} \br{\tphi} + V_{f} \tphi + W_{f} \br{\tphi} 
\end{split}
\end{equation*} and that $\tphi^{app}$ satisfies $[i\rd_t + \calL_{f}]\tphi^{app} = \err_{\tphi}.$  We compute\footnote{Here, since $\tphi^{app}$ is smooth and compactly supported away from the zeroes of $f$ at each $t$, there are no issues whatsoever in justifying the computation below.} \begin{equation*}
\begin{split}
	\frac{\ud}{\ud t} \Re\brk{\tphi,\tphi^{app}}_f &= \Re \left(\int  {- 2 \sgm_{c} \abs{f}^{- 2 \sgm_{c} -1}} \rd_{t} \abs{f} \tphi \overline{\tphi^{app}} +  \int i|f|^{ {- 2 \sgm_{c}}} (\calL_f[\tphi]  {-} Q_f[\tphi])\overline{\tphi^{app}} - \int i|f|^{ {- 2 \sgm_{c}}}\tphi(\overline{\calL_f[\tphi^{app}]  {-} \err_{\tphi}} ) \right).
\end{split}
\end{equation*} Using the estimates for $|\rd_t|f||, Q_f[\tphi]$, and $\err_{\tphi}$, we can bound \begin{equation*}
\begin{split}
	\left|  \Re \int  {\abs{f}^{- 2 \sgm_{c} - 1}} \rd_{t} \abs{f} \tphi \overline{\tphi^{app}}  \right|  \lesssim \nrm{f}_{L^\infty_tC^{1,1}}^2 \nrm{\tphi}_{L^2_f} \nrm{\tphi^{app}}_{L^2_f}  ,
\end{split}
\end{equation*} \begin{equation*}
\begin{split}
	\left| \Re \int i|f|^{ {- 2 \sgm_{c}}}Q_f[\tphi] \overline{\tphi^{app}}  \right| \lesssim \left(\nrm{f}_{L^\infty_tC^{1,1}}^2+\nrm{\tphi}_{L^\infty_t C^{1,1}}^2\right) \nrm{\tphi}_{L^2_f}\nrm{\tphi^{app}}_{L^2_f}
\end{split}
\end{equation*} and \begin{equation*}
\begin{split}
	 \left| \Re\int i|f|^{ {- 2 \sgm_{c}}}\tphi \, \overline{ \err_{\tphi}} \right| \le C_{f,\dlt} \nrm{\tphi^{app}}_{L^2_f} A^{- \sgm_{c}+3} \nrm{g_{0}}_{H^2_{(x_1)}}. 
\end{split}
\end{equation*}
We now consider the remaining expression
\begin{equation*}
	\int \abs{f}^{- 2 \sgm_{c}} \Re(i \calL_{f}[\tphi] \br{\tphi^{app}}) \, \ud x
	- \int \abs{f}^{- 2 \sgm_{c}} \Re(i \tphi \br{\calL_{f}[\tphi^{app}]}) \, \ud x.
\end{equation*}
For the contribution of the principal term $\abs{f}^{2} \rd_{xx}$, we obtain
\begin{align*}
	& \int \abs{f}^{- 2 \sgm_{c}} \Re(i \abs{f}^{2} \rd_{xx} \tphi \br{\tphi^{app}}) \, \ud x
	- \int \abs{f}^{- 2 \sgm_{c}} \Re(i \tphi \br{\abs{f}^{2} \rd_{xx} \tphi^{app}}) \, \ud x \\
	& = - \int \abs{f}^{2 - 2 \sgm_{c}} \Re(i \rd_{x} \phi \br{\rd_{x} \tphi^{app}}) \, \ud x
	- \int (2 - 2 \sgm_{c}) \abs{f}^{1 - 2\sgm_{c}} \rd_{x} \abs{f} \Re(i \rd_{x} \phi \br{\tphi^{app}})  \, \ud x \\
	&\relphantom{=}
	+ \int \abs{f}^{2 - 2\sgm_{c}} \Re(i \rd_{x} \phi \br{\rd_{x} \tphi^{app}})  \, \ud x
	+ \int (2 - 2 \sgm_{c}) \abs{f}^{1 - 2 \sgm_{c}} \rd_{x} \abs{f} \Re(i \phi \br{\rd_{x} \tphi^{app}})  \, \ud x \\
	&= - \int (2 - 2 \sgm_{c}) \abs{f}^{1 - 2 \sgm_{c}} \rd_{x} \abs{f} \left( \Re(i \rd_{x} \phi \br{\tphi^{app}}) - \Re(i \phi \br{\rd_{x} \tphi^{app}})\right)  \, \ud x =: I.
\end{align*}
Since $\abs{f} \rd_{x} \abs{f} = \Re(\br{f} \rd_{x} f)$, this term cancels with some of the first-order terms, i.e., 
\begin{align*}
	\int (\alp_{1} + 2 \bt_{1}) \abs{f}^{- 2 \sgm_{c}} \Re(\br{f} \rd_{x} f) \Re(i \rd_{x} \tphi \br{\tphi^{app}}) \, \ud x
	- \int (\alp_{1} + 2 \bt_{1}) \abs{f}^{- 2 \sgm_{c}} \Re(\br{f} \rd_{x} f) \Re(i \tphi \br{\rd_{x} \tphi^{app}}) \, \ud x 
	= - I.
\end{align*}
For the remaining first-order terms, we have, after integrating by parts,
\begin{align*}
	&- \int (- \alp_{1} + 2 \bt_{1}) \abs{f}^{- 2 \sgm_{c}} \Im(\br{f} \rd_{x} f) \Re(\rd_{x} \tphi \br{\tphi^{app}}) \, \ud x
	-  \int (- \alp_{1} + 2 \bt_{1}) \abs{f}^{- 2 \sgm_{c}} \Im(\br{f} \rd_{x} f) \Re(\tphi \br{\rd_{x} \tphi^{app}}) \, \ud x \\
	&= (-\alp_{1} + 2 \bt_{1}) \int \left( \rd_{x} (\abs{f}^{- 2 \sgm_{c}} \Im(\br{f} \rd_{x} f)) \right) \Re(\tphi \br{\tphi^{app}}) \, \ud x, \\
	&\int \alp_{1} \abs{f}^{- 2 \sgm_{c}}  \Re(i f \rd_{x} f \rd_{x} \br{\tphi} \br{\tphi^{app}}) \, \ud x
	- \int \alp_{1} \abs{f}^{- 2 \sgm_{c}} \Re(i \phi  \br{f \rd_{x} f \rd_{x} \br{\tphi^{app}}}) \, \ud x \\
	&= \frac{\alp_{1}}{2}\int \abs{f}^{- 2 \sgm_{c}}  \left( i f \rd_{x} f \rd_{x} \br{\tphi} \br{\tphi^{app}}
	- i \br{f \rd_{x} f} \rd_{x} \tphi \tphi^{app}\right) \, \ud x 	- \frac{\alp_{1}}{2} \int \abs{f}^{- 2 \sgm_{c}} \left( i \br{f \rd_{x} f} \phi  \rd_{x} \tphi^{app}
	- i f \rd_{x} f \br{\phi} \rd_{x} \br{\tphi^{app}} \right) \, \ud x \\
	&= - \frac{\alp_{1}}{2}\int \left( \rd_{x} (i \abs{f}^{- 2 \sgm_{c}}  f \rd_{x} f) \right) \br{\tphi} \br{\tphi^{app}} \, \ud x
	+ \frac{\alp_{1}}{2}\int \left( \rd_{x} (i \abs{f}^{- 2 \sgm_{c}}  \br{f \rd_{x} f}) \right) \tphi \tphi^{app} \, \ud x.
\end{align*}
Both expressions may be bounded from the above by $C \nrm{f}_{C^{1, 1}}^{2} \nrm{\tphi}_{L^{2}_{f}} \nrm{\tphi^{app}}_{L^{2}_{f}}$. Finally, for the zeroth order terms, we easily have
\begin{align*}
	\left| \int \abs{f}^{- 2 \sgm_{c}} \Re(i (V_{f} \tphi + W_{f} \br{\tphi}) \br{\tphi^{app}}) \, \ud x
	- \int \abs{f}^{- 2 \sgm_{c}} \Re(i \tphi \br{(V_{f} \tphi^{app} + W_{f} \br{\tphi^{app}})}  \, \ud x \right|
	\aleq \nrm{f}_{C^{1, 1}}^{2} \nrm{\tphi}_{L^{2}_{f}} \nrm{\tphi^{app}}_{L^{2}_{f}}.
\end{align*} 
This gives \eqref{eq:DS-gei}, which concludes the proposition.	
\end{proof}

\subsection{Proof of Theorem~\ref{thm:ill-posed-unbounded-gDS}}\label{subsec:DS-unbounded}
We are now in a position to conclude the proof of Theorem~\ref{thm:ill-posed-unbounded-gDS} for the equation \eqref{eq:gDS}. To begin with, let $f$ satisfy the assumptions of the theorem with $f_0 = f(t=0)$. We may assume that $f_{0}(0)=0$ and $f_{0}'(0) =: A > 0$ by translation and phase rotation if necessary. We also fix $x_{1}$ as in Proposition \ref{prop:DS-WP}. 

Now let $\eps>0$, $  {s_{0}}\ge s_{c}$ and $0<\dlt' \le \dlt$ be given. We take some $\lmb\le -1$ and $g_{0}$ satisfying the assumptions of Proposition \ref{prop:DS-WP}, and define $\tphi_{0}$ by \begin{equation*}
\begin{split}
\tphi_{0} = \eps \, c( {s_{0}})|\lmb|^{- {s_{0}}} \tphi^{app}_{(\lmb)}(t=0)[g_0;f]. 
\end{split}
\end{equation*} Here, $\tphi^{app}_{(\lmb)} = \tphi^{app}_{(\lmb)}[g_0;f]$ is the degenerating wave packet constructed in Proposition \ref{prop:DS-WP} using $g_0$. It is not difficult to check  that $\tphi^{app}_{(\lmb)}(t=0) \in C^\infty_{c}$ and $\nrm{\tphi^{app}_{(\lmb)}(t=0)}_{C^{ {s_{0}}}}\lesssim_{ {s_{0}}} |\lmb|^{ {s_{0}}}$. Hence, by taking a sufficiently small $c( {s_{0}})>0$, we can ensure that $\nrm{\tphi_0}_{C^{ {s_{0}}}} \le \eps$ uniformly for all $\lmb\le-1$, as required by the statement of the theorem. We observe that \begin{equation}\label{eq:DS-inflation-initial-obs}
\begin{split}
\Re\brk{\tphi_0,\tphi^{app}_{(\lmb)}(t=0)}_f \ge c_0 \nrm{\tphi_0}_{L^2_f}\nrm{\tphi^{app}(t=0)}_{L^2_f}
\end{split}
\end{equation} for some $c_0>0$ independent of $\lmb$. To proceed, let us assume that the first option in the theorem does not hold; namely, there exists a solution $\phi$ to \eqref{eq:gDS} satisfying $\nrm{\phi-f}_{L^\infty([0,\dlt'];C^{ s_{c}})}<+\infty$ and $\phi(t=0)= f_0 + \tphi_0$. On $[0,\dlt']$, we write $\tphi = \phi-f$ and set \begin{equation*}
	\begin{split}
		M_{ {2}} = \sup_{t\in[0,\dlt']} \left( \nrm{f(t)}_{C^{1,1}} + \nrm{\tphi}_{C^{1,1}}\right).
	\end{split}
\end{equation*} We shall now establish the claimed norm inflation statement for $\tphi$, by taking $|\lmb|$ sufficiently large but in a way depending only on $f$ and $\dlt'$. 

\medskip

On the time interval $[0,\dlt']$, using Proposition \ref{prop:mee-Schro} and \eqref{eq:ds-wp-reg} we obtain that \begin{equation*}
\begin{split}
\nrm{\tphi(t)}_{L^2_f}\le \exp(CM_{ {2}}^2t)\nrm{\tphi_0}_{L^2_f}, \qquad 	\nrm{\tphi^{app}(t)}_{L^2_f} \le C_{f,\dlt} A^{- \sgm_{c}+1} \nrm{g_0}_{L^2}. 
\end{split}
\end{equation*}
In particular, we note that $\nrm{\tphi_0}_{L^2_f} < +\infty$ since $\tphi_{0}$ is supported away from the zeroes of $f$, and as discussed in \S\ref{subsec:modified-energy-estimate}, $\tphi(t)$ vanishes sufficiently fast at the zeroes of $f$ (ultimately due to Lemma~\ref{lem:gDS-vanishing-propagate}) so that $\nrm{\tphi(t)}_{L^2_f}$ is well-defined and obeys the above bound.
Applying \eqref{eq:DS-gei}, integrating in time on the interval $[0, \min\{ \dlt', cM_{2}^{-2}, A^{-2}|\lmb|^{-\frac12}  \}]$ for a sufficiently small $c>0$ and using \eqref{eq:DS-inflation-initial-obs}, we have \begin{equation}\label{eq:gDS-L2-lowerbound}
\begin{split}
\Re\brk{\tphi(t),\tphi^{app}(t)}_f \ge \frac{c_0}{2}\nrm{\tphi_0}_{L^2_f}\nrm{g_0}_{L^2} \quad \hbox{ for } \abs{t} \leq \min\{ \dlt', cM_{2}^{-2}, A^{-2}|\lmb|^{-\frac12} \}. 
\end{split}
\end{equation} Next, applying \eqref{eq:ds-wp-deg}--\eqref{eq:ds-wp-deg-small} with $\gmm' = - \sgm_{c}$ and $s = s_{c}$, we have
\begin{align*}
&\Re\brk{\tphi(t),\tphi^{app}(t)}_f - C_{f,\dlt} A^{- \sgm_{c}+1} |\lmb|^{-1} \nrm{\tphi_{0}}_{L^{2}_{f}}  \nrm{g_{0}}_{H^{\abs{s_{c}}  }_{(x_{1})} }  \\
& \leq \Re\brk*{\tphi(t), \rd_{x}^{s_{c}} \left(  \frac{\abs{f}^{- \sgm_{c}+s_{c}-\frac{1}{2}}}{i^{s_{c}} \lmb^{s_{c}} (1+\abs{f} \rd_{x} S)^{s_{c}}} \psi^{app}(t)\right)}  \\
& \leq \nrm{\rd_{x}^{s_{c}} \tphi(t)}_{L^{\infty}}  |\lmb|^{-s_{c}} \nrm{(1+\abs{f} \rd_{x} S)^{-s_{c}} \abs{f}^{s_{c}-\sgm_{c}-\frac{1}{2}} \psi^{app}(t)}_{L^{1}} \\
& \leq C_{f,\dlt}^{1-\sgm_{c}} A^{s_{c} -\sgm_{c} + \frac12 } |\lmb|^{-s_{c}} e^{- 2 \abs{\lmb} (s_{c} -\sgm_{c} + \frac{1}{2}) A^{2} t}\nrm{\rd_{x}^{s_{c}} \tphi(t)}_{L^{\infty}} \nrm{g_{0}}_{L^{2}}.
\end{align*}
Taking $\abs{\lmb}$ sufficiently large, we may ensure that  
\begin{align}\label{eq:lmb-req}
C_{f,\dlt} A^{-\sgm_{c}+1} |\lmb|^{-1} \nrm{g_{0}}_{H_{(x_{1})}^{\abs{s_{c}}}} \le \frac{c_{0}}{4} \nrm{g_{0}}_{L^{2}} \qquad \mbox{and} \qquad A^{-2}|\lmb|^{-\frac12} < \dlt', 
\end{align}
which gives, after combining the previous two inequalities with \eqref{eq:gDS-L2-lowerbound}, 
\begin{align*}
	\frac{c_{0}}{4}  C_{f,\dlt}^{-\sgm_{c}+1} A^{-(s_{c} - \sgm_{c} + \frac12 )} |\lmb|^{s_{c}} e^{2 \abs{\lmb}(s_{c} -\sgm_{c} + \frac{1}{2}) A^{2} t} \nrm{\tilde{\phi}_{0}}_{L^{2}_{f}} \leq \nrm{\rd_{x}^{s_{c}} \tphi(t)}_{L^{\infty}} \quad \hbox{ for } \abs{t} \leq \min\{ cM_{2}^{-2}, A^{-2}|\lmb|^{-\frac12} \}.
\end{align*} 
For each $|\lmb|$ satisfying \eqref{eq:lmb-req}, there are two cases; either (i) $cM_{2}^{-2} < A^{-2}|\lmb|^{-\frac12}$ or (ii) $cM_{2}^{-2} \ge A^{-2}|\lmb|^{-\frac12}$. In the case (i), we obtain that $M_{2} \gtrsim_{A} |\lmb|^{\frac14} \gtrsim_{A} (\dlt')^{-\frac12}$ using \eqref{eq:lmb-req}. Here, we could have assumed that $|\lmb|$ is sufficiently large from the beginning so that $\sup_{t \in [0,\dlt']} \nrm{f(t)}_{C^{1,1}} \ll_{A} |\lmb|^{\frac14}.$ Then, $M_{2} \simeq \sup_{t \in [0,\dlt']} \nrm{\tphi(t)}_{C^{1,1}} $ and the desired norm inflation follows simply from our assumption in \eqref{eq:gDS-sc} that $s_{c} \ge 2$. In the case (ii), we simply take  $t = A^{-2}\abs{\lmb}^{-\frac{1}{2}}$ in \eqref{eq:gDS-L2-lowerbound}, which gives the claimed norm inflation (actually, in this case we obtain a much stronger growth in terms of $1/\dlt'$) using that $s_{c} > \sgm_{c} - \frac{1}{2}$. This finishes the proof. 
\hfill \qedsymbol 

\begin{remark} \label{rem:sobolev-instab}
At the end of the above proof, observe that we could have followed the same argument but have used \eqref{eq:ds-wp-deg}--\eqref{eq:ds-wp-deg-small} with $\gmm' = - \sgm_{c}$, $s = \sgm$ and $p = 2$ to derive
\begin{equation*}
\frac{c_{0}}{4}  C_{f,\dlt}^{-\sgm_{c}+1} A^{-(\sgm - \sgm_{c})} |\lmb|^{s_{c}} e^{2 \abs{\lmb}(\sgm -\sgm_{c} ) A^{2} t} \nrm{\tilde{\phi}_{0}}_{L^{2}_{f}} \leq \nrm{\tphi(t)}_{H^{\sgm}},
\end{equation*}
for $t \leq \min\set{c M_{2}^{-2}, A^{-2} \abs{\lmb}^{-\frac{1}{2}}}$. This can be used to prove the inflation of the $H^{\sgm}$-norm for any $\sgm > \sgm_{c}$ in the second alternative of Theorem~\ref{thm:ill-posed-unbounded-gDS}.
\end{remark}

\subsection{Proof of Theorem \ref{thm:ill-posed-nonexist-gDS}}\label{subsec:DS-nonexist}

Let us divide the proof of Theorem \ref{thm:ill-posed-nonexist-gDS} into several steps. 

\medskip

\noindent \textit{Choice of background solution}. Towards a contradiction, we shall assume that  there exist $\eps>0$ and $s_{0} \ge s_{c} + 2$ such that for any $\phi_0\in C^\infty(\bbT)$ satisfying $\nrm{\phi_0}_{C^{s_{0}}}<\eps$, there exist $\dlt = \dlt(\phi_0)>0$ and a solution $\phi \in L^\infty([0,\dlt];{  C^{s_{c}+1,1}})$ to  {\eqref{eq:gDS}} with initial data $\phi(t=0)=\phi_0$.

Under this assumption, let us fix a function $\mathring{f}_0 \in C^\infty(\bbT)$ which is supported in $(-\frac12,\frac12)$ and $\mathring{f}_0(x) = x$ in $[-\frac14,\frac14]$. We then set \begin{equation*}
\begin{split}
f_0 := \sum_{k=k_{0}}^{\infty} f_{k,0} := \sum_{k=k_{0}}^{\infty}  A_k 2^{-k} \mathring{f}_0(2^k(x-x_k)) ,\quad A_k = 2^{-k^2}, \quad x_k = 2^{-\frac{k}{2}}
\end{split}
\end{equation*} where $k_{0}=k_{0}(s_{0},\eps,\mathring{f}_0) \ge 1$ is taken sufficiently large to achieve $\nrm{f_{0}}_{C^{s_{0}}}<\frac{\eps}{2}$. It is not difficult to see that $f_0 \in C^\infty(\bbT)$. Furthermore, since the supports of $f_{k,0}$ are disjoint from each other, for each $k\ge k_{0}$, we may choose $\chi_k \in C^\infty(\bbT)$ to be a cutoff function satisfying $\chi_{k} = 1$ on $\supp(f_{k,0})$ and $\chi_{k} = 0$ on $\supp(f_{k',0})$ for any $k'\ne k$.   From the contradiction hypothesis, we have a solution $f(t) \in L^\infty([0,\dlt];{C^{{s_{c}+1},1}})$ to \eqref{eq:gDS} with initial data ${f}_0$, for some $\dlt>0$. The estimate $|\rd_{t} {f}| \lesssim |{f}|$ shows that $\supp({f}(t)) = \supp({f}_0)$ on $[0,\dlt]$, and since $\chi_{k}$ equals either 0 or 1 on $\supp({f}(t)) = \supp({f}_0)$, we have that \begin{equation*}
\begin{split}
	\chi_k  = \chi_{k}^{3}, \qquad \rd_{x} \chi_{k} = 0 
\end{split}
\end{equation*} on $\supp({f}(t))$ for any $k\ge k_{0}$. Using these observations, it follows that for each $k\ge k_{0}$, we have that $f_{k} := \chi_{k} f $ is again a solution to \eqref{eq:gDS} with initial data $f_{k}(t=0) = f_{k,0}$. Furthermore, the $L^\infty([0,\dlt];{C^{{s_{c}+1},1}})$--norm of $f_{k}$ is bounded uniformly in $k$. 

In the following, we consider the weight $w(t,x) =  {|f(t,x)|^{- 2 \sgm_{c}}}$ and define \begin{equation*}
\begin{split}
\brk{a,b}_w (t) = \int_{\bbT} |w(t,x)| a(t,x)\overline{b(t,x)} \,\ud x, \quad \nrm{a}_{L^{2}_{w}}^{2} = \brk{a,a}_w. 
\end{split}
\end{equation*} 


\medskip

\noindent \textit{Choice of wave packet solutions}. We now fix some nonzero function $g_0 \in C^\infty$ supported in $(\frac{1}{8}, \frac{1}{4})$  and take $g_{k}(x) :=  {2^{\frac{k}{2}}} g_0(2^{k}(x-x_{k}))$. For some  sequence $\{\lmb_k\}_{k\ge k_{0}}$ to be determined (for now, we take $-\lmb_{k} \ge A_{k}^{10}$), we consider the sequence of wave packet solutions 
\begin{equation*}
	\begin{split}
		 {\tphi_{k}^{app} := \tphi^{app}_{(\lmb_{k})}[g_{k}; f_{k}], }
	\end{split}
\end{equation*} where $\tphi^{app}_{(\lmb_{k})}[g_{k}; f_{k}]$ is the wave packet solution from Proposition \ref{prop:DS-WP} with data $g_{k}, \lmb_{k}$, adapted to the linearly degenerate solution $f_{k}$, with $A = A_{k}$ and $x_{1} = 2^{-k-2}$.  We define the corresponding error by \begin{equation}\label{eq:DS-nonexist-error}
\begin{split}
[i\rd_t + \calL_{f_{k}}] \tphi^{app}_k = {\err}_k 
\end{split}
\end{equation}  {where the operator $[i\rd_t + \calL_{f_{k}}]$ is obtained from \eqref{eq:deg-Schro-pert} by replacing $f$ with $f_{k}$.} Applying Proposition \ref{prop:DS-WP}, we obtain the following bounds: with $\dlt_{k} := \min\{  \dlt, A_{k}^{-2}|\lmb_{k}|^{-\frac12}  \} = A_{k}^{-2}|\lmb_{k}|^{-\frac12}  $ (by our choice of $-\lmb_{k}$ in the above), \begin{itemize}
	\item $\nrm{\tphi^{app}_k}_{L^\infty([0,\dlt_{k}];L^2_w)} \le C_{f_{k},\dlt_{k}}\nrm{\tphi^{app}_k}_{L^2_w}(t=0) \le C_{f_{k},\dlt_{k}} A_{k}^{-\sgm_{c}+1} \nrm{g_{0}}_{L^{2}} $;
	\item $\nrm{\err_k}_{L^\infty([0,\dlt_{k}],L^2_{w})} \le C_{f_{k},\dlt_{k}} A_{k}^{-\sgm_{c}+3} \nrm{g_{k}}_{H^{2}_{(2^{-k-2})}} \le C_{f_{k},\dlt_{k}} A_{k}^{- \sgm_{c}+3} \nrm{g_{0}}_{H^{2}} $; 
\end{itemize} and \begin{equation*}
\begin{split}
	|f|^{-2\sgm_{c}} \tphi^{app}_{k} = \rd_{x}^{s_{c}} \left(  \frac{|f|^{-\sgm_{c} + s_{c} - \frac12}}{ (i \lmb_{k})^{s_{c}} (1 + |f| \rd_{x}S)^{s_{c}} } \psi^{app}_{k}  \right) + |f|^{-2\sgm_{c}} \tphi^{small}_{k}
\end{split}
\end{equation*} 
with \begin{equation*}
\begin{split}
	\left\Vert \frac{|f|^{-\sgm_{c} + s_{c} - \frac12}}{ (i \lmb_{k})^{s_{c}} (1 + |f| \rd_{x}S)^{s_{c}} } \psi^{app}_{k}   \right\Vert_{L^{1}} \le C_{f_{k},\dlt_{k}}^{1-\sgm_{c}} A_{k}^{-\sgm_{c}+ s_{c}+\frac12}  |\lmb_{k}|^{-s_{c}} \exp\left( -2|\lmb_{k}|(-\sgm_{c} + s_{c} + \frac12) A_{k} t \right)  \nrm{g_0}_{L^2}
\end{split}
\end{equation*} for $ 0 \le t \le \dlt_{k}$ and \begin{equation*}
\begin{split}
	\nrm{  \tphi_{k}^{small} }_{L^\infty([0,\dlt_{k}]; L^{2}_{w})} \le C_{f_{k},\dlt_{k}} A_{k}^{-\sgm_{c}+1}|\lmb_{k}|^{-1} \nrm{g_{k}}_{H^{s_{c}}_{ (2^{-k-2}) } } \le C_{f_{k},\dlt_{k}} A_{k}^{-\sgm_{c}+1} |\lmb_{k}|^{-1}   \nrm{g_{0}}_{H^{s_{c}}}. 
\end{split}
\end{equation*} 
From \eqref{eq:Cf-bound} we see that \begin{equation*}
\begin{split}
	 C_{f_{k},\dlt_{k}} \lesssim (1 + A_{k}^{-1}M)^{N_{0}} \exp(C_{0} M^{2} \dlt_{k} ), \qquad M = \sup_{t \in [0, \dlt]} \nrm{f(t,\cdot)}_{C^{s_{c}+1,1}}
\end{split}
\end{equation*} where the implicit constant and $N_{0}$ depends  on $g_{0}, \alp_{1}, \bt_{1}, \mu_{1}, s_{c}$, but not on $k$ and $\lmb_{k}$. Then, simply using $\dlt_{k}\le \dlt$ and recalling $A_{k} = 2^{-k^{2}}$, we see that $
	C_{f_{k},\dlt_{k}} \lesssim 2^{N_{0} k^{2}}$ holds, where the implicit constant depends further on $M, \dlt$ but not on $k$ and $\lmb_{k}$. In turn, this gives an upper bound on the constants in the estimates above; for instance \begin{equation*}
		\begin{split}
			C_{f_{k},\dlt_{k}}^{1-\sgm_{c}} A_{k}^{-\sgm_{c}+ s_{c}+\frac12} \lesssim 2^{N_{1} k^{2}}
		\end{split}
	\end{equation*} with some $N_{1} > 0 $ depending additionally on $\sgm_{c}, s_{c}$. \textbf{In the following, we shall write $\lesssim$ as long as the implicit constant does not depend on $k$ and $\lmb_{k}$.}

\medskip

\noindent \textit{Choice of initial data.} We now take  \begin{equation}\label{eq:gDS-initial-data}
\begin{split}
\tphi_0(x) = \sum_{k=k_{0}}^{\infty} \tphi_{k,0}(x) := \sum_{k=k_{0}}^{\infty} \exp(-|\lmb_k|^{\frac14}) \tphi_k^{app}(t=0,x), 
\end{split}
\end{equation} which belongs to $C^\infty(\bbT)$. By taking $k_{0}$ even larger if necessary, we can guarantee that $\nrm{\tphi_0}_{C^{m}}<\frac{\eps}{2}$. Then we set $\phi_0 = f_0 + \tphi_0,$ which verifies $\nrm{\phi_{0}}_{C^{m}}<\eps$. Again from the contradiction hypothesis, we have a $L^\infty_tC^{ {s_{c}+1},1}$--solution $\phi(t)$ to \eqref{eq:gDS} with initial data $\phi_0$ on some time interval $[0,\dlt']$. We may assume that $0<\dlt'\le\dlt$, and define \begin{equation*}
\begin{split}
\tphi(t) := \phi(t) - f(t), \quad \tphi_k(t) := \chi_k\tphi(t)
\end{split}
\end{equation*} for all $k\ge k_{0}$. We have that $\sum_{k=k_{0}}^{\infty} \tphi_k = \tphi$; this follows from $\rd_t |f+\tphi| \lesssim |f+\tphi|$ and the uniform pointwise estimate $
|f+\tphi|(t,x) \lesssim |f_0+\tphi_0|(x)\lesssim |f_0|(x).$ Then we see that $\tphi_k$ solves \begin{equation*}
\begin{split}
[i\rd_t+\calL_{f_{k}}]\tphi_k = Q_{f_k}[\tphi_k]  , 
\end{split}
\end{equation*} (which is \eqref{eq:plugin2} with $f$ and $\tphi$ replaced with $f_k$ and $\tphi_k$, respectively). We note that the $L^\infty([0,\dlt']; C^{s_0+1,1})$ norm is uniformly bounded for $\{ f_{k} \}_{k\ge k_0}$ and $\{ \tphi_{k} \}_{k \ge k_{0}}$. Therefore, from Proposition \ref{prop:mee-Schro}, we obtain the estimate   \begin{equation}\label{eq:gee0}
\begin{split}
\nrm{\tphi_k}_{L^\infty([0,\dlt'];L^2_{w})} \lesssim \nrm{\tphi_{k,0}}_{L^2_{w_0}}
\end{split}
\end{equation} uniformly in $k\ge k_{0}$. Now, combining this with the generalized energy estimate \eqref{eq:DS-gei} for $\tphi_{k}$ and $\tphi^{app}_{k}$, we obtain that \begin{equation}\label{eq:gee}
\begin{split}
\left|\frac{\ud}{\ud t} \, \Re \brk{\tphi_k, \tphi^{app}_k}_{w}  \right| \lesssim    2^{N_{1} k^{2}} \nrm{\tphi_{k,0}}_{L^2_{w_0}} \nrm{g_{0}}_{H^{2}} 
\end{split}
\end{equation} for $t \in [0,\dlt_{k}]$. We shall now take $|\lmb_{k}|$ larger so that $\dlt_{k} = A_{k}^{-2}|\lmb_{k}|^{-\frac12}$ satisfies $2^{N_{1} k^{2}}\dlt_{k}$ is very small with respect to the implicit constants in \eqref{eq:gee0} and \eqref{eq:gee}. Then, since at $t = 0$ we have \begin{equation*}
\begin{split}
\Re \brk{\tphi_k, \tphi^{app}_k}_{w} (t = 0)  \ge \frac14 \nrm{\tphi_{k,0}}_{L^2_{w_0}}\nrm{\tphi^{app}_{k,0}}_{L^2_{w_0}},
\end{split}
\end{equation*} by integrating \eqref{eq:gee} in time from $t = 0$ to $\dlt_{k}$, we obtain \begin{equation}\label{eq:gee2}
\begin{split}
\Re \brk{\tphi_k, \tphi^{app}_k}_{w} (\dlt_{k}) \ge \frac{1}{8} \nrm{\tphi_{k,0}}_{L^2_{w_0}}\nrm{\tphi^{app}_{k,0}}_{L^2_{w_0}} . 
\end{split}
\end{equation} At $t = \dlt_{k}$ we write 
\begin{equation*}
	\begin{split}
		\Re \brk{\tphi_k, \tphi^{app}_k}_{w}(\dlt_{k})  = (-1)^{s_{c}}\Re \brk{ \rd_{x}^{s_{c}}\tphi_{k},  \frac{|f|^{-\sgm_{c} + s_{c} - \frac12}}{ (i \lmb_{k})^{s_{c}} (1 + |f| \rd_{x}S)^{s_{c}} } \psi^{app}_{k} }  + \Re \brk{ \tphi_{k}, \tphi^{small}_{k}}_{w} 
	\end{split}
\end{equation*} and then combining the estimates of the right hand side with \eqref{eq:gee2}, we get \begin{equation*}
\begin{split}
	\nrm{\tphi_{k,0}}_{L^2_{w_0}}\nrm{\tphi^{app}_{k,0}}_{L^2_{w_0}} \lesssim 2^{N_{1} k^{2}} \left( \nrm{\rd_{x}^{s_{c}} \tphi_{k} }_{L^\infty}   |\lmb_{k}|^{-s_{c}} \exp\left( -2|\lmb_{k}|^{\frac12}(-\sgm_{c} + s_{c} + \frac12) \right)  \nrm{g_0}_{L^2}  + |\lmb_{k}|^{-1} \nrm{g_{0}}_{H^{s_{c}}}  \right) .
\end{split}
\end{equation*}
By taking $|\lmb_{k}|$ even larger if necessary, we can guarantee that  $|\lmb_{k}|^{-1} 2^{ N_{1} k^{2} } \nrm{g_{0}}_{H^{s_{c}}} \ll \nrm{\tphi^{app}_{k,0}}_{L^2_{w_0}}$ holds ($\ll$ is defined in terms of the implicit constant in the previous inequality), so that we deduce \begin{equation*}
\begin{split}
	  \nrm{\tphi_{k,0}}_{L^2_{w_0}}\nrm{\tphi^{app}_{k,0}}_{L^2_{w_0}}  \lesssim 2^{N_{1} k^{2}} \nrm{\rd_{x}^{s_{c}} \tphi_{k} }_{L^\infty}|\lmb_{k}|^{-s_{c}} \exp\left( -2|\lmb_{k}|^{\frac12}(-\sgm_{c} + s_{c} + \frac12) \right) \nrm{g_{0}}_{L^{2}} 
\end{split}
\end{equation*} and then recalling the form of $\tphi_{k,0}$ from \eqref{eq:gDS-initial-data},  \begin{equation*}
\begin{split}
	\nrm{\rd_{x}^{s_{c}} \tphi_{k} }_{L^\infty}( t= \dlt_{k} ) \gtrsim 2^{-N_{1} k^{2}} |\lmb_{k}|^{s_{c}} \exp\left( 2|\lmb_{k}|^{\frac12}(-\sgm_{c} + s_{c} + \frac12) - |\lmb_{k}|^{\frac14} \right) .
\end{split}
\end{equation*} Note that the left hand side is bounded by $\nrm{f}_{L^\infty([0,\dlt'];C^{s_{c}+1,1})} + \nrm{\phi}_{L^\infty([0,\dlt'];C^{s_{c}+1,1})} $ for all $k$ sufficiently large. This is a contradiction since the right hand side diverges as $k\to\infty$. The proof is now complete. \hfill \qedsymbol

\section{KdV-type equations}\label{sec:KdV}
This section is organized as follows. After setting up some pieces of notation in Section~\ref{subsec:KdV-prelim}, we study the properties of regular cubically degenerate solutions -- typically denoted by $f$ -- in Section~\ref{subsec:KdV-bg}. Then in Section~\ref{subsec:KdV:deg-wp}, we carry out the key construction of degenerating wave packets for the linearized equation around $f$, and in Section~\ref{subsec:gdkdv-en}, we establish a modified energy estimate for the perturbation (solving the nonlinear difference equation) around $f$. Finally, in Sections~\ref{subsec:unbounded-gdkdv} and \ref{subsec:nonexist-gdkdv}, we prove Theorems~\ref{thm:ill-posed-unbounded-gdkdv} and \ref{thm:ill-posed-nonexist-gdkdv}, respectively.

 \subsection{Preliminaries} \label{subsec:KdV-prelim}
  
We introduce the following quantity defined for a $C^{2, 1}$ function $f$  {on an interval $I$}:  \begin{equation*}
\begin{split}
{\nrm{f}_{Y(I)} = \nrm{f^{-\frac{2}{3}}\rd_x f}^3_{L^\infty(I)} + \nrm{f^{-\frac{1}{3}}\rd_{xx}f}_{L^\infty(I)}^{\frac{3}{2}} + \nrm{f}_{L^\infty(I)} + \nrm{\rd_{xxx}f}_{L^\infty(I)} .}
\end{split}
\end{equation*} We shall write $f \in {Y(I)}$ if $\nrm{f}_{Y(I)}$ is finite. This quantity is appropriate to handle solutions with degeneracies of order  {at least} 3. For convenience we set \begin{equation*}
\begin{split}
 {\nrm{f}_{\widetilde{C}^{k,\alpha}(I)} = \nrm{f}_{C^{k,\alpha}(I)} + \nrm{f}_{Y(I)}. }
\end{split}
\end{equation*} For $f$ depending on time, we say $f \in L^\infty([0,\dlt];\widetilde{C}^{k,\alpha}(I))$ {(resp.~$f \in L^\infty([0,\dlt];Y(I)$)} if  \begin{equation*}
\begin{split}
\nrm{f}_{L^\infty([0,\dlt];\widetilde{C}^{k,\alpha}(I))} := \sup_{t\in[0,\dlt]} \nrm{f(t)}_{\widetilde{C}^{k,\alpha}(I)}  < + \infty. \\ {(\hbox{resp.~} \nrm{f}_{L^\infty([0,\dlt];Y(I))} := \sup_{t\in[0,\dlt]} \nrm{f(t)}_{Y(I)}  < + \infty. )}
\end{split}
\end{equation*} It is easy to see using the Taylor expansion, that any $C^{3,\alpha}$ function which vanishes cubically at its zeroes must belong to $Y$. However, \textit{propagation} of $Y$--boundedness for \eqref{eq:gdkdv} in general requires higher regularity, e.g. $C^{4,1}$ (see Proposition \ref{prop:C41}). 
 
For later use, we introduce the notation 
\begin{equation}\label{eq:f-gdkdv}
\begin{split}
\brk{a,b}_f (t) = \int f(t,x)^{- \frac{2}{3} \sgm_{c}} a(t,x) b(t,x)\,\ud x, \quad \nrm{a}_{L^2_f} (t) = \sqrt{\brk{a,a}_f (t)}.
\end{split}
\end{equation}  
For the motivation behind the power $f^{- \frac{2}{3} \sgm_{c}}$, see Section~\ref{subsec:gdkdv-en}.

\subsection{Properties of a {regular cubically degenerate} solution} \label{subsec:KdV-bg} 

We first discuss a few basic properties of a regular cubically degenerate solution $f$ to \eqref{eq:gdkdv}, which shall serve as the background for our illposedness mechanism.

Under the assumption $f\in L^\infty_t\widetilde{C}^{3,\alpha}(I)$ with any $\alp > 0$, we can propagate the information that $f$ vanishes cubically on {an endpoint of $I$} and compute the coefficient. 

\begin{lemma}\label{lem:KdV-C3-solution}
	Let $f \in L^\infty([0,\dlt]; \widetilde{C}^{3,\alpha}(I))$ be a solution of \eqref{eq:gdkdv} with initial data $f_0$ {that is positive on $I \setminus \rd I$ and vanishes to order at least $3$ on each point in $\rd I$}, where $0 < \alp \leq 1$. Then the following statements hold.
	\begin{enumerate}
\item The zeroes and the sign of $f(t, x)$ are preserved in time, i.e., $f(t, x)$ vanishes on $\rd I$ and $f(t, x) > 0$ for $x \in I \setminus \rd I$ for all $t \in [0, \dlt]$.
\item Let $I=[a,b]$. Then, the set of $t$-dependent functions $\set{\rd_{x}^{k} f(t, a)}_{k=0}^{3}$ for $t \in [0, \dlt]$ is determined by the initial data at $x = a$, i.e., $\set{\rd_{x}^{k} f(0, a)}_{k=0}^{3}$. In particular,
\begin{equation}\label{eq:KdV-cubic-profile}
	\begin{split}
	f(t,x) = (\beta(t)(x-a))^3 + O(\nrm{f}_{L^{\infty}_{t} ([0, \dlt]; C^{3, \alp}(I))} |x-a|^{3+\alp}) ,\quad x\rightarrow a^+ 
	\end{split}
	\end{equation} where $\beta(t)$ is the solution of \begin{equation}\label{eq:KdV-beta-ODE}
	\begin{split}
	\dot{\beta}(t) = - {(2+6 \alp_{1})} \beta^4(t), \quad 6\beta^3(0) = f_{0,xxx}(a),
	\end{split}
	\end{equation}
and the implicit constant in $O(\cdot)$ is universal. The same statement applies to $b \in \rd I$. 
\end{enumerate}
\end{lemma}
\begin{proof}
	Since we are assuming that $f(t,\cdot) \in C^3$, from \eqref{eq:gdkdv}, we have \begin{equation*}
	\begin{split}
	|\rd_tf| \le C  {(|f_{xxx}||f| + |f_x||f_{xx}| + \abs{f_{x}} \abs{f}^{m-1})}. 
	\end{split}
	\end{equation*} From the assumption that $f(t,\cdot)\in Y$, we have the pointwise estimate \begin{equation*}
	\begin{split}
	|\rd_tf|\le C \left( \nrm{f}_{Y} + \nrm{f}_{C^3} {(1+ \nrm{f}_{L^{\infty}}^{m-2})} \right)|f|.
	\end{split}
	\end{equation*} This shows that \begin{equation}\label{eq:KdV-f-pointwise}
	\begin{split}
	f_0(x) \exp\left( -Ct\nrm{f}_{L^\infty_t\tilde{C}^{3,\alpha}}  {(1+ \nrm{f}_{L^{\infty}}^{m-2})}  \right) \le f(t,x) \le f_0(x) \exp\left( Ct\nrm{f}_{L^\infty_t\tilde{C}^{3,\alpha} }  {(1+ \nrm{f}_{L^{\infty}}^{m-2})}  \right) 
	\end{split}
	\end{equation} for any $x\in I$, which proves the first statement. The second statement follows from simply evaluating the equation \eqref{eq:gdkdv} at $x = a, b$ and carrying out a minor modification of the proof of Lemma~\ref{lem:gDS-vanishing-propagate}. Here, the fact that $f$ vanishes  {at least} cubically at $x = a$ ensures that no $\rd_{x}^{k} f(t, a)$ with $k > 3$ occurs in the ODEs for the Taylor coefficients.  We omit the details. \qedhere
\end{proof}

Before we proceed further, let us note that the assumptions of  {Theorem}~\ref{thm:ill-posed-unbounded-gdkdv} on the solution $f$ are automatically satisfied for any sufficiently smooth solutions of \eqref{eq:gdkdv}. 

\begin{proposition}\label{prop:C41}
	 {Consider an interval $I = [a, b] \subseteq \bbT$. Let $f_0 \in C^{4,1}(\bbT)$ satisfy $f_0>0$ on $I \setminus \set{a}$ (resp.~$I \setminus \set{b}$) and vanishes at least cubically at $a$ (resp.~$b$), so that $f_{0} \in Y(I)$.  Then there exists $\dlt > 0$ depending on $\nrm{f_{0}}_{Y(I)}$ such that, if $f$ is a solution to \eqref{eq:gdkdv} with initial data $f_{0}$ satisfying $f \in L^\infty([0,\dlt];C^{4,1}(\bbT))$,} then $f$ satisfies $f|_{I} \in L^\infty([0,\dlt];  {Y(I)})$  {with the bound $\nrm{f}_{L^{\infty}([0, \dlt]; Y(I))} \leq C C_{0} \exp(C M \dlt)$.}
	 
	 Furthermore, for this value of $\dlt$, let $u$ be another solution to \eqref{eq:gdkdv} belonging to $L^\infty([0,\dlt];C^{4,1}(\bbT))$ with initial data $u_{0}$ satisfying $u_{0} \in Y(I)$ and \begin{equation}\label{eq:phi-f-ratio-initial}
	 	\begin{split}
	 		|u_{0}(x)|  \le C_{1}f_{0}(x)
	 	\end{split}
	 \end{equation} for some $C_{1}>0$ uniformly for $x\in I$. Then, for some $0<\dlt'\le\dlt$ depending only on $\nrm{u_0}_{Y}, \nrm{f_0}_{Y}$, \begin{equation}\label{eq:phi-f-ratio}
	 \begin{split}
	 	{|u(t,x)|} \le C_{1}(1+CC_{0}t) \exp(CMt) f(t,x), \qquad t \in [0,\dlt']
	 \end{split}
 \end{equation} uniformly for $x \in I$, where $C_{0} = C_{0}(\nrm{f_0}_Y,\nrm{u_0}_{Y})$ and $M=M(\nrm{f}_{L^\infty_tC^{4,1}},\nrm{u}_{L^\infty_tC^{4,1}}  )$.
\end{proposition}
\begin{proof}  Without loss of generality, we consider the case $f_{0} > 0$ on $I \setminus \set{a}$ with $f_{0}$ vanishing at least cubically at $a$.  We compute that \begin{equation*}
	\begin{split}
	\rd_t f & = -\mu_{1} f^{m-1} f_{x} - \alp_{1} f_xf_{xx} - f f_{xxx},\\
	\rd_tf_x & = - \mu_{1} (f^{m-1} f_{x})_{x} - \alp_{1} (f_{xx})^2 - (\alp_{1}+1) f_{x}f_{xxx} - f f_{xxxx},\\
	\rd_t f_{xx} & = - \mu_{1} (f^{m-1} f_{x})_{xx} - (3 \alp_{1} +1) f_{xx}f_{xxx} - (\alp_{1}+2) f_{x}f_{xxxx} - f f_{xxxxx}. 
	\end{split}
	\end{equation*} Upon $f \in L^\infty_t C^{4,1}$, we have the pointwise estimate \begin{equation*}
	\begin{split}
	\frac{\ud}{\ud t} \left(|f|^2 + |f_x|^3 + |f_{xx}|^6\right)  \le CM\left(|f|^2 + |f_x|^3 + |f_{xx}|^6 \right),
	\end{split}
	\end{equation*}  where we introduce the shorthand $M = 1 + \nrm{f}_{L^\infty_t C^{4,1}}^{4} +  \nrm{f}_{L^{\infty}_{t} C^{2,1}}^{ {6m-2}} $ for simplicity. By Gronwall's inequality,  we have the pointwise estimate  \begin{equation}\label{eq:Km2-cubic-pointwise}
	\begin{split}
	(|f|^2 + |f_x|^3 + |f_{xx}|^6)(t,x) \le \exp(C  {M} t)(|f_0|^2 + |f_{0,x}|^3 + |f_{0,xx}|^6)(x). 
	\end{split}
	\end{equation} From the assumptions on the initial data, we have that \begin{equation}\label{eq:Km2-initial-condition}
	\begin{split}
	|f_{0,x}(x)|^3 + |f_{0,xx}(x)|^6  \le C_{0}^{2} (f_{0}(x))^2
	\end{split}
	\end{equation} holds pointwise on $I$, with some $C_{0}>0$ depending  {only on $\nrm{f_{0}}_{Y}$}.  Returning to \eqref{eq:Km2-cubic-pointwise} and applying Young's inequality, we deduce the pointwise bound $$|f_x(t,x)f_{xx}(t,x) | \le C C_{0}\exp(C  {M} t) f_0(x)$$ for all $x \in I$. In turn, using this bound in the equation for $\rd_tf$, we obtain  {for all $x \in I$ that} \begin{equation}\label{eq:rdtf}
	\begin{split}
	|\rd_t f(t,x)| \le C  {M} |f(t,x)|+ C C_{0}\exp(C  {M} t)f_{0}(x).
	\end{split}
	\end{equation} 
 Dividing by $f_{0}$ and applying Gronwall's inequality to $\frac{f}{f_{0}} - 1$, we obtain
	\begin{equation*}
	\abs*{\frac{f}{f_{0}} - 1} \leq C C_{0} t \exp(C M t)
\end{equation*}
which, after some simplification, implies   
	\begin{equation}\label{eq:f-f_0-pointwise}
	\begin{split}
	 f(t,x) \le  {(1 + C C_{0} t) \exp(C M t)} f_0(x)
	\end{split}
	\end{equation} as well as \begin{equation}\label{eq:f-f_0-pointwise2}
	\begin{split}
	 f(t,x) \ge  {(1 - C C_{0} t) \exp(- C M t)} f_0(x).
	\end{split}
\end{equation}This guarantees that  {$f \in L^{\infty}_{t}([0, \dlt]; Y)$, provided that $\dlt$ is sufficiently small depending on $C_{0} = C_{0}(\nrm{f_{0}}_{Y})$}.

For the second statement, we note that the assumption \eqref{eq:phi-f-ratio-initial} implies \begin{equation*}
	\begin{split}
		|u_{0}(x)|^2 + |u_{0,x}(x)|^{3} + |u_{0,xx}(x) |^{6} \le C_{1}^{2}(1+ C\nrm{u_{0}}_{Y})|f_{0}(x)|^{2}
	\end{split}
\end{equation*} for some absolute constant $C>0$. With this bound, we may apply the above argument to $u$ instead of $f$ and obtain the bound \begin{equation*}
	\begin{split}
		|u(t,x)| \le  {(1 + C C_{0} t) \exp(C M t)} |u_0(x)| \le C_{1}(1+C_{0}t) \exp(C M t) f_{0}(x), \qquad t \in [0,\dlt'],
	\end{split}
\end{equation*} for some $0<\dlt'\le \dlt$ depending only on $C_{0} = C_{0}(\nrm{u_{0}}_{Y})$. Here, $M = 1 + \nrm{u}_{L^\infty_t C^{4,1}}^{4} +  \nrm{u}_{L^{\infty}_{t} C^{2,1}}^{ {6m-2}}$. Using this bound together with \eqref{eq:f-f_0-pointwise2}, we obtain the desired estimate \eqref{eq:phi-f-ratio}, by taking $\dlt'$ smaller in a way depending only on $ \nrm{u_{0}}_{Y} , \nrm{f_0}_{Y}$ if necessary. This finishes the proof.
\end{proof}  

\subsection{Degenerating wave packets for the linearized equation}\label{subsec:KdV:deg-wp}

In this subsection, our goal is to construct degenerating wave packets for the linearization of \eqref{eq:gdkdv} around a (possibly hypothetical) regular cubic degenerate solution; see Proposition~\ref{prop:KdV-deg-wp} below.

\subsubsection{Linearized equation and degenerating wave packets}

In the following, we fix some function $f$ that satisfies all the assumptions from  {Theorem}~\ref{thm:ill-posed-unbounded-gdkdv}, and further assume for simplicity that the interval is given by $I = [0,b]$ for some $b>0$.  {We} fix some $0< x_1<b$ that 
\begin{equation} \label{eq:gdkdv-x1}
\frac{1}{2}f_{0,xxx}(x) < f_{0,xxx}(0) < 2f_{0,xxx}(x) \hbox{ for all } x \in [0,x_1]. 
\end{equation}

We now write $u = f +\phi$, where $u$ is a solution to \eqref{eq:gdkdv}. Then, we have that $\phi$ must solve \begin{equation}\label{eq:KdV-pert}
\begin{split}
\rd_t\phi + \calL_f  {\phi}  = Q[\phi],
\end{split}
\end{equation} 
with 
\begin{equation} \label{eq:gdkdv-Lf}
	\calL_{f} \phi = f \phi_{xxx} + \alp_{1} f_{x} \phi_{xx} + (\alp_{1} f_{xx} + \mu_{1} f^{m-1}) \phi_{x} + (f_{xxx} + (m-1) \mu_{1} f^{m-2} f_{x}) \phi,
\end{equation}
\begin{equation} \label{eq:gdkdv-Q}
	Q[\phi] = - \phi \phi_{xxx} - \alp_{1} \phi_{x} \phi_{xx} -  {\frac{\mu_{1}}{m}} \left(  (f + \phi)^{m} - f^{m} - m f^{m-1} \phi \right)_{x}.
\end{equation}
We are concerned with constructing wave packets to the following linearized equation: \begin{equation}\label{eq:KdV-pert-lin-re}
\begin{split}
\rd_t\phi + \calL_f  {\phi} = 0. 
\end{split}
\end{equation} Recall the notation $\nrm{g}_{W^{s,p}_{(L)}} = \sum_{j = 0}^{s} \nrm{(L \rd_{x})^{j} g}_{L^{p}_{x}}$, $H^{s}_{(L)} = W^{s,2}_{(L)}$ from the previous section. Our aim is to prove the following result.
\begin{proposition}\label{prop:KdV-deg-wp}
	 Let $f \in L^{\infty}_{t}([0, \dlt]; \tld{C}^{s-1, 1}(I))$ be a solution to \eqref{eq:gdkdv} with initial data $f_{0}$ satisfying $f_{0} > 0$ on $I \setminus \set{0}$, vanishing cubically at $0$ and $f_{0} \in \tld{C}^{s_{0}-1, 1}(I)$, where $4 \leq s \leq s_{0}$. Let $A = \frac{1}{6} f_{0, xxx}(0)$ and fix $0 < x_{1} < 1$ so that \eqref{eq:gdkdv-x1} holds. 
	Then given $\lmb\in\bbN$ and $g_0 \in C^\infty_c$ supported in $(\frac{x_1}{2},x_1)$, we may associate a function  {$\phi^{app}_{(\lmb)}[g_{0}, f]$} defined in $[0, \dlt] \times I$ satisfying the following properties: \begin{itemize}
		\item (linearity) the map $g_{0} \mapsto \phi^{app}_{(\lmb)}[g_{0}, f]$ is linear;
		\item (support property) $\supp (\phi^{app}_{(\lmb)} [g_{0}, f] (t, \cdot)) \subset  {(0, x_{1}) \cap (0, C_{\tld{f}} \,  x_{1} e^{- 3\bt(t)A^{\frac23} \lmb^{2} t} } )$;
		\item  {(initial data estimates) for $0 \leq n \leq s_{0}$ and $1 \leq p \leq \infty$, we have}
	 		\begin{align}
	\frac{1}{C_{\tld{f}_{0}}} \left( \nrm{g_{0}}_{L^{2}} - \lmb^{-1} \nrm{g_{0}}_{H^{1}_{(x_{1})}} \right)  \leq \nrm{f_{0}^{-\frac{\sgm_{c}}{3}}(x) \phi^{app}_{\lmb}(0, x)}_{L^{2}} &\leq C_{\tld{f}_{0}} \nrm{g_{0}}_{L^{2}}, \label{eq:KdV-deg-wp-id-L2} \\
		\nrm{f_{0}^{-\frac{\sgm_{c}}{3}}(x)(A^{-\frac{1}{3}} f_{0}^{\frac{1}{3}}(x)\rd_x)^n \phi^{app}_{(\lmb)}(0,x)}_{L^p} &\leq C_{\tld{f}_{0}} x_{1}^{\frac{1}{p} - \frac{1}{2}}\abs{\lmb}^{n} \nrm{g_{0}}_{W^{n, p}_{(x_{1})}}; \label{eq:KdV-deg-wp-id}
		-nrm
		\end{align}
		\item (regularity) for $t\in[0,\dlt]$ and $0 \leq n \leq s$, \begin{equation}\label{eq:KdV-deg-wp-norm-estimate}
		\begin{split}
		\nrm{f^{-\frac{\sgm_{c}}{3}}(A^{-\frac{1}{3}}  f^{\frac{1}{3}}(t,x)\rd_x)^n \phi^{app}_{(\lmb)}(t,x)}_{L^2} \leq C_{\tld{f}} \abs{\lmb}^n \nrm{g_0}_{H^{n}_{ {(x_{1})}}};
		\end{split}
		\end{equation}
		\item (degeneration)  
for any $1 \leq p \leq 2$, a nonnegative even integer $s' \leq s$, and $\gmm' \geq -s'-\frac{1}{p}+\frac{1}{2}$, we have
		\begin{equation} \label{eq:KdV-wp-deg}
		f^{\frac{-\sgm_{c}+\gmm'}{3}} \phi^{app}_{(\lmb)} = \rd_{x}^{s'}\left( \frac{f^{\frac{-\sgm_{c}+\gmm'+s'}{3}}}{(-1)^{\frac{s'}{2}}  {A^{\frac{s'}{3}}} \lmb^{s'}}  \phi^{app}_{(\lmb)}\right) + f^{\frac{-\sgm_{c}+\gmm'}{3}} \phi^{small}_{(\lmb)},
\end{equation}
		where for $0 \leq j \leq 1$ and $t \in[0, \dlt]$, we have 
		\begin{align}
		\left\| \rd_{x}^{j} \left(  {A^{-\frac{s'}{3}}} \lmb^{-s'} f^{\frac{-\sgm_{c}+\gmm'+s'}{3}} \phi^{app}_{(\lmb)}\right) (t, x) \right\|_{L^{p}} 
		&\leq C_{\tld{f}}^{1+\gmm'}  {x_{1}^{(\gmm' + (s' - j) + \frac{1}{p} - \frac{1}{2})} A^{\frac{\gmm'}{3}} }\lmb^{-(s'-j)} \notag \\
		&\phantom{\leq} \times e^{-3 \bt(t)  {A^{\frac{2}{3}}} \lmb^{2}(\gmm'+(s'-j)+\frac{1}{p}-\frac{1}{2}) t} \nrm{g_{0}}_{ {H^{1}_{(x_{1})}}}, \label{eq:KdV-wp-deg-estimate} \\
		\nrm*{f^{-\frac{\sgm_{c}}{3}} \tphi^{small}_{(\lmb)}(t, x)}_{L^{2}} &\leq C_{\tld{f}} \lmb^{-1} \nrm{g_{0}}_{H^{s'}_{ {(x_{1})}}}; \label{eq:KdV-wp-deg-small}
		\end{align}
		\item (error estimate) letting 
		\begin{equation*}
		\err[\phi^{app}_{(\lmb)}] = (\rd_{t} + \calL_{f}) \phi^{app}_{(\lmb)},
\end{equation*}
		for $t\in[0,\dlt]$, we have \begin{equation} \label{eq:KdV-wp-err}
		\begin{split}
		\nrm{f^{-\frac{\sgm_{c}}{3}} \err[\phi^{app}_{(\lmb)}](t, x)}_{L^2} \le  {C_{\tld{f}} \left( 1 + \nrm{f}_{L^{\infty}_{t} C^{0, 1}}^{m-2} \right)} (1+ |\lmb|^2t ) |\lmb| \nrm{g_0}_{H^3_{ {(x_{1})}}}. 
		\end{split}
		\end{equation}
	\end{itemize} In the above properties, each constant referred to as $C_{\tld{f}}$ (resp.~$C_{\tld{f}_{0}}$) obeys the estimate
\begin{align*}
	C_{\tld{f}} &\leq C_{s} \exp \left( N_{s} A^{-1} \nrm{f}_{L^{\infty}_{t} \tld{C}^{s-1, 1}} \right), \\
	\bb(\hbox{resp.~} C_{\tld{f}_{0}} &\leq C_{s} \exp \left( N_{s} A^{-1} \nrm{f_{0}}_{\tld{C}^{s_{0}-1, 1}} \right), \bb)
\end{align*}
for some $C_{s} > 0$ and $N_{s} \in \bbN$ independent of $f$ and $x_{1}$ (but possibly dependent on $s$). 
\end{proposition} When $f$ or $g_{0}$ are clear from the context, we shall often simply omit them in $\phi^{app}_{(\lmb)}[g_{0}, f]$. 

\subsubsection{Renormalization and conjugation}
For the construction, we introduce the normalization
\begin{equation*}
	\tld{f}(t, x) := \frac{f(t, x)}{A}.
\end{equation*}By this normalization, we have $\tld{f}(0, x) = x^{3} + o_{f}(x^{3})$. Using $\tld{f}$, we define for  $t\in[0,\dlt]$ and $x \in (0,x_1]$ \begin{equation*}
\begin{split}
y(t,x) = -\int_x^{x_1} \frac{1}{{\tld{f}(t,x')}^{\frac{1}{3}}} \,\ud x' \le 0. 
\end{split}
\end{equation*} 
Note that $\rd_{x}y>0$. Then, we compute from $\rd_y = \tld{f}^{\frac{1}{3}}\rd_x$ that 
\begin{equation*}
\begin{split}
\tld{f}^{\frac{2}{3}} \rd_{xx} & = \rd_{yy} - \frac{1}{3} f^{-1}f_y \rd_y, \qquad
\tld{f}\rd_{xxx} = \rd_{yyy} - f^{-1}f_y \rd_{yy} +\left( -\frac{1}{3} f^{-1}f_{yy} + \frac{5}{9} f^{-2}(f_y)^2 \right)\rd_y .
\end{split}
\end{equation*} Furthermore, in the time derivative of $\phi$, \begin{equation*}
\begin{split}
	\rd_{t} \phi(t,x) = \rd_{t} \phi(t,y) + (\rd_{t}y ) \rd_{y}\phi(t,y)
\end{split}
\end{equation*} we set $q :=\rd_{t}y$ for simplicity. Note that in the $(t,x)$-coordinates, we have 
\begin{equation}\label{eq:q-in-x}
	A^{-1} q = \frac{1}{3} \int_{x}^{x_{1}} \frac{- \tld{f} \tld{f}_{xxx} - \alp_{1} \tld{f}_{x} \tld{f}_{xx} - \frac{\mu_{1}}{m} (f^{m-2} \tld{f}^{2})_{x}}{\tld{f}(t, x')^{\frac{4}{3}}}  \, \ud x'.
\end{equation}
Then, in the $(t,y)$-coordinates, \eqref{eq:KdV-pert-lin-re} transforms into 
\begin{equation} \label{eq:KdV-pert-lin-renorm}
\begin{aligned}
	&  {A^{-1}} \rd_{t} \phi + \phi_{yyy} + (\alp_{1} - 1) \tld{f}^{-1} \tld{f}_{y} \phi_{yy} \\
	&= -  {A^{-1}} q \phi_{y} + \left(\left(\frac{1}{3}-\alp_{1}\right) \tld{f}^{-1} \tld{f}_{yy} + \left(- \frac{5}{9} + \frac{2}{3} \alp_{1}  \right) \tld{f}^{-2} (\tld{f}_{y})^{2} + \mu_{1} f^{m-2}  {\tld{f}^{\frac{2}{3}}} \right) \phi_{y}\\
	&\relphantom{=}
	- \left( (\tld{f}^{-\frac{1}{3}} \rd_{y})^{3} \tld{f} + (m-1) \mu_{1} f^{m-3} f_{y}  {\tld{f}^{\frac{2}{3}}}  \right) \phi.
\end{aligned}
\end{equation}
We shall regard the expressions on the right-hand side of \eqref{eq:KdV-pert-lin-renorm} as error terms. To remove the last term on the left-hand side, we introduce the conjugated variable
\begin{equation} \label{eq:gdkdv-varphi}
	\varphi = e^{-G} \phi,
\end{equation}
where $G$ shall be determined below. We compute
\begin{align*}
	\rd_{y} \phi &= e^{G} \left( \rd_{y} \varphi + \rd_{y} G \varphi \right), \\
	\rd_{y}^{2} \phi &= e^{G} \left( \rd_{y}^{2} \varphi + 2 \rd_{y} G \rd_{y} \varphi + (\rd_{y}^{2} G + (\rd_{y} G)^{2}) \varphi \right), \\
	\rd_{y}^{3} \phi &= e^{G} \left( \rd_{y}^{3} \varphi + 3 \rd_{y} G \rd_{y}^{2} \varphi + (3 \rd_{y}^{2} G + 3 (\rd_{y} G)^{2}) \rd_{y} \varphi
	+ (\rd_{y}^{3} G + 3 \rd_{y} G \rd_{y}^{2} G + (\rd_{y} G)^{3}) \varphi \right).
\end{align*}
Hence, the left-hand side of \eqref{eq:KdV-pert-lin-renorm}, after factoring out $e^{G}$, becomes
\begin{align*}
	& {A}^{-1} \rd_{t} \varphi + \varphi_{yyy} + \left(3 G_{y} + (\alp_{1} -1) f^{-1} f_{y} \right) \varphi_{yy} \\
	& + (3 G_{yy} + 3 G_{y}^{2} + 2 (\alp_{1}-1) G_{y} f^{-1} f_{y}) \varphi_{y} \\
	&+ \left( {A}^{-1} G_{t} + (G_{yyy} + 3 G_{y} G_{yy} + G_{y}^{3} + (\alp_{1}-1) (G_{yy} + G_{y}^{2}) f^{-1} f_{y} \right) \varphi.
\end{align*}
The right-hand side, after factoring out $e^{G}$, becomes
\begin{align*}
	&-  {A}^{-1} q \varphi_{y} + \left(\left(\frac{1}{3}-\alp_{1}\right) f^{-1} f_{yy} + \left(- \frac{5}{9} + \frac{2}{3} \alp_{1}  \right) f^{-2} (f_{y})^{2} + \mu_{1}  {A}^{-1} f^{m-1}  {\tld{f}^{-\frac{1}{3}}} \right) \varphi_{y}\\
	&
	- \left( q  - \left(\left(\frac{1}{3}-\alp_{1}\right) f^{-1} f_{yy} + \left(- \frac{5}{9} + \frac{2}{3} \alp_{1}  \right) f^{-2} (f_{y})^{2} + \mu_{1}  {A^{-1}} f^{m-1}  {\tld{f}^{-\frac{1}{3}}} \right) \right) G_{y} \varphi \\
	& - \left((f^{-\frac{1}{3}} \rd_{y})^{3} f + (m-1) \mu_{1}  {A^{-1}} f^{m-2} f_{y}  {\tld{f}^{-\frac{1}{3}}}  \right) \varphi.
\end{align*}
To remove the second-order term, we are motivated to choose
\begin{equation*}
	G_{y} = - \frac{\alp_{1} - 1}{3} f^{-1} f_{y} = \frac{\sgm_{c}-\frac{1}{2}}{3} f^{-1} f_{y}.
\end{equation*}
Noting that $f^{-1} f_{y} = (\ln f)_{y}$, we see that $G(t, y) =  \left( \sgm_{c} - \tfrac{1}{2}\right) \frac{1}{3} \ln f(t, y) + C$ for some choice of $C$. We choose $C = \frac{1}{6} \ln A$ so that
\begin{equation} \label{eq:KdV-G}
	e^{G(t, y)} =  f^{\frac{\sgm_{c}}{3}} \tld{f}^{-\frac{1}{6}},
\end{equation}
in view of the weight $f^{-\frac{2}{3} \sgm_{c}} $ in the modified energy estimate we shall prove later (see also the definition of $L^{2}_{f}$ in Section~\ref{subsec:KdV-prelim}).

In conclusion, \eqref{eq:KdV-pert-lin-renorm}, after factoring out $e^{G}$, may be rewritten as
\begin{equation} \label{eq:KdV-varphi-err}
\begin{aligned}
	& {A^{-1}} \rd_{t} \varphi + \varphi_{yyy} \\
	&= \left(-  {A^{-1}} q + C_{1, 1} f^{-1} f_{yy} + C_{1, 2} f^{-2} (f_{y})^{2} + \mu_{1}  {A^{-1}} f^{m-1} { \tld{f}^{-\frac{1}{3}}} \right) \varphi_{y} \\
	&\relphantom{=}+ \left(-  {A^{-1}} G_{t} + C_{0, 1} f^{-1} f_{yyy} + C_{0, 2} f^{-2} f_{y} f_{yy} + C_{0, 3} f^{-3} (f_{y})^{3} + C_{0, 4} \mu_{1}  {A^{-1}} f^{m-2} f_{y}  {\tld{f}^{-\frac{1}{3}}}\right) \varphi,
\end{aligned}
\end{equation}
where $C_{j, k} \in \bbR$ are constants that depend on $\alp_{1}$, $\mu_{1}$ and $m$.

\subsubsection{Specification of the wave packet and the Proof of Proposition~\ref{prop:KdV-deg-wp}}

Given $g_{0}(x)$, we set $h_0(y) ={x_{1}^{\frac{1}{2}}} g_{0}(x(0, y))$, and for each $\lmb\in\bbN$, we first take $\varphi^{app}_{(\lmb)}[g_{0},f]$ to be the standard wave packet for the Airy equation {with time rescaled by $A$ and} with frequency $\lmb$, i.e.,
\begin{equation} \label{eq:gdkdv-wp-varphi}
	\varphi^{app}_{(\lmb)}[g_{0}, f](t, y) = \Re(e^{i \lmb (y +  {A}  \lmb^{2}  t)} h_{0}(y + 3  {A}  \lmb^{2} t))
	= \cos(\lmb (y +  {A} \lmb^{2}   t)) h_{0}(y + 3  {A} \lmb^{2} t).
\end{equation}
Then, we define the degenerating wave packet $\phi^{app}_{(\lmb)}[g_{0},f]$ by $e^{G}\varphi^{app}_{(\lmb)}[g_{0}, f]$. Explicitly, we have
\begin{equation} \label{eq:gdkdv-wp}
	\phi^{app}_{(\lmb)}[g_{0}, f] =  f(t, y)^{\frac{\sgm_{c}}{3}} \tld{f}(t,y)^{-\frac16} \cos(\lmb (y +  {A} \lmb^{2} t)) h_{0}(y + 3  {A} \lmb^{2} t).
\end{equation}

\begin{proof}[Proof of Proposition~\ref{prop:KdV-deg-wp}]
	Now that we have specified the construction of $\phi^{app}_{(\lmb)}[g_{0}, f]$, we verify its properties claimed in Proposition~\ref{prop:KdV-deg-wp}. {In what follows, the dependence of constants on $f$, $A$ and $x_{1}$ has been made explicit.} {Moreover, we shall use the notation $C_{\tld{f}}$ introduced in Proposition~\ref{prop:KdV-deg-wp}.} 
	
	\medskip
	
	\noindent \textit{Linearity and support property}. To begin with, the linearity property is clear. To prove the support property, we first note from \eqref{eq:KdV-cubic-profile} and the positivity of $\tld{f}$ that $y < 0$ implies $x < x_{1}$, and vice versa. Note also that
\begin{equation*}
	y(t, x) = \frac{1}{\tld{\bt}(t)} \left(\ln \frac{x}{x_{1}} + B(t, x) \right),
\end{equation*} {where
\begin{equation*}
	\tld{\bt}(t) := \frac{\bt(t)}{A^{\frac{1}{3}}}, \qquad \abs{B(t, x)} \leq C x_{1} \nrm{\tld{f}}_{L^{\infty}_{t} C^{3, 1}}.
\end{equation*}
Here, $\bt(t)^{3} = \frac{1}{6} f_{xxx}(t, 0)$ as in Lemma~\ref{lem:KdV-C3-solution}; observe that $\tld{\bt}(0) = 1$ by definition.} The above formula for $y$ gives
\begin{equation} \label{eq:gdkdv-xy}
\begin{split}
x(t,y) = x_1 e^{\tld{\bt}(t)y - B} ,
\end{split}
\end{equation} 
from which the rest of the support property follows. 

 From \eqref{eq:gdkdv-xy}, it follows that 
\begin{equation} \label{eq:gdkdv-fy}
\begin{aligned}
\tld{f}(t,y) &= x_1^3 \tld{\bt}^3(t)  e^{3\tld{\bt}(t)y-3B} \left( 1  +  {O\left( x_1 e^{\tld{\bt}(t)y-B} \nrm{\tld{f}}_{L^{\infty}_{t} C^{3, 1}} \right) }\right) \\
	& \leq \exp(C \nrm{f}_{L^{\infty}_{t} \tld{C}^{3, 1}}) x_1^3 \tld{\bt}^3(t)  e^{3\tld{\bt}(t)y}. 
\end{aligned}
\end{equation} 
Using the control of $\nrm{f}_{L^{\infty}_{t} Y}$, we furthermore have
\begin{equation} \label{eq:gdkdv-dfy}
	\abs{f_{y}} \leq C \nrm{\tld{f}}_{L^{\infty}_{t} Y}^{\frac{1}{3}} \abs{f}, \quad
	\abs{f_{yy}} \leq C \nrm{\tld{f}}_{L^{\infty}_{t} Y}^{\frac{2}{3}} \abs{f}, \quad
	\abs{f_{yyy}} \leq C \nrm{\tld{f}}_{L^{\infty}_{t} Y}  \abs{f}.
\end{equation}
For higher derivatives, it is straightforward to verify by induction that
\begin{equation} \label{eq:gdkdv-dkfy}
	\abs{\rd_{y}^{k} f} \leq C_{k} \nrm{\tld{f}}_{L^{\infty}_{t} \tld{C}^{k-1, 1}}^{\frac{k}{3}} \abs{f} \quad \hbox{ for } k \geq 4.
\end{equation} We furthermore claim that, by for any integer $0 \leq k \leq s+1$,
\begin{equation} \label{eq:gdkdv-wp-id}
	\nrm{h_{0}}_{H^{k}_{y}} \leq C_{k} \left( 1 + \nrm{\tld{f}}_{L^{\infty}_{t} \tld{C}^{s-1, 1}}^{\frac{k-1}{3}}\right) \nrm{g_{0}}_{H^{k}_{(x_{1})}}.
\end{equation}
Indeed, arguing via induction in a similar fashion as above, we may verify that, for any $0 \leq k \leq s+1$,
\begin{equation*}
	\sum_{j' = 0}^{k} \abs{\rd_{y}^{j'} h_{0}} \leq x_{1}^{\frac{1}{2}} \left( \abs{g_{0}} + C_{k} \sum_{j' = 1}^{k} \sum_{j = 1}^{j'} \nrm{\tld{f}_{0}}_{L^{\infty}_{t} \tld{C}^{s-1, 1}}^{\frac{j'-j}{3}} \tld{f}_{0}^{\frac{j}{3}} \abs{\rd_{x}^{j} g_{0}} \right).
\end{equation*}
Note furthermore that, by \eqref{eq:gdkdv-x1}, $C^{-1} x_{1}^{3} \leq \tld{f}(x) \leq C x_{1}^{3}$ on $\supp g_{0}$. Taking the $L^{2}_{y}$ norm of both sides and changing variables, we are led to \eqref{eq:gdkdv-wp-id}.

\medskip

\noindent \textit{Initial data and regularity estimates}. Let us now verify the initial data and regularity estimates. We begin by noting that
\begin{align*}
	\int f^{-\frac{2}{3} \sgm_{c}} \phi^{app}(t, x)^{2} \, \ud x
	&= \int  {(f^{-\frac{\sgm_{c}}{3}} \tld{f}^{\frac{1}{6}} \phi^{app}(t, y))^{2}} \, \ud y  = \int \varphi^{app}(t, y)^{2} \, \ud y 
	=\nrm{h_{0}}_{L^{2}_{y}}^{2},
\end{align*}
from which the regularity estimate in the case $n = 0$ follows. Moreover, from this identity it is clear that $\nrm{f^{-\frac{\sgm_{c}}{3}} \phi^{app}(0, x)  }_{L^2} \le C_{\tld{f}_{0}} \nrm{g_0}_{L^2}$. To obtain the claimed lower bound, first note that \begin{equation*}
	\begin{split}
		\nrm{\cos(\lmb y) h_{0}}_{L^2_{y}}^{2} = \int  \frac1{\lmb} \rd_{y} ( \sin(2\lmb y) ) h_{0}^{2}(y)  \ud y + \frac12 \nrm{h_{0}}_{L^{2}_y}^{2}
	\end{split}
\end{equation*} holds, and then one can integrate by parts in the first term on the right hand side, with $\nrm{h_{0}}_{L^{2}_y} \gtrsim_{\tld{f}_{0}} \nrm{g_{0}}_{L^{2}}$.

Next, when $n = 1$, we note that 
\begin{align*}
	\rd_{y} \phi^{app}(t, x)
	&= \Re\left( i \lmb   {f(t, y)^{\frac{\sgm_{c}}{3}} \tld{f}^{\frac{1}{6}}} e^{i \lmb (y +  {A} \lmb^{2} t)} h_{0}(y+3  {A} \lmb^{2} t)\right)\\
	&\relphantom{=}
	+\frac{\sgm_{c}-\frac{1}{2}}{3} \frac{\rd_{y} f}{f} {f(t, y)^{\frac{\sgm_{c}}{3}} \tld{f}^{\frac{1}{6}}} \Re\left( e^{i \lmb (y +  {A} \lmb^{2} t)} h_{0}(y+3 {A} \lmb^{2} t)\right) \\
	&\relphantom{=}
	+ {f(t, y)^{\frac{\sgm_{c}}{3}} \tld{f}^{\frac{1}{6}}} \Re\left( e^{i \lmb (y + {A} \lmb^{2} t)} \rd_{y} h_{0}(y+3  {A} \lmb^{2} t) \right).
\end{align*}
From this expression, the regularity estimate follows by the earlier computation; the power $\abs{\lmb}$ arises from the first term, the second term is estimated using \eqref{eq:gdkdv-dfy}, and the need for the $H^{1}_{ {(x_{1})}}$ norm of $g_{0}$ is due to the third term. The case of higher $n$ can be handled similarly; we omit the details. Lastly, the initial data estimate can be proved simply by taking $t = 0$; see the computations below for $s'=0$. 


\medskip

\noindent \textit{Degeneration property}. When $s' = 0$, we simply set $\phi^{small}_{(\lmb)} = 0$. Arguing as in the proof of the regularity property, we have
\begin{align*}
	\nrm{f(t, x)^{\frac{-\sgm_{c}+\gmm'}{3}} \phi^{app}_{(\lmb)}(t, x)}_{L^{p}_{x}}^{p}
	&= \int \abs{\varphi^{app}_{(\lmb)}(t, y)}^{p}  {f(t, y)^{\frac{p \gmm'}{3}} \tld{f}(t, y)^{\frac{-\frac{p}{2}+1}{3}}} \, \ud y \\
	&\leq  {A^{\frac{p \gmm'}{3}}} \left( \int \abs{\varphi^{app}_{(\lmb)}(t, y)}^{2} \, \ud y\right)^{\frac{p}{2}} \left(\int_{\supp \varphi^{app}_{(\lmb)}(t, \cdot)} \tld{f}(t, y)^{\frac{\frac{1}{\frac{1}{p}-\frac{1}{2}} \gmm' + 1}{3}} \, \ud y\right)^{p\left(\frac{1}{p}-\frac{1}{2}\right)} \\
	&\leq  \nrm{h_{0}}_{L^{2}}^{p}  {A^{\frac{p \gmm'}{3}}} ( {C_{\tld{f}} x_{1} \tld{\bt}})^{p\left(\gmm'+\frac{1}{p}-\frac{1}{2}\right)} \left(\int_{-\infty}^{-3  {A} \lmb^{2} t} e^{\left( \frac{1}{\frac{1}{p}-\frac{1}{2}} \gmm' + 1 \right)  {\tld{\bt}}(t) y} \, \ud y\right)^{p\left(\frac{1}{p}-\frac{1}{2}\right)} \\
	&\leq \nrm{h_{0}}_{L^{2}}^{p}  {(C_{\tld{f}} x_{1})}^{p\left(\gmm'+\frac{1}{p}-\frac{1}{2}\right)}  {(A^{\frac{1}{3}}\tld{\bt})}^{p \gmm'} e^{- 3 p  {\tld{\bt}}(t)  {A} \lmb^{2} (\gmm'+\frac{1}{p}-\frac{1}{2}) t},
\end{align*} where we have simply used \eqref{eq:gdkdv-fy} to bound $\tld{f}(t,y)$. This proves \eqref{eq:KdV-wp-deg-estimate} in the case $s = 0$.

To handle the case $s > 0$, it is convenient to introduce the following notation (as in the Schr\"odinger case): given some function $r=r(t,y)$, 
\begin{equation*}
	H = rO_{k}(h_{0}) \quad \impmi \quad \sup_{t \in [0, \dlt]} \nrm*{\tld{f}^{\frac{1}{6}} \frac{H}{r}}_{L^{2}_{y}} \leq  {C_{\tld{f}}} \nrm{h_{0}}_{H^{k}_{y}}.
\end{equation*}
In this case, note that $\nrm{\tld{f}^{\frac{1}{6}} (\cdot)}_{L^{2}_{y}} = \nrm{\cdot}_{L^{2}_{x}}$ for each $t$. We shall also freely use \eqref{eq:gdkdv-wp-id} to relate the right-hand side with $\nrm{g_{0}}_{H^{k}_{(L)}}$. In what follows, the expression abbreviated as $\frac{1}{\lmb} O_{k}(h_{0})$ constitutes $f^{-\frac{\sgm_{c}}{3}} \phi^{small}_{(\lmb)}$; the desired estimate \eqref{eq:KdV-wp-deg-small} would be an immediate consequence of the $L^{2}$ boundedness property embedded in the $O_{k}(\cdot)$ notation.

We treat the case $s = 2$. We begin with the identity
\begin{equation*}
	\cos (\lmb(y +  {A} \lmb^{2} t)) = - \frac{f^{\frac{2}{3}}}{ {A^{\frac{2}{3}}}\lmb^{2}} \left( {\tld{f}}^{-\frac{1}{3}} \rd_{y} \right)^{2} \cos (\lmb(y +  {A} \lmb^{2} t))  {- \frac{1}{3\lmb} f^{-1}\rd_{y}f \sin(\lmb(y+ A \lmb^{2}t)) }. 
\end{equation*}
Plugging this identity into the expression \eqref{eq:gdkdv-wp} for $\phi^{app}_{(\lmb)}[g_{0}, f]$ and commuting $\left(\tld{f}^{-\frac{1}{3}} \rd_{y}\right)^{2}$ (which equals $\rd_{x}^{2}$ in the $(t, x)$-coordinates) outside, we have
\begin{equation*}
	f(t, y)^{\frac{-\sgm_{c}+\gmm'}{3}} \phi^{app}_{(\lmb)}
	= \left(\tld{f}^{-\frac{1}{3}} \rd_{y}\right)^{2} \left( \frac{f(t, y)^{\frac{-\sgm_{c}+\gmm'+2}{3}}}{(-1)  {A^{\frac{2}{3}}} \lmb^{2}} \phi^{app}_{(\lmb)} \right)
	+ \frac{f(t, y)^{\frac{\gmm'}{3}}}{\lmb} O_{2}(h_{0}).
\end{equation*}
Arguing as in the case $s = 0$, the expression inside the parentheses can be shown to obey the degeneration bound \eqref{eq:KdV-wp-deg-estimate}. The cases $s > 2$ are handled similarly. 

\medskip

\noindent \textit{Error bound}. We begin by noticing that, by our construction, we have
\begin{equation*}
	\nrm{f^{-\frac{\sgm_{c}}{3}} \err[\phi^{app}_{(\lmb)}]}_{L^{2}_{x}}
	\leq \nrm{(\rd_{t} + \rd_{yyy}) \varphi^{app}_{(\lmb)}}_{L^{2}_{y}} + \nrm{(\hbox{RHS of \eqref{eq:KdV-varphi-err}})}_{L^{2}_{y}}.
\end{equation*}
The first term is the error for the standard wave packet for the Airy equation with frequency $\lmb$; it is easily bounded by $C \abs{\lmb} \nrm{h_{0}}_{H^{3}_{y}}$, which is acceptable. Now, it only remains to estimate the $L^{2}_{y}$ norm of each term on the right-hand side of \eqref{eq:KdV-varphi-err}. The worst contribution turns out to be $- q \varphi_{y}$, which we turn to first. By \eqref{eq:gdkdv-fy}, \eqref{eq:gdkdv-dfy}, and the definition of $q$, it follows that 
\begin{equation}\label{eq:q}
	\abs{ {A^{-1} q(t,y)}} \leq C_{\tld{f}} (1+\nrm{f}_{L^{\infty}_{t} C^{0, 1}}^{m-2}) \int_{y}^{0} \, \ud y' \leq C_{\tld{f}} (1+\nrm{f}_{L^{\infty}_{t} C^{0, 1}}^{m-2}) \abs{y}.
\end{equation}
By the support property of $\varphi^{app}_{(\lmb)}$, we have 
\begin{equation*}
	\nrm{A^{-1} q (\varphi^{app}_{(\lmb)})_{y}}_{L^{2}} \le  {C_{\tld{f}} (1+\nrm{f}_{L^{\infty}_{t} C^{0, 1}}^{m-2})} (1+  {A} \lmb^{2} t) \abs{\lmb} \nrm{h_{0}}_{H^{1}_{y}},
\end{equation*}
which is acceptable. The remaining terms on the right-hand side of \eqref{eq:KdV-varphi-err} involving $\varphi_{y}$ are bounded by  {$C_{\tld{f}} (1+\nrm{f}_{L^{\infty}_{t} C^{0, 1}}^{m-2}) \abs{\lmb} \nrm{h_{0}}_{H^{1}_{y}}$}, which are strictly better. Next, since  {$\abs{A^{-1} \rd_{t} f} \leq C_{\tld{f}} (1+\nrm{f}_{L^{\infty}_{t} C^{0, 1}}^{m-2}) \abs{f}$ (as in the estimate for $q$)}, we have
\begin{equation*}
	 {\abs{A^{-1} \rd_{t} G} \leq C_{\tld{f}} (1+\nrm{f}_{L^{\infty}_{t} C^{0, 1}}^{m-2}).}
\end{equation*}
Using this bound, as well as \eqref{eq:gdkdv-fy} and \eqref{eq:gdkdv-dfy}, the terms on the right-hand side of \eqref{eq:KdV-varphi-err} involving $\varphi$ are bounded by  {$C_{\tld{f}} (1+\nrm{f}_{L^{\infty}_{t} C^{0, 1}}^{m-2}) \nrm{\varphi}_{L^{2}_{y}}$}, which is good. This completes the proof of the \eqref{eq:KdV-wp-err}. 
\end{proof}

\subsection{Modified energy estimate for the perturbation} \label{subsec:gdkdv-en}

Recall the equation  {satisfied} by $\phi$: 
 {\begin{equation}\label{eq:gdkdv-cubic-phi}
\begin{aligned}
&\rd_{t} \phi + f \phi_{xxx} + \alp_{1} f_{x} \phi_{xx} + (\alp_{1} f_{xx} + \mu_{1} f^{m-1}) \phi_{x} + (f_{xxx} + (m-1) \mu_{1} f^{m-2} f_{x}) \phi \\
&= - \phi \phi_{xxx} - \alp_{1} \phi_{x} \phi_{xx} - \mu_{1} \left(  (f + \phi)^{m} - f^{m} - m f^{m-1} \phi \right)_{x}.
\end{aligned}
\end{equation} }
Regarding a solution $\phi$ of the above, we have the following modified energy estimate. 
\begin{proposition}\label{prop:KdV-mee}
	 {Assume that $f$ is a  solution to \eqref{eq:gdkdv} satisfying $f \in L^{\infty}([0, \dlt]; \tld{C}^{3, \alp}(I))$ with initial data $f_{0}$ that is positive on $I \setminus \rd I$ and vanishes to order at least $3$ on each point in $\rd I$.} Moreover, assume that $\phi\in  L^\infty([0,\dlt]; C^{3, \alp}(I))$  is a solution to \eqref{eq:gdkdv-cubic-phi} satisfying \begin{equation*}
	\begin{split}
	f + \phi \in L^{\infty}([0, \dlt]; \tld{C}^{3, \alp}(I)), \qquad f^{-1}(f+\phi) \in L^{\infty}([0, \dlt]; L^{\infty}(I)).
	\end{split}
	\end{equation*} Then we have the estimate \begin{equation*}
	\begin{split}
	\int_{ {I}} \phi^2(t,x)  {f(t,x)^{-\frac{2}{3}\sgm_{c}}}\,\ud x \le  {\exp\left(C_{f, f+\phi} t \right)}  \int_{ {I}} \phi_{0}^2(x)  {f(t,x)^{-\frac{2}{3}\sgm_{c}}}\,\ud x
	\end{split}
	\end{equation*} for $t\in[0,\dlt]$, where $\phi_{0}(x) = \phi(0, x)$ and
	\begin{equation}\label{eq:Cf}
C_{f, f+\phi} \leq C \sup_{t\in[0,\dlt]}  (\nrm{f}_{Y}+  (1+ \nrm{ f^{-1}(f+\phi) }_{L^{\infty}}^{\frac12})  \nrm{f + \phi}_{Y} + (\nrm{f}_{C^{0, 1}}+ \nrm{f + \phi}_{C^{0, 1}})^{m-1} ),
\end{equation}
	with $C>0$ an absolute constant.
\end{proposition}


\begin{proof}
	In what follows, we shall simply present a formal computation without worrying about the validity of the expressions and manipulations. {Also, all integrals are taken over $I$.} As in \S\ref{subsec:modified-energy-estimate}, the assumption $\abs{\phi_{0}(x)} \leq C \abs{f_{0}(x)}$ and the finiteness of the right-hand side of $\int \phi^2_0(x)  {f(t,x)^{-\frac{2}{3} \sgm_{c}}}\,\ud x  < +\infty$ would imply, via Lemma~\ref{lem:KdV-C3-solution}, the vanishing property of $\phi(t, \cdot)$ on $\rd I$ that is necessary to justify the computation; we shall leave the routine details to the reader. 
	
	To prove the proposition, we compute \begin{equation*}
	\begin{split}
	\frac{\ud}{\ud t}\int \phi^{2}  {f(t,x)^{-\frac{2}{3}\sgm_{c}}} \,\ud x & = \int 2\phi\rd_t\phi  {f(t,x)^{-\frac{2}{3}\sgm_{c}}} \,\ud x + \int \phi^{2} \rd_t  f(t,x)^{-\frac{2}{3}\sgm_{c}} \,\ud x ,
	\end{split}
	\end{equation*} and the last term can be bounded as in the proof of \eqref{eq:KdV-f-pointwise}; we have \begin{equation*}
	\begin{split}
	\left|\int \phi^{2} \rd_t  {f(t,x)^{-\frac{2}{3}\sgm_{c}}} \,\ud x   \right|\le C \nrm{f}_{Y} \int \phi^{2}   {f(t,x)^{-\frac{2}{3}\sgm_{c}}} \,\ud x. 
	\end{split}
	\end{equation*} We decompose the other term in the right hand side as follows, up to a factor of 2: \begin{equation*}
	\begin{split}
	I &= -\int \phi \left( f\phi_{xxx} + \phi f_{xxx} +  {\alp_{1}} f_x \phi_{xx} +  {\alp_{1}} f_{xx}\phi_{x} \right)  {f^{-\frac{2}{3}\sgm_{c}}} \,\ud x, \\
	II & = - \mu_{1} \int \phi \left( (m-1) \phi f^{m-2} f_x + f^{m-1} \phi_x \right)  {f^{-\frac{2}{3}\sgm_{c}}} \,\ud x , \\
	III& = - \int \phi \, {Q[\phi]}  {f^{-\frac{2}{3}\sgm_{c}}} \,\ud x . 
	\end{split}
	\end{equation*}  To estimate $I$, we observe the following chain of inequalities: \begin{equation*}
	\begin{split}
	\left| \int \phi^2f_{xxx}  {f^{-\frac{2}{3}\sgm_{c}}} \,\ud x\right|\le C \nrm{f}_{Y} \int \phi^2  {f^{-\frac{2}{3}\sgm_{c}}}\,\ud x, 
	\end{split}
	\end{equation*} \begin{equation*}
	\begin{split}
	\left| \int  {\alp_{1}}\phi\phi_x f_{xx}  {f^{-\frac{2}{3}\sgm_{c}}} \,\ud x  \right| = \left| \frac{ {\alp_{1}}}{2}\int \phi^2 \rd_x(f_{xx}  {f^{-\frac{2}{3}\sgm_{c}}}) \,\ud x  \right|,
	\end{split}
	\end{equation*} \begin{equation*}
	\begin{split}
	\int  {\alp_{1}} \phi\phi_{xx} f_{x}  {f^{-\frac{2}{3}\sgm_{c}}} \,\ud x = -  {\alp_{1}} \int (\phi_x)^2 f_x  {f^{-\frac{2}{3}\sgm_{c}}} \,\ud x + \frac{ {\alp_{1}}}{2}\int \phi^2\rd_{xx}(f_x  {f^{-\frac{2}{3}\sgm_{c}}}) \,\ud x , 
	\end{split}
	\end{equation*} and lastly \begin{equation*}
	\begin{split}
	\int \phi\phi_{xxx}f  {f^{-\frac{2}{3}\sgm_{c}}} \,\ud x 
	&= -\int\phi_{xx}\phi_x f  {f^{-\frac{2}{3}\sgm_{c}}} - \phi_{xx}\phi\rd_x(f  {f^{-\frac{2}{3}\sgm_{c}}})  \,\ud x \\
	&= \frac{3}{2} \int(\phi_x)^2\rd_x(f {f^{-\frac{2}{3}\sgm_{c}}}) \,\ud x - \frac{1}{2}\int\phi^2\rd_{xxx}(f {f^{-\frac{2}{3}\sgm_{c}}}) \,\ud x. 
	\end{split}
	\end{equation*} From \begin{equation*}
	\begin{split}
	\rd_x(f {f^{-\frac{2}{3}\sgm_{c}}}) = \left( 1-\frac{2}{3}\sgm_{c} \right)  {f^{-\frac{2}{3}\sgm_{c}}} f_x
	 {= \frac{2}{3} \alp_{1} f^{-\frac{2}{3}\sgm_{c}} f_{x}}, 
	\end{split}
	\end{equation*} we obtain a cancellation of terms involving $(\phi_x)^2$ and then we observe \begin{equation*}
	\begin{split}
	|\rd_x(f_{xx} {f^{-\frac{2}{3}\sgm_{c}}})| + |\rd_{xx}(f_x {f^{-\frac{2}{3}\sgm_{c}}})| + |\rd_{xxx}(f {f^{-\frac{2}{3}\sgm_{c}}})| \le C \nrm{f}_{Y}   {f^{-\frac{2}{3}\sgm_{c}}}
	\end{split}
	\end{equation*} to conclude the estimate \begin{equation*}
	\begin{split}
	|I|\le C \nrm{f}_{Y}  \left| \int \phi^2  {f^{-\frac{2}{3}\sgm_{c}}} \,\ud x  \right| . 
	\end{split}
	\end{equation*} Next, to treat $II$ we simply integrate by parts: \begin{equation*}
	\begin{split}
	\left|II\right| 
	&= {\left| \mu_{1} \int \phi^2 \left((m-1) f^{m-2} f_{x} - \frac{-\frac{2}{3}\sgm_{c} + m-1}{2} f^{m-2} f_{x} \right) f^{-\frac{2}{3}\sgm_{c}}\,\ud x \right|}\le C  {\nrm{f}_{C^{0, 1}}^{m-1}}  \int \phi^2  {f^{-\frac{2}{3}\sgm_{c}}}\,\ud x. 
	\end{split}
	\end{equation*} Finally, we turn to $III$. Observe that we are free to use $\nrm{\phi}_{C^{3, \alp}}$, since it is controlled by $\nrm{f}_{C^{3, \alp}} + \nrm{f + \phi}_{C^{3, \alp}}$; similarly for $\nrm{\phi}_{C^{0, 1}}$. Recall the expression for $Q[\phi]$ given in \eqref{eq:gdkdv-Q}. We first estimate \begin{equation*}
	\begin{split}
	\left| \int \phi (- \phi \phi_{xxx} ) f^{-\frac{2}{3}\sgm_{c}} \, \ud x\right|
	\leq C \nrm{\phi}_{C^{2, 1}} \int \phi^{2} f^{-\frac{2}{3}\sgm_{c}} \, \ud x,
	\end{split}
	\end{equation*}  
	and
	\begin{align*}
	&\left| \int \phi (- \alp_{1} \phi_{x} \phi_{xx}  ) f^{-\frac{2}{3}\sgm_{c}} \, \ud x\right| \\
	&\leq \left| \int \frac{\alp_{1}}{2} \phi^{2} \phi_{xxx} f^{-\frac{2}{3}\sgm_{c}} \, \ud x \right| 
	+ \left| \int \frac{\alp_{1} \sgm_{c}}{3} \phi^{2} (\phi+f)_{xx} f^{\frac{-2\sgm_{c}-3}{3}} f_{x} \, \ud x \right| 
	+ \left| \int \frac{\alp_{1} \sgm_{c}}{3} \phi^{2} f_{xx} f^{\frac{-2\sgm_{c}-3}{3}} f_{x} \, \ud x \right|\\
	&\leq \left| \int \frac{\alp_{1}}{2} \phi^{2} \phi_{xxx} f^{-\frac{2}{3}\sgm_{c}} \, \ud x \right| 
	+ \left| \int \frac{\alp_{1} \sgm_{c}}{3} \phi^{2} (\phi+f)_{xx} f^{\frac{-2\sgm_{c}-3}{3}} f_{x} \, \ud x \right| 
	+ \left| \int \frac{\alp_{1} \sgm_{c}}{3} \phi^{2} f_{xx} f^{\frac{-2\sgm_{c}-3}{3}} f_{x} \, \ud x \right|\\
	&\leq C (\nrm{\phi}_{C^{2, 1}} +   \nrm{f^{-1}(f+\phi)}_{L^{\infty}}^{\frac13} \nrm{f+\phi}_{Y}^{\frac{2}{3}} \nrm{f}_{Y}^{\frac{1}{3}} + \nrm{f}_{Y})  \int \phi^{2} f^{0\frac{2}{3}\sgm_{c}} \, \ud x.
	\end{align*}   
	The remaining terms in $Q[\phi]$ are easier to treat, and collecting the estimates, we conclude \begin{equation*}
		\begin{split}
			\left| \frac{\ud}{\ud t} \int \phi^2 f^{-\frac{2}{3}\sgm_{c}} \,\ud x \right|\le C_{f,f+\phi} \int \phi^2 f^{-\frac{2}{3}\sgm_{c}} \,\ud x
		\end{split}
	\end{equation*} for some $C_{f,f+\phi}$ satisfying \eqref{eq:Cf}. The proposition follows by integrating in time.	 \end{proof}

Henceforth, we shall use the notation $\brk{\cdot, \cdot}_{f}$ and $\nrm{\cdot}_{L^{2}_{f}}$ introduced in Section~\ref{subsec:KdV-prelim} for the weighted $L^{2}$-inner product and norm on $I$ in Proposition~\ref{prop:KdV-mee}.

\subsection{Generalized energy estimate and proof of Theorem~\ref{thm:ill-posed-unbounded-gdkdv}} \label{subsec:unbounded-gdkdv}
Let $f$ satisfy the assumptions of Theorem~\ref{thm:ill-posed-unbounded-gdkdv}. Without loss of generality, we may assume that $a = 0$ and $0 < \eps \leq 1$. For simplicity, we shall focus on the case $\bt(0) = (\frac{1}{6} f_{0, xxx}(0))^{\frac{1}{3}} = 1$, the general case being analogous. Fix $0 < x_{1} < b$ so that \eqref{eq:gdkdv-x1} holds. Fix $g_{0} \in C^{\infty}_{c}$ supported in $(\frac{x_{1}}{2}, x_{1})$ with normalization $\nrm{g_{0}}_{L^{2}} = 1$. 
In what follows, {\bf we shall suppress the dependence of constants on  {$f$ and} $g_{0}$, in addition to $\alp_{1}$, $\mu_{1}$ and $m$ as before.}  {Also, {\bf we write $C(M)$ for a positive strictly increasing function of $M \in (0, \infty)$ such that $C(M) \to \infty$ as $M \to \infty$, which may vary from line to line.}}

Let $\phi^{app}_{(\lmb)} =\phi^{app}_{(\lmb)}[g_{0}, f]$ according to Proposition~\ref{prop:KdV-deg-wp}. We shall take
\begin{equation*}
	\phi_{0} = c_{0} \eps \lmb^{-m_{0}} x_{1}^{\frac{1}{2}} \phi^{app}_{(\lmb)}(0)
\end{equation*}
where $c_{0}$ is chosen so that $\nrm{\phi_{0}}_{C^{m_{0}}} \leq \eps$ using \eqref{eq:KdV-deg-wp-norm-estimate}, and $\lmb$ is to be determined below. Furthermore, by the normalization $\nrm{g_{0}}_{L^{2}} = 1$, we have 
\begin{equation} \label{eq:ill-posed-unbounded-gdkdv-dual0}
	\frac{1}{C_{0}} \abs{\lmb}^{-m_{0}}  \leq \nrm{\phi_{0}}_{L^{2}_{f_{0}}} \leq C_{0} \abs{\lmb}^{-m_{0}}, \quad 
	\brk{\phi_{0}, \phi^{app}_{(\lmb)}(0)}_{f_{0}}
	\geq \frac{1}{C_{0}} \nrm{\phi_{0}}_{L^{2}_{f_{0}}},
\end{equation}
 {for some constant $C_{0} > 0$ (which, in fact, depends on $\nrm{f}_{\tld{C}^{s_{0}-1, 1}}$).}
At this point, it is easy to ensure that \eqref{eq:phi-f-ratio-initial} is satisfied with $C_1=2$, where $u_0=f_0+\phi_0$.

Fix also $0 < \dlt' \leq \dlt$. To prove the theorem, we assume that the first alternative does not hold, i.e., there exists a solution $f + \phi \in L^{\infty}([0, \dlt']; C^{s-1, 1}(I))$ to \eqref{eq:gdkdv}. By Proposition~\ref{prop:C41} (and since $s \geq 5$), there exists $0 < t_{0} \leq \dlt'$ depending only on $\nrm{f(0)}_{Y}$ and $\nrm{(f + \phi)(0)}_{Y}$ such that $f, f+\phi \in L^{\infty}([0, t_{0}]; Y(I))$, and  $f^{-1}(f+\phi)\in L^\infty([0, t_{0}]; L^{\infty}(I))$, Moreover, by the same proposition, we have
 {\begin{equation} \label{eq:ill-posed-unbounded-gdkdv-Cs}
\sup_{0< t < t_{0} } \left( \nrm{f(t)}_{\tld{C}^{4, 1}(I)} + \nrm{(f+\phi)(t)}_{\tld{C}^{4, 1}(I)} + \nrm{f^{-1}(f+\phi)(t)}_{L^{\infty}(I)}\right) \leq C(M_{5}),
\end{equation}
where $M_{5} := \nrm{\phi}_{L^{\infty}([0, t_{0}]; C^{4, 1}(I))}$. (Here, we remind the reader our convention of omitting the dependence on $f$ in this proof.) By Proposition~\ref{prop:KdV-mee} and the preceding bound, we have 
\begin{equation} \label{eq:ill-posed-unbounded-gdkdv-mee}
	\nrm{\phi(t)}_{L^{2}_{f}} \leq \exp(C_{1}(M_{5}) t)\nrm{\phi_{0}}_{L^{2}_{f_{0}}},
\end{equation}
for some positive strictly increasing function $C_{1}(\cdot)$ that diverges at infinity. We emphasize that this function is \emph{independent} of $\lmb$, although $M_{5} = \nrm{\phi}_{L^{\infty}([0, t_{0}]; C^{4, 1}(I))}$ might be dependent on $\lmb$. }

Now using that $\phi$ is a solution to \eqref{eq:gdkdv-cubic-phi} and $(\rd_{t} + \calL_{f}) \phi^{app}_{(\lmb)} = \err[\phi^{app}_{(\lmb)}]$, we compute that \begin{equation*}
\begin{split}
\frac{\ud}{\ud t} \brk{\phi,\phi^{app}_{(\lmb)}}_{f} 
& = -\brk{\phi, \calL_{f}[\phi^{app}_{(\lmb)}]}_{f} + \brk{ \phi,\err[\phi^{app}_{(\lmb)}]}_{f} \\
&\relphantom{=} - \brk{\calL_{f}[\phi], \phi^{app}_{(\lmb)}}_{f} + \brk{Q_{f}[\phi],\phi^{app}_{(\lmb)}}_f \\
&\relphantom{=} - \frac{2}{3} \sgm_{c} \brk{ f^{-1}\rd_tf \phi ,\phi^{app}_{(\lmb)}}_{f} .
\end{split}
\end{equation*} 
We first uncover some cancellations between the two terms involving the linearized operator $\calL_{f}$, which resemble those in the proof of Proposition~\ref{prop:KdV-mee}. We write
\begin{align*}
	&-\brk{\phi, \calL_{f}[\phi^{app}_{(\lmb)}]}_{f}
	-\brk{\calL_{f}[\phi], \phi^{app}_{(\lmb)}}_{f} \\
	&= -\int \phi \left(f \phi^{app}_{(\lmb)xxx} + \alp_{1} f_{x} \phi^{app}_{(\lmb) xx} + \alp_{1} f_{xx} \phi^{app}_{(\lmb) x} + f_{xxx} \phi^{app}_{(\lmb)} \right)f^{-\frac{2}{3}\sgm_{c}} \, \ud x \\
	&\relphantom{=} 
	- \mu_{1} \int \phi \left(f^{m-1} \phi^{app}_{(\lmb) x} + (m-1) f^{m-2} f_{x} \phi^{app}_{(\lmb)} \right) f^{-\frac{2}{3}\sgm_{c}} \, \ud x \\
	&\relphantom{=} 
	-\int \left(f \phi_{xxx} + \alp_{1} f_{x} \phi_{xx} + \alp_{1} f_{xx} \phi_{x} + f_{xxx} \phi \right) \phi^{app}_{(\lmb)} f^{-\frac{2}{3}\sgm_{c}} \, \ud x \\
	&\relphantom{=} 
	- \mu_{1} \int \left(f^{m-1} \phi_{x} + (m-1) f^{m-2} f_{x} \phi \right) \phi^{app}_{(\lmb)} f^{-\frac{2}{3}\sgm_{c}} \, \ud x \\
	&=: I + II + III + IV.
\end{align*}
We begin with $I + III$. The zeroth-order terms (in both $\phi$ and $\phi^{app}_{(\lmb)}$ are not dangerous, but we need to perform some integration by parts for the higher order terms. For the third-order terms, we have
\begin{align*}
	&- \int \left( \phi \phi^{app}_{(\lmb) xxx} + \phi_{xxx} \phi^{app}_{(\lmb)} \right)f^{1-\frac{2}{3}\sgm_{c}} \, \ud x \\
	&=  \int \left(\phi_{x} \phi^{app}_{(\lmb) xx} + \phi_{xx} \phi^{app}_{(\lmb) x} \right) f^{1-\frac{2}{3}\sgm_{c}} \, \ud x \\
	&\relphantom{=}
	+ \left(1 - \frac{2}{3} \sgm_{c}\right) \int \left(\phi \phi^{app}_{(\lmb) xx} + \phi_{xx} \phi^{app} \right) f^{-1} f_{x} f^{1-\frac{2}{3}\sgm_{c}} \, \ud x \\
	&=  - 3 \left(1 - \frac{2}{3} \sgm_{c} \right) \int \phi_{x} \phi^{app}_{(\lmb) x} f_{x} f^{-\frac{2}{3}\sgm_{c}} \, \ud x \\
	&\relphantom{=}
	- \left(1 - \frac{2}{3} \sgm_{c} \right) \int \left(\phi \phi^{app}_{(\lmb) x} + \phi_{x} \phi^{app} \right) \left( f_{x} f^{-\frac{2}{3}\sgm_{c}}\right)_{x}  \, \ud x \\
	&=  - 3 \left(1 - \frac{2}{3} \sgm_{c} \right) \int \phi_{x} \phi^{app}_{(\lmb) x} f_{x} f^{-\frac{2}{3}\sgm_{c}} \, \ud x 
	+ \left(1 - \frac{2}{3} \sgm_{c} \right) \int \phi \phi^{app}_{(\lmb)}  \left( f_{x} f^{-\frac{2}{3}\sgm_{c}}\right)_{xx}  \, \ud x;
\end{align*}
for the second-order terms, we have
\begin{align*}
	&- \int \left( \alp_{1} \phi \phi^{app}_{(\lmb) xx} f_{x} + \alp_{1} \phi_{xx} \phi^{app}_{(\lmb)} f_{x} \right)f^{-\frac{2}{3}\sgm_{c}} \, \ud x \\
	&= 2 \alp_{1} \int  \phi_{x} \phi^{app}_{(\lmb) x} f_{x} f^{-\frac{2}{3}\sgm_{c}} \, \ud x 
	+ \alp_{1} \int \left(\phi \phi^{app}_{(\lmb) x} + \phi_{x} \phi^{app}_{(\lmb)}\right)\left(f_{x} f^{-\frac{2}{3}\sgm_{c}} \right)_{x} \\
	&= 2 \alp_{1} \int  \phi_{x} \phi^{app}_{(\lmb) x} f_{x} f^{-\frac{2}{3}\sgm_{c}} \, \ud x 
	- \alp_{1} \int \phi \phi^{app}_{(\lmb)}\left(f_{x} f^{-\frac{2}{3}\sgm_{c}} \right)_{xx};
\end{align*}
and for the first-order terms, we have
\begin{align*}
	- \int \left( \alp_{1} \phi \phi^{app}_{(\lmb) x} f_{xx} + \alp_{1} \phi_{x} \phi^{app}_{(\lmb)} f_{xx} \right)f^{-\frac{2}{3}\sgm_{c}} \, \ud x 
	= \alp_{1} \int \phi \phi^{app}_{(\lmb)} \left( f_{xx} f^{-\frac{2}{3}\sgm_{c}} \right)_{x} \, \ud x.
\end{align*}
In particular, since $3(1- \frac{2}{3} \sgm_{c}) = 2 \alp_{1}$, integrals that involve $\phi_{x} \phi^{app}_{(\lmb)x}$ cancel and we are left with
\begin{equation*}
	\abs{I + III} \leq C \nrm{f}_{Y} \nrm{\phi}_{L^{2}_{f}} \nrm{\phi^{app}_{(\lmb)}}_{L^{2}_{f}}.
\end{equation*}
Next, $II + IV$ consist of first- and zeroth-order terms, where the former may be treated as above and the latter already acceptable. We have
\begin{equation*}
	\abs{II + IV} \leq C \nrm{f}_{C^{0, 1}}^{m-1} \nrm{\phi}_{L^{2}_{f}} \nrm{\phi^{app}_{(\lmb)}}_{L^{2}_{f}}.
\end{equation*}
For the remaining terms in $\frac{\ud}{\ud t} \brk{\phi, \phi^{app}_{(\lmb)}}_{f}$, we have 
\begin{align*}
\nrm{\err[\phi^{app}_{(\lmb)}]}_{L^2_f} &\le C (1+\lmb^{2} t) \abs{\lmb},\\
{\nrm{Q[\phi]}_{L^2_f}} &\le C \nrm{\phi}_{L^2_{f}}, \\
\nrm{f^{-1} \rd_{t} f}_{L^{\infty}} &\le C.
\end{align*}
We conclude that
\begin{equation} \label{eq:KdV-gen-en-conc}
\begin{aligned}
	\left| \frac{\ud}{\ud t} \brk{\phi, \phi^{app}_{(\lmb)}}_{f} \right|
	&\leq C \left(1 + (1+ \lmb^{2} t) \abs{\lmb} \right) \nrm{\phi}_{L^{2}_{f}} \\
	&\leq C \left(1 + (1+ \lmb^{2} t) \abs{\lmb} \right) \exp (C_{1}(M_{5}) t) \frac{\nrm{\phi_{0}}_{L^{2}}}{4 C_{0}}.
\end{aligned}
\end{equation}
Integrating this estimate in time and using \eqref{eq:ill-posed-unbounded-gdkdv-dual0}, we have
 {\begin{equation} \label{eq:ill-posed-unbounded-gdkdv-dual}
	\brk{\phi, \phi^{app}_{(\lmb)}}_{f}(t) \geq \frac{3}{4 C_{0}} \nrm{\phi_{0}}_{L^{2}} \quad \hbox{ for } \abs{t} \leq \min \set{t_{0}, C_{1}(M_{5})^{-1}, c \abs{\lmb}^{-\frac{3}{2}}}.
\end{equation}
for some $c > 0$.}

To proceed, let $m$ be the smallest even integer greater than or equal to $s'$, and define $j = m-s'$. Applying Proposition~\ref{prop:KdV-deg-wp} with $\gmm' = -\sgm_{c}$ and $s' = m$, we have
\begin{align}
	\brk{\phi, \phi^{app}_{(\lmb)}}_{f}
	\leq \int \phi \rd_{x}^{s'+j} \left( \frac{f^{\frac{-2\sgm_{c}+s'+j}{3}}}{(-1)^{s_{0}+j} \lmb^{s_{0}+j}} \phi^{app}_{(\lmb)} \right) \, \ud x + \nrm{\phi}_{L^{2}_{f}} \nrm{\phi_{(\lmb)}^{small}}_{L^{2}_{f}}. \label{eq:ill-posed-unbounded-gdkdv-dual-test}
\end{align}
 {By \eqref{eq:KdV-wp-deg-small} and \eqref{eq:ill-posed-unbounded-gdkdv-mee}, the last term may be bounded as follows for some $C_{2} > 0$:
\begin{align}
\nrm{\phi}_{L^{2}_{f}} \nrm{\phi_{(\lmb)}^{small}}_{L^{2}_{f}}
&\leq C_{2} \lmb^{-1} \exp(C_{1}(M_{5}) t) \frac{1}{4 C_{0}} \nrm{\phi_{0}}_{L^{2}} \notag \\
&\leq \frac{1}{4 C_{0}} \nrm{\phi_{0}}_{L^{2}} \quad \hbox{ if } \abs{t} \leq C_{1}(M_{5})^{-1} \hbox{ and } \abs{\lmb} > C_{2}.  
 \label{eq:ill-posed-unbounded-gdkdv-dual-small}
\end{align}}

We are now ready to conclude the proof of the theorem. For each $\lmb$, there are two possible cases: (i)~$C_{1}(M_{5})^{-1} \leq c \abs{\lmb}^{-\frac{3}{2}}$, or (ii) $C_{1}(M_{5})^{-1} > c \abs{\lmb}^{-\frac{3}{2}}$. In case~(i), we have $c^{-1} \abs{\lmb}^{\frac{3}{2}} \leq C_{1}(M_{5})$, so $M_{5} > (\dlt')^{-\frac{1}{2}}$ if $\abs{\lmb}$ is chosen large enough depending on $C_{1}(\cdot)$ and $\dlt'$. Since $s' \geq s_{c} \geq 5$, the desired norm inflation follows. Hence, it only remains to consider case~(ii). Then, by \eqref{eq:ill-posed-unbounded-gdkdv-dual}, \eqref{eq:ill-posed-unbounded-gdkdv-dual-test} and \eqref{eq:ill-posed-unbounded-gdkdv-dual-small}, we have
\begin{align*}
	\frac{1}{2 C_{0}} \nrm{\phi_{0}}_{L^{2}_{f}}
	\leq \int \phi \rd_{x}^{s'+j} \left( \frac{f^{\frac{-2\sgm_{c}+s'+j}{3}}}{(-1)^{s_{0}+j} \lmb^{s_{0}+j}} \phi^{app}_{(\lmb)} \right) \, \ud x \quad \hbox{ for } \abs{t} \leq \min\set{t_{0}, c \abs{\lmb}^{-\frac{3}{2}}}.
\end{align*}
Using duality and applying \eqref{eq:KdV-wp-deg-estimate}, we arrive at
\begin{align*}
	\frac{1}{2 C_{0}} \nrm{\phi_{0}}_{L^{2}_{f}} 
	& \leq  \nrm{\rd_{x}^{s'} \phi}_{L^{\infty}} \nrm{\rd_{x}^{j} \left( \lmb^{-s_{0}-j} f^{\frac{-2\sgm_{c}+s'+j}{3}} \phi^{app}_{(\lmb)}\right)}_{L^{1}}  \\
	& \leq C \lmb^{-s'} e^{-3 \bt(t) \lmb^{2} (-\sgm_{c}+s'+\frac{1}{2}) t} \nrm{\rd_{x}^{s'} \phi}_{L^{\infty}} .
\end{align*}
Rearranging the factors, we finally arrive at
\begin{equation*}
	\nrm{\rd_{x}^{s'} \phi(t)}_{L^{\infty}} \geq \frac{1}{C} \lmb^{s'-m_{0}} \exp \left(3 \bt(t) \lmb^{2} (-\sgm_{c} + s' + \tfrac{1}{2}) t \right) \quad \hbox{ for } 0 \leq t \leq \min \set{t_{0}, c \abs{\lmb}^{-\frac{3}{2}}}.
\end{equation*}
Fix $t = c \abs{\lmb}^{-\frac{3}{2}}$; by taking $\lmb$ sufficiently large, we may clearly ensure that $t \leq t_{0}$. Since $s' \geq s_{c} > \sgm_{c} - \frac{1}{2}$ and $\lmb^{2} t = c \abs{\lmb}^{\frac{1}{2}}$, by taking $\lmb$ larger we may also ensure that the right-hand side is at least $(\dlt')^{-\frac{1}{2}}$ for each $s'$. This completes the proof. \hfill \qedsymbol

\subsection{Proof of Theorem \ref{thm:ill-posed-nonexist-gdkdv}}\label{subsec:Kmn-nonexist} \label{subsec:nonexist-gdkdv}

We are now in a position to complete the proof of Theorem \ref{thm:ill-posed-nonexist-gdkdv}. As we shall see, the argument is parallel to that for Theorem \ref{thm:ill-posed-nonexist-gDS} in the Schr\"odinger case. 

Let $s_{0} \geq s \geq s_{c}$ and $\eps > 0$ be given as in the statement of Theorem~\ref{thm:ill-posed-nonexist-gdkdv}. Suppose, for contradiction, that for every $u_{0} \in C^{\infty}(\bbT)$ satisfying $\nrm{u_{0}}_{C^{s_{0}}} < \eps$, there exists $\dlt = \dlt(u_{0}) > 0$ and a corresponding solution $u$ to \eqref{eq:gdkdv} belonging to $L^{\infty} ([0, \dlt]; C^{s}(\bbT))$. 

We shall fix a function $\mathring{f}_{0} \in C^{\infty}(\bbT)$ supported in $[-4x_{1},4x_{1}]$ which satisfies $\mathring{f}_{0}(x)=x^{3}$ in $[-x_{1},x_{1}]$, $\mathring{f}_{0}>0$ on $(0,2x_{1})$ and cubic vanishing property at $2x_{1}$ for some small $0<x_{1}<\frac{1}{100}$; in what follows, we omit the dependence of constants on $\mathring{f}_{0}$. Then, we take \begin{equation*}
	\begin{split}
		f_{0} := \sum_{k=k_{0}}^{\infty} f_{k,0} := \sum_{k=k_{0}}^{\infty} A_{k} 2^{-3k} \mathring{f}_{0}(2^{k}(x-x_{k})), \quad x_{k} = 2^{-\frac{k}{2}}, \quad A_{k} = 2^{-k^{2}}, 
	\end{split}
\end{equation*} where $k_{0}$ shall be fixed below. We see that $f_{0} \in C^{\infty}$ and $\nrm{f_{0}}_{C^{s_{0}}}<\frac{\eps}{2}$ provided that $k_{0}$ is sufficiently large depending on $s_{0}$ and $\eps$. Moreover, we can check that there exists a constant $C_{0} = C_{0}(f_{0})>0$ such that \eqref{eq:Km2-initial-condition} holds for $f_{0}$. For each $k \in \bbN$, we take cutoff functions $\chi_{k}$ which is equal to 1 on the support of $f_{k,0}$  and vanishes on the support of $f_{k',0}$ for all $k'\ne k$. 

Now, let $ {f} \in L^\infty([0,\dlt];C^{s}(\bbT))$ be a solution to \eqref{eq:gdkdv} with $ {f}(t=0)= {f}_{0}$. Using the equation and $C^{4,1}$--regularity, we have the pointwise estimate $|\rd_{t}  {f}| \lesssim |f_{0}| + | {f}|$, which guarantees that the support of $ {f}(t,\cdot)$ is preserved in time. For simplicity, we set $M_{f} = \nrm{f}_{ L^\infty([0,\dlt];C^{4,1})}$ and replace $\dlt$ with $\min\{ \dlt, c \}$ with some small constant $c = c(M_{f})>0$, so that we have uniformly \begin{equation}\label{eq:kdv-uniform-f}
	\begin{split}
		\frac12 f_{0}(x) \le f(t,x) \le \frac32 f_{0}(x), \qquad t \in [0,\dlt]
	\end{split}
\end{equation} whenever $f_{0}(x) > 0$. The existence of such a constant $c$ follows from \eqref{eq:f-f_0-pointwise} and \eqref{eq:f-f_0-pointwise2}. 

Since $\chi_{k}$ is either 0 or 1 on $\supp(f_{0}) = \supp( f(t,\cdot) )$, we have $\rd_{x}\chi_{k} \equiv 0$ on the support of $f$ and $\chi_{k}f = \chi_{k}^{2}f$. From these observations, it follows that $\chi_{k}f =: f_{k}$ provides a solution to \eqref{eq:gdkdv} for any $k\ge k_{0}$ with initial data $f_{k,0}$. Let $I_{k} = [x_{k}, x_{k} + 2^{-k} 2 x_{1}]$, and observe that 
\begin{equation} \label{eq:KdV-fk0-Y}
\nrm{f_{k, 0}}_{Y(I_{k})} \leq A_{k} (\nrm{\mathring{f}_{0}}_{Y([0, 2 x_{1}])} + 2^{-3k} \nrm{\mathring{f}_{0}}_{L^{\infty}([0, 2 x_{1}]}) \to 0 \hbox{ as } k \to \infty.  
\end{equation}
Taking $k_{0}$ larger and $\dlt$ smaller if necessary, by Proposition~\ref{prop:C41}, we have, for every $k \geq k_{0}$,
\begin{equation*}
\nrm{f_{k}}_{L^{\infty}([0, \dlt]; Y(I_{k})} \leq C(M_{f}).
\end{equation*}

Let us fix a $C^\infty$--smooth profile $g_{0}$ supported in $(\frac{x_1}{2},x_{1})$ and normalized in $L^{2}$; in what follows, we omit the dependence of constants on $g_{0}$. We take $g_{k}(x) = 2^{\frac{k}{2}} g_{0}( 2^{k}(x-x_{k}) )$. For a strictly increasing sequence $\{ \lmb_{k} \}_{k\ge k_{0}}$ ($\lmb_{k}\gg1$) to be determined, we consider the wave packets \begin{equation*}
	\begin{split}
		\phi_{k}^{app} (t,x) := \phi^{app}_{(\lmb_{k})}[g_{k}, f_{k}], 
	\end{split}
\end{equation*} where $\phi^{app}_{(\lmb_{k})}[g_{k}, f_{k}]$ denotes the wave packet constructed in Proposition \ref{prop:KdV-deg-wp} using the solution $f_{k}$ with profile $g_{k}$ and frequency $\lmb_{k}$. We define the corresponding error by \begin{equation}\label{eq:kdv-nonexist-error}
	\begin{split}
		[\rd_t + \calL_{f_{k}}] \phi^{app}_k = \err_{\phi_{k}}
	\end{split}
\end{equation} where $\calL_{f_{k}}$ is simply \eqref{eq:gdkdv-Lf} with $f$ replaced with $f_{k}$. Recall also the definition of $L^2_{f}$ norm from \eqref{eq:f-gdkdv}, and observe that $f = f_{k}$ on the support of $\phi^{app}_{k}$. Applying Proposition \ref{prop:KdV-deg-wp}, we obtain the following properties of $\phi^{app}_{k}$ for all $t \in [0, \dlt]$:  \begin{itemize}
	\item $\nrm{\phi^{app}_k(t,x)}_{L^2_{f}} \le C_{\tld{f}_{k}} \nrm{ \phi^{app}_{k,0}}_{L^2_{f_{k}}} \le C_{\tld{f}_{k}} \nrm{g_{0}}_{L^{2}},$
	\item $\nrm{(\tld{f}^{\frac{1}{3}} \rd_{x})^{n}\phi^{app}_k(t,x)}_{L^2_{f}} \le C_{\tld{f}_{k}} \nrm{ \phi^{app}_{k,0}}_{L^2_{f_{k}}} \le C_{\tld{f}_{k}},$
	\item $\nrm{\err_{\phi_{k}}(t,x)}_{L^2_{f}} \le C_{\tld{f}_{k}} (1+A_{k}^{m-2})\lmb_{k} (1 + \lmb_{k}^{2} t)$,
\end{itemize} and, since $s$ is even, we have
\begin{equation}
	f^{-\frac{2}{3}\sgm_{c}} \phi^{app}_{k} = \rd_{x}^{s}\left( \frac{f^{\frac{-2\sgm_{c}+s}{3}}}{(-1)^{\frac{s}{2}} A_{k}^{\frac{s}{3}}\lmb_{k}^{s}}  \phi^{app}_{k}\right) + f^{-\frac{2}{3}\sgm_{c}} \phi^{small}_{k},
\end{equation} with \begin{itemize}
\item $\left\| \left(A_{k}^{-\frac{s}{3}}\lmb^{-s} f^{\frac{-2\sgm_{c}+s}{3}} \phi^{app}_{k}\right) (t, \cdot ) \right\|_{L^{1}} \leq C_{\tld{f}_{k}} A_{k}^{-\frac{\sgm_{c}}{3}} \lmb_{k}^{-s} e^{-3 \bt_{k}(t) A_{k}^{\frac{2}{3}} \lmb_{k}^{2}(-\sgm_{c}+s+\frac{1}{2}) t},$
\item $\nrm*{f^{-\frac{\sgm_{c}}{3}} \tphi^{small}_{k}(t, \cdot )}_{L^{2}} \leq C_{\tld{f}_{k}} \lmb_{k}^{-1}$,
\end{itemize}
where $\bt_{k}(t)$ is the solution to \eqref{eq:KdV-beta-ODE} with $f$ replaced by $f_{k}$. Observe that $C_{\tld{f}_{k}}$ depends on $M_{f}$, $k$ and $s$, but \emph{not} on $\lmb_{k}$. In view of \eqref{eq:KdV-deg-wp-id-L2}, note that $\lmb_{k}$ should be sufficiently large depending on $k$ and $g_{0}$ to ensure that the first inequality in the first item holds. Define \begin{equation*}
	\begin{split}
		\phi_0(x) := \sum_{k=k_{0}}^{\infty} \phi_{k,0}(x) := \sum_{k=k_{0}}^{\infty} \exp(- \lmb_k^{\frac18}) \phi_{k, 0}^{app}(x).
	\end{split}
\end{equation*} By ensuring some growth of $\lmb_{k}$ (e.g., $\lmb_{k} \geq 2^{k}$) and by taking $k_{0}$ even larger if necessary, we can guarantee that $\nrm{\phi_0}_{C^{s}}<\frac{\eps}{2}$, so that $u_0 := f_0 + \phi_0,$  is $C^{\infty}$--smooth and satisfies $\nrm{u_{0}}_{C^{s_{0}}}<\eps$. From the contradiction hypothesis, we have a $L^\infty_tC^{s}$--solution $u(t,x)$ to \eqref{eq:gdkdv} with initial data $u_0$ on some time interval $[0,\dlt']$. By shrinking either $\dlt$ or $\dlt'$, we may assume that $0<\dlt' = \dlt$.
We now set \begin{equation*}
	\begin{split}
		\phi(t) := u(t) - f(t), \quad \phi_k(t) := \chi_k\phi(t)
	\end{split}
\end{equation*} for all $k\ge k_{0}$. Moreover, taking $k_{0}$ larger if necessary, we can easily arrange that $f_{0}(x)$ and $u_{0}(x)$ are uniformly comparable; for all $x$, \begin{equation*}
\begin{split}
	\frac{7}{8} u_{0}(x) \le f_{0}(x) \le \frac98 u_{0}(x). 
\end{split}
\end{equation*} 
From the conservation of the support in time, we have that $\sum_{k=k_{0}}^{\infty} \phi_k = \phi$. We now introduce
\begin{equation} \label{eq:uniform}
	\begin{split}
		M = 1+ \sup_{t \in [0,\dlt]} ( \nrm{f(t)}_{C^{s}} + \nrm{\phi(t)}_{C^{s}} ) 
		\ge M_{f},
	\end{split}
\end{equation} which is finite by the contradiction hypothesis, and further replace $\dlt$ with $\min\{ \dlt, c \}$, where $c=c(M)>0$ is large so that we have \begin{equation}\label{eq:kdv-uniform-u}
\begin{split}
	\frac12 u_{0}(x) \le u(t,x) \le \frac32 u_{0}(x), \qquad t\in[0,\dlt], 
\end{split}
\end{equation} whenever $u_{0}(x) > 0$. Note also that, by a computation similar to \eqref{eq:KdV-fk0-Y}, we have $\nrm{u_{0}}_{Y(I_{k})} \to 0$ as $k \to \infty$. Choosing $k_{0}$ sufficiently large and $\dlt$ small enough, by Proposition~\ref{prop:C41}, we have,  for every $k \geq k_{0}$,
\begin{equation*}
\nrm{u_{k}}_{L^{\infty}([0, \dlt]; Y(I_{k}))} \leq C(M).
\end{equation*}
We see that $\chi_{k} u$ is a solution to \eqref{eq:gdkdv} and it follows that $\phi_{k}$ solves \begin{equation*}
\begin{split}
	[\rd_{t} + \calL_{f_{k}}]\phi_{k} = Q_{f_{k}}[\phi_{k}], 
\end{split}
\end{equation*}  where $Q_{f_{k}}$ is simply \eqref{eq:gdkdv-Q} with $f$ replaced with $f_{k}$. 

From Proposition~\ref{prop:KdV-mee}, we have the modified energy estimate for $\phi_{k}$:  \begin{equation*}
\begin{split}
	\nrm{\phi_k(t)}_{L^2_{f}(I_{k})} \le C(M) \nrm{\phi_{k,0}}_{L^2_{f_0}(I_{k})}. 
\end{split}
\end{equation*} Proceeding as in the proof of \eqref{eq:KdV-gen-en-conc}, we obtain \begin{equation}\label{eq:gdkdv-gee}
	\begin{split}
		\left|\frac{\ud}{\ud t} \brk{\phi_{k}, \phi^{app}_{k}}_{f}  \right| \le C(M, k, s) \left( 1 + \lmb_k (1 + \lmb_{k}^2 t) \right)\nrm{\phi_{k,0}}_{L^2_{f_{0}}(I_{k})}, \qquad t \in [0,\dlt].
	\end{split}
\end{equation} 
We shall take  $t_{k} := \lmb_{k}^{-\frac53} \ll 1$, and make sure that $k_{0}$ is large enough so that $t_{k} \le \dlt$. Integrating \eqref{eq:gdkdv-gee} from $t = 0$ to $t_{k}$, \begin{equation*}
\begin{split}
	\brk{\phi_k, \phi^{app}_k}_{f} (t_{k}) \ge \left( 1 - C(M, k, s) \left( 1 + \lmb_k (1 + \lmb_{k}^2 t_{k}) \right) t_{k} \right)\nrm{\phi_{k,0}}_{L^2_{f_{0}}}  \ge  \frac12\nrm{\phi_{k,0}}_{L^2_{f_{0}}},
\end{split}
\end{equation*}  by taking $\lmb_{k}$ larger if necessary. On the other hand, we write \begin{equation*}
\begin{split}
	\brk{\phi_k, \phi^{app}_k}_{f} (t_{k})& = \frac{1}{(-1)^{\frac{s}{2}} A_{k}^{\frac{s}{3}} \lmb_{k}^{s}} \brk{ \phi_{k},  \rd_{x}^{s}\left( f^{\frac{-2\sgm_{c}+s}{3}}  \phi^{app}_{k}\right)  }(t_{k}) +  \brk{ f^{-\frac{\sgm_{c}}{3}} \phi_{k},  f^{-\frac{\sgm_{c}}{3}} \phi^{small}_{k} }(t_{k}).
\end{split}
\end{equation*} 
Using the above estimates for $\phi_{k}$ and $\phi^{small}_{k}$ at $t = t_{k}$, for $\lmb_{k}$ sufficiently large, the last term on the right-hand side is bounded by $\frac{1}{4} \nrm{\phi_{k, 0}}_{L^{2}_{f_{0}}}$. We may therefore obtain
\begin{equation*}
\begin{split}
	\frac14\nrm{\phi_{k,0}}_{L^2_{f_{0}}} &\le A_{k}^{-\frac{s}{3}} \lmb_{k}^{-s} \nrm{\rd_{x}^{s} \phi_{k}(t_{k})}_{L^{\infty}} \nrm{f^{\frac{-2\sgm_{c} + s}{3}} \phi_{k}^{app}(t_{k})}_{L^{1}} \\
	&\le C(M, k, s) A_{k}^{-\frac{\sgm_{c}}{3}} \lmb_{k}^{-s} \exp \left(- 3 \bt_{k}(t) A_{k}^{\frac{2}{3}} (-\sgm_{c} + s + \tfrac{1}{2}) \lmb_{k}^{2-\frac{5}{3}} \right) \nrm{\rd_{x}^{s} \phi_{k}(t_{k})}_{L^{\infty}}
\end{split}
\end{equation*} 
Recalling that $\nrm{\phi_{k, 0}}_{L^{2}_{f_{0}}} \geq c(k) e^{-\lmb_{k}^{\frac{1}{8}}}$, we arrive at the lower bound \begin{equation*}
\begin{split}
	\nrm{\phi_{k}(t_{k})}_{C^{s}} \ge c(M, k, s) \lmb_{k}^{s} \exp \left(3 \bt_{k}(t) A_{k}^{\frac{2}{3}} (-\sgm_{c} + s + \tfrac{1}{2}) \lmb_{k}^{\frac{1}{3}} - \lmb_{k}^{-\frac{1}{8}}\right).
\end{split}
\end{equation*}
Finally choosing $\lmb_{k}$ to be sufficiently large, we may guarantees that 
\begin{equation*}
M \geq \sup_{t \in [0, \dlt'] } \nrm{\phi(t)}_{C^{s}} \geq \nrm{\phi_{k}(t_{k})}_{C^{s}} \to \infty \hbox{ as } k \to \infty.
\end{equation*}
This contradicts the finiteness of $M$ in \eqref{eq:uniform}, which completes the proof.\hfill \qedsymbol

\appendix
\section{Takeuchi--Mizohata illposedness via duality} \label{sec:Mz}
In this appendix, we show how application of the duality (or generalized energy) argument from \cite{JO1} and this paper leads to simple proofs of quantitative illposedness results for first-order perturbations of the free Schr\"odinger equation related to the Takeuchi--Mizohata condition, including Proposition~\ref{prop:Mz-d-norm-infl} (see Appendix~\ref{subsec:Mz-d}). 

\subsection{One-dimensional case} \label{subsec:Mz-1}
We begin with the one-dimensional first-order perturbation of the free Schr\"odinger equation,
\begin{equation} \label{eq:Mz-1}
i \rd_{t} u + \rd_{xx} u + b(x) \rd_{x} u = 0.
\end{equation}
Fix $x_{0} \in \bbR$. For $T > 0$ and $\mu \geq 1$, we define the weight $$w(x) = \exp\left(\int_{0}^{x} \Re \frac{b(x')}{2} \, \ud x' \right)$$ and the growth factor 
\begin{equation*}
	M(T, \mu) = \inf_{(y, y_{0}) : \abs{y} \leq \mu^{-1}, \abs{y_{0}} \leq \mu^{-1}}\exp\left( \int_{x_{0} + y_{0}}^{x_{0} + 2 T + y} \Re \frac{b(x')}{2} \, \ud x' \right).
\end{equation*}
Fix also $\psi_{1} \in C^{\infty}(\bbR)$ with $\supp \psi_{1} \subseteq \set{x : \abs{x} < 1}$ and $\nrm{\psi_{1}}_{L^{2}} = 1$. Given $\mu \geq 1$ (inverse spatial scale), define $\psi_{\mu, x_{0}} = \mu^{\frac{d}{2}} \psi_{1}(\mu (x-x_{0}))$. Given also $\lmb \geq 1$ (frequency, or inverse semi-classical parameter) -- which in practice would be much larger than $\mu$ -- define
\begin{equation}\label{eq:Mz-1D-wp}
\tld{u} (t, x) = w^{-1}(x) e^{i \lmb x - i \lmb^{2} t} \exp \left( - \int_{0}^{\lmb t} i \Im b(x - 2s) \, \ud s \right) \psi_{\mu, x_{0}}(x - 2 \lmb t).
\end{equation}
This is a wave packet that approximately solves \eqref{eq:Mz-1}; see \eqref{eq:Mz-1-wp}--\eqref{eq:Mz-1-wp-err} below.

\begin{proposition} \label{prop:Mz-1}
Let $\tld{u}$ be as in \eqref{eq:Mz-1D-wp} and let $u_{0}$ satisfy 
\begin{equation*}
\int \Re(u_{0} \br{\tld{u}(0)}) w^{2} \, \ud x = 1, \quad \supp u_{0} \subseteq [x_{0} - \mu^{-1}, x_{0} + \mu^{-1}].
\end{equation*}
Then there exists at least one corresponding solution $u$ of \eqref{eq:Mz-1} belonging to $L^{\infty}_{loc,t}(\bbR; L^{2}_{w})$. Assume that it furthermore satisfies $u \in L^{\infty}_{t}([0, t_{f}]; L^{2})$ with 
\begin{equation} \label{eq:Mz-1-hyp}
	t_{f} \leq c \mu^{-1},
\end{equation}
where $c$ is a constant depending only on $\nrm{b}_{C^{1, 1}}$ and $\nrm{\psi_{1}}_{H^{2}}$. Then, $u(t)$ necessarily satisfies the pointwise lower bound
\begin{equation} \label{eq:Mz-1-lower}
	\nrm{u(t)}_{L^{2}} \geq \frac{1}{2} M(\lmb t, \mu) \nrm{u_{0}}_{L^{2}} \quad \hbox{ for all } 0 \leq t \leq t_{f}.
\end{equation}
\end{proposition}

\begin{proof}
As discussed in Section~\ref{subsec:discussion-pf}, we consider the conjugation $v = w u$. Introducing the formally self-adjoint operator
\begin{equation*}
	\tld{\calL} = \lap + i \Im \, b(x) \rd_{x} + \frac{i}{2}  \Im b_{x},
\end{equation*}
we have the conjugation identity
\begin{equation*}
	(i \rd_{t} + \tld{\calL}) v = \left( w^{-1} \rd_{x}^{2} w + b w^{-1} \rd_{x} w + \frac{i}{2} \Im b_{x} \right) v + w (i \rd_{t} + \lap + b(x) \rd_{x} ) u.
\end{equation*} Under the assumption that $v(t) \in L^{2}$, we have for $t\ge0$
\begin{equation} \label{eq:Mz-1-v}
	\nrm{v(t)}_{L^{2}} \leq e^{C_{0} t} \nrm{v_{0}}_{L^{2}},
\end{equation}
where $C_{0}$ depends only on $\nrm{b}_{C^{1, 1}}$ and $\nrm{\psi_{1}}_{H^{2}}$. Moreover, we also have
\begin{equation} \label{eq:Mz-1-v-err}
	\nrm{(i \rd_{t} + \tld{\calL}) v(t)}_{L^{2}} \leq C_{0} e^{C_{0} t} \nrm{v_{0}}_{L^{2}},
\end{equation}
where $C_{0}$ depends only on $\nrm{b}_{C^{1, 1}}$ and $\nrm{\psi_{1}}_{H^{2}}$. 

Next, we consider the standard wave packet $\tld{v}$ for $i \rd_{t} + \tld{\calL}$ obtained by solving \eqref{eq:HJ-T-sch} with $a = 1$, $\Re b = 0$, $\bfPhi(0, x) = \lmb x$ and $\bfa(0, x) = \psi_{\mu, x_{0}}(x)$. It is given explicitly by
\begin{equation*}
	\tld{v} = e^{i \lmb x - i \lmb^{2} t} \exp \left( - \int_{0}^{\lmb t} i \Im \, b(x - 2s) \, \ud s \right) \psi_{\mu, x_{0}}(x - 2 \lmb t).
\end{equation*}
(Note, furthermore that $\tld{u} = w^{-1} \tld{v}$.) By the definition, we clearly have, for all $t$, 
\begin{equation} \label{eq:Mz-1-wp}
	\nrm{\tld{v}(t)}_{L^{2}} = 1.
\end{equation}
Moreover, by the support property of $\tld{v}$, it follows that
\begin{equation} \label{eq:Mz-1-wp-w}
	\nrm{w \tld{v}(t)}_{L^{2}} \leq \sup_{y : \abs{y} \leq \mu^{-1}} w(x - 2 \lmb t + y).
\end{equation}
Finally, we consider the error incurred by $\tld{v}$. A straightforward computation gives the following:
\begin{lemma}\label{lem:error} 
	We have 
	\begin{equation} \label{eq:Mz-1-wp-err}
		\nrm{(i \rd_{t} + \tld{\calL}) \tld{v}}_{L^{2}} \leq \tld{C}_{0} \mu^{2}.
	\end{equation}
	where $\tld{C}_{0}$ depends only on $\nrm{b}_{C^{1, 1}}$ and $\nrm{\psi_{1}}_{H^{2}}$. 
\end{lemma}
Given the above lemma, by the self-adjointness of $\tld{\calL}$, we have
\begin{align*}
\frac{\ud}{\ud t} \brk{v, \tld{v}} = - \brk{(i \rd_{t} + \tld{\calL}) v, i \tld{v}} + \brk{iv, i (i \rd_{t} + \tld{\calL}) \tld{v}}.
\end{align*}
By \eqref{eq:Mz-1-v}, \eqref{eq:Mz-1-v-err}, \eqref{eq:Mz-1-wp} and \eqref{eq:Mz-1-wp-err}, we have
\begin{align*}
\abs*{\frac{\ud}{\ud t} \brk{v, \tld{v}}} \leq (C_{0}  + \tld{C}_{0} \mu^{2}) e^{C_{0} t} \nrm{v_{0}}_{L^{2}}.
\end{align*}
Provided that \eqref{eq:Mz-1-hyp} holds with $c$ sufficiently small compared to $C_{0}$ and $\tld{C}_{0}$, we see that
\begin{equation*}
\brk{v, \tld{v}}(t) \geq \frac{1}{2} \nrm{w u_{0}}_{L^{2}} \quad \hbox{ for all } 0 \leq t \leq t_{f}.
\end{equation*}
On the one hand, by duality (i.e., Cauchy--Schwartz) and \eqref{eq:Mz-1-wp-w},
\begin{align*}
\brk{v, \tld{v}}(t) = \brk{u, w \tld{v}}(t) \leq \sup_{y : \abs{y} \leq \mu^{-1}} w(x - 2 \lmb t + y) \nrm{u(t)}_{L^{2}} 
\end{align*}
On the other hand, by the support property of $u_{0}$, we have
\begin{align*}
\nrm{u_{0}}_{L^{2}} \leq \sup_{y : \abs{y} \leq \mu^{-1}} w^{-1}(x_{0} - 2 \lmb t + y) \nrm{w u_{0}}_{L^{2}}.
\end{align*}
Combining the preceding three inequalities, we arrive at \eqref{eq:Mz-1-lower}. \qedhere
\end{proof}

\begin{proof}[Proof of Lemma \ref{lem:error}]
	We compute, with $\psi = \psi_{\mu, x_{0}}$, \begin{equation*}
		\begin{split}
			i \rd_{t}  \tld{v} = \lmb^{2} \tld{v} + \lmb \tilde{v}  \Im\, b(x-2\lmb t) - 2i\lmb \tld{v} \frac{\psi_{x}(x-2\lmb t)}{\psi(x-2\lmb t)} , \quad 		\rd_{x} \tld{v} = \left( i\lmb - \int_{0}^{\lmb t} i \Im \, b_{x} (x-2s) \, \ud s +  \frac{\psi_{x}(x-2\lmb t)}{\psi(x-2\lmb t)} \right) \tld{v}  ,
		\end{split}
	\end{equation*} and \begin{equation*}
\begin{split}
	\rd_{xx} \tld{v} &= -\left(  \lmb -  \int_{0}^{\lmb t}   \Im \, b_{x} (x-2s) \, \ud s \right)^{2} \tld{v} + 2\left( i\lmb - \int_{0}^{\lmb t} i \Im \, b_{x} (x-2s) \, \ud s \right) \tld{v} \frac{\psi_{x}(x-2\lmb t)}{\psi(x-2\lmb t)} \\
	&\quad  -  \tld{v} \int_{0}^{\lmb t} i \Im \, b_{xx} (x-2s) \, \ud s + \tld{v} \frac{\psi_{xx}(x-2\lmb t)}{\psi(x-2\lmb t)} . 
\end{split}
\end{equation*}
Then, after several direct cancellations, we have  \begin{equation*}
	\begin{split}
		(i \rd_{t} + \tld{\calL}) \tld{v} & =  \lmb \left(  \Im\, b(x-2\lmb t)   -   \Im \, b(x)   + 2 \int_{0}^{\lmb t}  \Im \, b_{x} (x-2s) \, \ud s\right) \tld{v} \\
		&\quad - \Im \, b \left( -\int_{0}^{\lmb t} i \Im \, b_{x} (x-2s) \, \ud s  +   \frac{\psi_{x}(x-2\lmb t)}{\psi(x-2\lmb t)} \right)\tld{v} + \left( \int_{0}^{\lmb t} i \Im \, b_{x} (x-2s) \, \ud s  \right)^{2} \tld{v} \\
		&\quad  +2\left( - \int_{0}^{\lmb t} i \Im \, b_{x} (x-2s) \, \ud s \right) \tld{v} \frac{\psi_{x}(x-2\lmb t)}{\psi(x-2\lmb t)} \\
		&\quad + \left(  - \int_{0}^{\lmb t} i \Im \, b_{xx} (x-2s) \, \ud s\right) \tld{v} +  \frac{\psi_{xx}(x-2\lmb t)}{\psi(x-2\lmb t)} \tld{v} + \frac{i}{2} \Im \, b_{x} \,  \tld{v}. 
	\end{split}
\end{equation*} Using \begin{equation*}
\begin{split}
	\int_{0}^{\lmb t}   \Im \, b_{x} (x-2s) \, \ud s = -\frac12 \Im  \left( b(x-2\lmb t) - b(x) \right)
\end{split}
\end{equation*} we get a cancellation of remaining terms of order $\lmb$. Moreover, the same identity eliminates all integrals on the domain $[0, \lmb t]$. Using $\mu \geq 1$ to bound all the other terms by $O(\mu^{2})$, the proof is complete. \end{proof}

\subsection{Multi-dimensional case} \label{subsec:Mz-d}
We consider the following equation on $\bbR^{d}$:
\begin{equation} \label{eq:Mz-d}
i \rd_{t} u + \lap u + b^{j}(x) \rd_{j} u = 0.
\end{equation}

Unfortunately, the proof of a pointwise lower bound in Proposition~\ref{prop:Mz-1} breaks down due to the lack of a simple physical space conjugation that removes $\Re b^{j} (x) \rd_{j}$. Instead, we shall prove two (conceptually) weaker statements using the duality method, including Proposition~\ref{prop:Mz-d-norm-infl}.

The first result is an unconditional \emph{integrated} lower bound that is valid for nontrivally long (i.e., $t \gg \lmb^{-1}$) timescales. To state this result, given $x_{0} \in \bbR^{d}$, $\omg_{0} \in \bbS^{d-1}$, $\mu \geq 1$ and $T > 0$, define
\begin{equation*}
	M_{x_{0}, \omg_{0}}(T, \mu) = \inf_{y : \abs{y} \leq \mu^{-1}} \exp\left( - \int_{0}^{T} \Re b^{j}(x_{0} + y- 2 s \omg_{0}) (\omg_{0})_{j} \, \ud s\right).
\end{equation*}
Fix $\psi_{1} \in C^{\infty}(\bbR^{d})$ with $\supp \psi_{1} \subseteq \set{x : \abs{x} < 1}$ and $\nrm{\psi_{1}}_{L^{2}} = 1$. Given $\mu \geq 1$, define $\psi_{\mu, x_{0}} = \mu^{\frac{d}{2}} \psi_{1}(\mu (x-x_{0}))$. Given also $\lmb \geq 1$, define (cf.~\cite[\S VII.2]{Mz2})
\begin{equation*}
\tld{u} (t, x) = e^{i \lmb \omg_{0} \cdot x - i \lmb^{2} t} \exp \left( - \int_{0}^{\lmb t} b^{j}(x - 2 s \omg_{0}) (\omg_{0})_{j} \ud s \right) \psi_{\mu, x_{0}}(x - 2 \lmb \omg_{0} t).
\end{equation*}
\begin{proposition} \label{prop:Mz-d}
Let $u \in L^{\infty}_{t}([0, t_{f}]; L^{2})$ be a solution to \eqref{eq:Mz-d} with initial data $u_{0}$ satisfying $\brk{u_{0}, \tld{u}(0)} = 1$, where $\tld{u}$ is determined from $\mu$, $\nu$ and $\psi_{1}$ as above. Then as long as
\begin{equation} \label{eq:Mz-d-hyp}
	\mu \leq c \lmb^{\frac{1}{3}}, \quad t_{f} \leq c \lmb^{-\frac{2}{3}},
\end{equation}
where $c$ is a constant depending only on $\nrm{b}_{C^{1, 1}}$ and $\nrm{\psi_{1}}_{H^{2}}$, $u(t)$ necessarily satisfies the averaged lower bound
\begin{equation} \label{eq:Mz-d-lower}
	\frac{1}{t_{f}} \int_{0}^{t_{f}} \frac{(1 + \mu^{-1} \lmb t)^{2}}{(1 + \mu^{-1} \lmb t_{f})^{2}} M_{x_{0}, \omg_{0}}(\lmb t, \mu)^{-1} \nrm{u(t)}_{L^{2}} \, \ud t \geq \frac{1}{6} \nrm{u_{0}}_{L^{2}}.
\end{equation}
\end{proposition}
An example of an initial data $u_{0}$ satisfying the above hypothesis is, of course, $u_{0} = \tld{u}(0)$, in which case $u$ is expected to behave like $\tld{u}$. An argument similar to the proof of \eqref{eq:Mz-d-v-L2} shows that
\begin{equation*}
	\nrm{\tld{u}(t)}_{L^{2}} \leq \sup_{y : \abs{y} \leq \mu^{-1}} \exp\left( - \int_{0}^{\lmb t} \Re b^{j}(x_{0} + y- 2 s \omg_{0}) (\omg_{0})_{j} \, \ud s\right) \nrm{\tld{u}(0)}_{L^{2}}.
\end{equation*}
Thence, provided we choose $\mu^{-1}$ to be sufficiently small depending on $b$, \eqref{eq:Mz-d-lower} is sharp for $\tld{u}$ up to a constant.

\begin{proof}
We introduce $\calL$ and its formal $L^{2}$-adjoint $\calL^{\ast}$ (in operator notation),
\begin{equation*}
	\calL = \lap + b^{j}(x) \rd_{j}, \quad \calL^{\ast} = \lap - \rd_{j} \br{b}^{j}(x).
\end{equation*}
The basis of the proof of Proposition~\ref{prop:Mz-d} is the following generalized energy identity:
\begin{equation} \label{eq:Mz-d-gen-en}
\frac{\ud}{\ud t} \brk{u_{1},u_{2}} = - \brk{(i \rd_{t} + \calL) u_{1}, i u_{2}} - \brk{i u_{1}, (i \rd_{t} + \calL^{\ast}) u_{2}},
\end{equation}
which is a consequence of the Leibniz rule for $\rd_{t}$ and $0 = \brk{i \calL u_{1}, u_{2}} + \brk{i  u_{1}, i \calL^{\ast} u_{2}}$.
A simple but important observation is that \eqref{eq:Mz-d-gen-en} holds even under the weak assumption $u_{1} = u \in L^{\infty}_{t}([0, t_{f}]; L^{2})$, provided that $u_{2}$ is nice enough, e.g., smooth in $t, x$ and compactly supported in space for each fixed time.

The identity \eqref{eq:Mz-d-gen-en} motivates us to consider \emph{not} a wave packet for $i \rd_{t} + \calL$, but rather its adjoint $i \rd_{t} + \calL^{\ast}$. Given $1 \leq \mu \leq \lmb$, consider 
\begin{equation*}
\tld{u}^{\ast} (t, x) = e^{i \lmb \omg_{0} \cdot x - i \lmb^{2} t} \exp \left( \int_{0}^{\lmb t} \br{b}^{j}(x - 2 s \omg_{0}) (\omg_{0})_{j} \ud s \right) \psi_{\mu, x_{0}}(x - 2 \lmb \omg_{0} t).
\end{equation*}
Observe that
\begin{equation*}
	\exp \left( \int_{0}^{\lmb t} \Re \br{b}^{j}(x - 2 s \omg_{0}) (\omg_{0})_{j} \ud s \right) \leq M(\lmb t, \mu)^{-1} \quad \hbox{ for } x \in \supp \psi_{\mu, x_{0}}(\cdot - 2 \lmb \omg_{0} t),
\end{equation*}
where we have introduced the abbreviation $M(\lmb t, \mu) = M_{x_{0}, \omg_{0}}(\lmb t, \mu)$. Hence, using also that $\Re \br{b}^{j} = \Re b^{j}$, it follows that
\begin{equation} \label{eq:Mz-d-v-L2}
	\nrm{\tld{u}^{\ast}(t)}_{L^{2}} \leq M(\lmb t, \mu)^{-1}.
\end{equation}
{The following lemma quantifies the error $\err[\tld{u}^{\ast}] = (i \rd_{t} + \lap) \tld{u}^{\ast} - \rd_{j} (\br{b}^{j}(x) \tld{u}^{\ast})$ incurred by $\tld{u}^{\ast}$. \begin{lemma}\label{lem:error2}
	There exists a constant $C_{0}$, which depends only on $\nrm{b}_{C^{1, 1}}$ and $\nrm{\psi_{1}}_{H^{2}}$, such that
	\begin{equation} \label{eq:Mz-d-v-err}
		\nrm{\err[\tld{u}^{\ast}](t)}_{L^{2}} \leq C_{0} (\mu + \lmb t)^{2} M(\lmb t, \mu)^{-1}.
	\end{equation}
\end{lemma}}
We are now ready to implement the duality method. Assume for the moment that, for some $B > 0$, we have
\begin{equation} \label{eq:Mz-d-B}
	\frac{1}{t_{f}} \nrm*{\frac{(1 + \mu^{-1} \lmb t)^{2}}{(1 + \mu^{-1} \lmb t_{f})^{2}} M(\lmb t, \mu)^{-1} u}_{L^{1}_{t}([0, t_{f}]; L^{2})} < B \nrm{u_{0}}_{L^{2}}.
\end{equation}
By \eqref{eq:Mz-d-gen-en} and \eqref{eq:Mz-d-v-err}, we then have
\begin{equation*}
	\abs*{\frac{\ud}{\ud t} \brk{u, \tld{u}^{\ast}}} \leq C_{0} (\mu + \lmb t)^{2} M(\lmb t, \mu)^{-1}  \nrm{u(t)}_{L^{2}} .
\end{equation*}
Integrating in $t$ and using the contradiction assumption, we arrive at
\begin{equation*}
	\brk{u, \tld{u}^{\ast}}(t) \geq \nrm{v_{0}}_{L^{2}}\left(1 - B C_{0} \mu^{2} (1+\mu^{-1} \lmb t_{f})^{2} t_{f} \right).
\end{equation*}
Suppose that
\begin{equation}  \label{eq:Mz-d-hyp-B}
	\mu \leq (\tfrac{1}{8 B C_{0}})^{\frac{1}{3}} \lmb^{\frac{1}{3}}, \quad t_{f} \leq (\tfrac{1}{8 B C_{0}})^{\frac{1}{3}} \lmb^{-\frac{2}{3}}.
\end{equation}
Dividing into two cases $t_{f} \leq \frac{\mu}{\lmb}$ and $t_{f} \geq \frac{\mu}{\lmb}$, it follows that 
\begin{equation*}
B C_{0} \mu^{2} (1+\mu^{-1} \lmb t_{f})^{2} t_{f} \leq \frac{1}{2}.
\end{equation*}
 Therefore,
\begin{equation*}
 \brk{u, \tld{u}^{\ast}}(t) \geq \frac{1}{2} \nrm{u_{0}}_{L^{2}} \quad \hbox{ for all } 0 \leq t \leq t_{f}.
\end{equation*}
Then, applying \eqref{eq:Mz-d-v-L2}, we obtain the lower bound 
\begin{equation} \label{eq:Mz-d-ptwise}
	\nrm{u(t)}_{L^{2}} \geq \frac{1}{2} M(\lmb t, \mu) \nrm{u_{0}}_{L^{2}}.
\end{equation}
We now multiply both sides by $(1 + \mu^{-1} \lmb t)^{2} M(\lmb t, \mu)^{-1}$ and integrate. Since
\begin{equation*}
	\int_{0}^{t_{f}} (1 + \mu^{-1} \lmb t)^{2} \, \ud t \geq \frac{1}{3} t_{f} (1+\mu^{-1} \lmb t_{f})^{2},
\end{equation*}
we arrive at
\begin{equation} \label{eq:Mz-d-1/6}
\frac{1}{t_{f}} \nrm*{\frac{(1 + \mu^{-1} \lmb t)^{2}}{(1 + \mu^{-1} \lmb t_{f})^{2}} M(\lmb t, \mu)^{-1}  u}_{L^{1}_{t}([0, t_{f}]; L^{2})} \geq \frac{1}{6} \nrm{u_{0}}_{L^{2}}.
\end{equation}

To complete the proof, we assume, for the purpose of contradiction, that \eqref{eq:Mz-d-B} holds with $B = \frac{1}{6}$.
Take $c = (\tfrac{3}{4 C_{0}})^{\frac{1}{3}}$ in \eqref{eq:Mz-d-hyp} so that \eqref{eq:Mz-d-hyp-B} is satisfied. Then by the preceding argument, we arrive at \eqref{eq:Mz-d-1/6}, which is a contradiction that establishes \eqref{eq:Mz-d}.  \qedhere
\end{proof}
\begin{proof}[Proof of Lemma \ref{lem:error2}]
	As in the proof of Lemma \ref{lem:error}, we compute with $\psi = \psi_{\mu, x_{0}}$ that  \begin{equation*}
		\begin{split}
			i \rd_{t}  \tld{u}^{\ast} = \lmb^{2} \tld{u}^{\ast} 
			+ i \lmb \br{b}^{j}(x-2\lmb t \omg_{0}) (\omg_{0})_{j} \tld{u}^{\ast} 
			- 2i\lmb\frac{ \omg_{0} \cdot\nb \psi}{\psi}  \tld{u}^{\ast}.
		\end{split}
	\end{equation*} Introducing for simplicity \begin{equation*}
	\begin{split}
		I_{k}(x) = - \int_{0}^{\lmb t} \,  \rd_{k} \br{b}^{j} (x-2 s \omg_{0}) (\omg_{0})_{j} \, \ud s , \qquad I(x) = (I_{1}, \cdots, I_{d}),
	\end{split}
\end{equation*} we have \begin{equation*}
		\begin{split}
			\rd_{k}  \tld{u}^{\ast} &= \left( i\lmb (\omg_{0})_{k} - I_{k}  +  \frac{\rd_{k}\psi}{\psi} \right)  \tld{u}^{\ast}, \\
			\lap \tld{u}^{\ast} &= - \lmb^{2} |\omg_{0}|^{2}  \tld{u}^{\ast}  - 2i\lmb  \omg_{0} \cdot I \, \tld{u}^{\ast} + 2i\lmb \frac{\omg_{0} \cdot \nb \psi}{\psi} \, \tld{u}^{\ast} + \calR[  \tld{u}^{\ast}], 
		\end{split}
	\end{equation*} where \begin{equation*}
	\begin{split}
		\calR[  \tld{u}^{\ast}] = \left( |I|^{2} - \frac{2 I\cdot \nb\psi}{\psi} + \frac{|\nb\psi|^{2}}{\psi^{2}} + \sum_{k} \rd_{k} \frac{\rd_{k}\psi}{\psi} - \nb \cdot I \right) \tld{u}^{\ast}. 
	\end{split}
\end{equation*}
	Then, after several direct cancellations, we have  \begin{equation*}
		\begin{split}
			(i \rd_{t} + \tld{\calL}) \tld{u}^{\ast}  - \rd_{j} (\br{b}^{j}(x)  \tld{u}^{\ast}) & = \left(-2i\lmb \omg_{0} \cdot I + i\lmb \br{b}^{j}(x-2\lmb t\omg_{0}) (\omg_{0} )_{j} - i\lmb \br{b}^{j}(x) (\omg_{0})_{j} \right) \tld{u}^{\ast} \\
			&\qquad + \calR[ \tld{u}^{\ast}] - \br{b}^{j}(x) I_{j}(x) \tld{u}^{\ast} + \frac{\br{b}^{j}(x) \rd_{j}\psi}{\psi} \tld{u}^{\ast} - (\rd_{j} \br{b}^{j}(x) )  \tld{u}^{\ast}
		\end{split}
	\end{equation*} and as in the one-dimensional case, we use that \begin{equation*}
	\begin{split}
		\frac12 \left( \br{b}^{j}(x-2\lmb t \omg_{0}) - \br{b}^{j} (x)  \right) = \int_{0}^{\lmb t} \frac12 \frac{\ud}{\ud s} \br{b}^{j}(x-2s\omg_{0}) \, \ud s = \sum_{k} I_{k}(x) (\omg_{0})_{k} 
	\end{split}
\end{equation*} to get cancellations among the $O(\lmb)$ terms\footnote{Unlike the one-dimensional case, however, we cannot eliminate the integral on the domain $[0, \lmb t]$ in $I$. Hence, we let the right-hand side of \eqref{eq:Mz-d-v-err} depend on $\lmb t$.}. It is now not difficult to see that the remaining terms are bounded in $L^{2}$ by the right hand side of \eqref{eq:Mz-d-v-err}. 
\end{proof}

As alluded to before, the second result we shall prove using essentially the same argument is Proposition~\ref{prop:Mz-d-norm-infl}, i.e., that the failure of the Takeuchi--Mizohata condition (see \eqref{eq:Mz-d-cond}) implies \emph{norm inflation} for \eqref{eq:Mz-d}.

\begin{proof}[Proof of Proposition~\ref{prop:Mz-d-norm-infl}]
Assume, for contradiction, that there exists $B_{0} < + \infty$ such that, for every $u_{0} \in L^{2}$,
\begin{equation} \label{eq:Mz-d-norm-bdd}
	\nrm{u}_{L^{\infty}([0, \dlt]; L^{2})} \leq B_{0} \nrm{u_{0}}_{L^{2}}.
\end{equation}

By \eqref{eq:Mz-d-cond}, there exists a sequence $(x_{n}, \omg_{n}, T_{n})$ such that
\begin{equation*}
	M_{x_{n}, \omg_{n}}(T_{n}) := \exp \left( \int_{0}^{T_{n}} \Re b^{j}(x_{n} - 2 s \omg_{n}) (\omg_{n})_{j} \, \ud x \right) \geq e^{2 n}.
\end{equation*}
By restarting from the point $x_{n} + 2 T \omg_{n}$ where $M_{x_{n}, \omg_{n}}(T) = 1$ if necessary, we may assume also that $M_{x_{n}, \omg_{n}}(T) \geq 1$ for all $0 \leq T \leq T_{n}$. Since $M_{x_{n}, \omg_{n}}(T_{n}, \mu) \to M_{x_{n}, \omg_{n}}(T_{n})$ as $\mu \to \infty$, we may choose $\mu_{n}$ so that
\begin{equation*}
	M_{x_{n}, \omg_{n}}(T_{n}, \mu_{n}) \geq e^{n}.
\end{equation*}

We shall apply the argument in the proof of Proposition~\ref{prop:Mz-d} with the parameters 
\begin{equation*}
t_{f} = \frac{T_{n}}{\lmb_{n}}, \quad x_{0} = x_{n}, \quad \omg_{0} = \omg_{n}, \quad \mu = \mu_{n}, \quad \lmb = \lmb_{n}, 
\end{equation*}
where $\lmb_{n}$ shall be determined below. We denote by $\tld{u}^{\ast}_{n}$ the wave packet for $i \rd_{t} + \tld{\calL}^{\ast}$ with these parameters, and $u_{n}$ the solution to \eqref{eq:Mz-d} with initial data $u_{0} = \tld{u}^{\ast}_{n}(0)$ satisfying \eqref{eq:Mz-d-norm-bdd}. Taking $\lmb_{n}$ to be large enough, we may guarantee that $t_{f} = \frac{T_{n}}{\lmb_{n}} \leq \dlt$. Then by the contradiction assumption and the bound $M(T) \geq 1$, it follows that \eqref{eq:Mz-d-B} is satisfied for $u = u_{n}$ with $B = C_{1} B_{0}$, where $C_{1}$ depends only on $\nrm{\rd b}_{L^{\infty}}$. Furthermore, choosing $\lmb_{n}$ sufficiently large depending on $B_{0}$, $C_{0}$, $C_{1}$ and $T_{n}$, we may ensure that \eqref{eq:Mz-d-hyp-B} holds (here, it is important that the power of $\lmb$ in the second inequality is greater than $-1$). Thence, it follows from \eqref{eq:Mz-d-ptwise} and our choices of parameters that
\begin{equation*}
	\nrm{u_{n}(\tfrac{T_{n}}{\lmb_{n}})}_{L^{2}} \geq \frac{1}{2} e^{n} \nrm{u_{n}(0)}_{L^{2}}.
\end{equation*}
Taking $n \to \infty$, we arrive at a contradiction. \qedhere
\end{proof}
\begin{remark}
We note that when $d = 1$, Proposition~\ref{prop:Mz-d-norm-infl} is essentially a consequence of Proposition~\ref{prop:Mz-1}, although pedantically the notion of solution is slightly different due to the presence of a conjugation when $d = 1$. On the other hand, the preceding proof applies to all $d \geq 1$.
\end{remark}

\bibliographystyle{abbrv}
\bibliography{degen_disp}

\end{document}